\documentclass[oneside,a4paper,12pt]{article}

\usepackage[utf8]{inputenc} 
\usepackage[T1]{fontenc}
\usepackage{color, xcolor} 

\usepackage{abstract} 
\usepackage[hidelinks]{hyperref} 
\usepackage{authblk} 
\usepackage{titlesec} 
\usepackage{cite} 
\usepackage{appendix} 
\usepackage[titles]{tocloft}

\usepackage{multicol} 
\usepackage{fancyhdr} 
\usepackage{listings} 
\usepackage{graphicx, subfig} 
\usepackage{float} 
\usepackage{enumerate} 
\usepackage{longtable} 
\usepackage{multirow} 
\usepackage{makecell} 
\usepackage{tabularx} 
\usepackage{array} 

\usepackage{amsfonts} 
\usepackage{amsmath} 
\usepackage{amssymb} 
\usepackage{gensymb} 
\usepackage{amsthm} 
\usepackage[ruled]{algorithm2e} 
\usepackage{tikz} 
\usepackage{tikz-cd}
\usetikzlibrary{decorations.markings}
\usetikzlibrary{arrows.meta}
\usetikzlibrary{positioning}
\usetikzlibrary{shapes.geometric}
\usetikzlibrary{fit}
\usetikzlibrary{backgrounds}

\usepackage{dsfont} 
\usepackage{ulem} 
\usepackage{float} 

\usepackage[english]{babel}
\newtheorem{remark}{Remark}
\newtheorem{definition}{Definition}

\usepackage[total={7in,9in}]{geometry}
\usepackage{lineno}

\colorlet{ZW_COLOR}{blue}
\newcommand\zw[1]{#1}

\begin{document}

\title{An $h$-adaptive Tetrahedral Spectral Element Method with Applications to Kohn-Sham Density Functional Theory}
\author[1]{Zeyu Wang}
\author[2]{Hongfei Zhan}
\author[3]{Guanghui Hu}
\affil[1]{Department of Mathematics, Faculty of Science, University of Macau, Macao SAR, China}
\affil[2]{Department of Mathematics, National University of Singapore, Singapore 119076}
\affil[3]{Laboratory of Internet of Things for Smart City and Department of Mathematics, University of Macau, Macao, China}
\date{}

\maketitle

\begin{abstract}

    High-order $h$-adaptive spectral element methods on tetrahedral meshes provide an effective framework for resolving localized singularities and multiscale structures in complex three-dimensional geometries. However, their development is often hindered by difficulties in maintaining $C^0$ continuity across refinement interfaces and efficiently transferring solutions between adaptive meshes. Such limitations are particularly relevant in demanding applications such as all-electron Kohn-Sham density functional theory, which place stringent requirements on the accurate resolution of both nuclear singularities and multiple physical scales. In this paper, we present an efficient $h$-adaptive tetrahedral spectral element framework. To address the continuity challenge, we develop an adaptive strategy that combines element orientation alignment with geometric red-green refinement, thereby eliminating the need for algebraic hanging-node constraints while preserving inter-element continuity. Furthermore, an efficient topology-based point-location algorithm is introduced to accelerate interpolation between adaptive meshes. Numerical experiments on Poisson and Laplacian eigenvalue problems confirm the spectral convergence of the proposed method. Applications to all-electron Kohn-Sham equations further demonstrate its capability to accurately resolve nuclear singularities. Moreover, parallel performance studies exhibit excellent scalability, with matrix assembly and adaptivity modules generally achieving speedups above 15 times and the proposed interpolation algorithm attaining speedups ranging from 25 to 35 on 64-core configurations compared to the single-core performance. These results indicate that the proposed framework provides an accurate, robust, and efficient solution for large-scale, high-resolution simulations.
    
    \textbf{Keywords:} $h$-adaptive method, red-green refinement, spectral element method, Kohn-Sham density functional theory
\end{abstract}





\section{Introduction}

High-precision and computationally efficient solvers remain central in scientific computing and engineering applications, encompassing fluid dynamics \cite{hu1999weighted, jiang1996efficient}, ground-state electronic structure simulations \cite{umemoto2006dissociation, umemoto2006namgf3}
, and high-order harmonic generation \cite{li1989multiple, bao2015real}.
In these domains, it is common for solutions to exhibit multi-scale characteristics or singularities, which pose significant challenges for traditional low-order methods \cite{visbal2002use, wang2007high}.

Given the demand for high-resolution simulations, high-order methods are often preferred for their ability to achieve high accuracy with fewer degrees of freedom, including the spectral method \cite{gottlieb1977numerical, gottlieb2001spectral, canuto2006spectral, shen2011spectral}, the finite element method \cite{zienkiewicz1977finite, hughes2003finite, brenner2008mathematical}, the discontinuous finite element method \cite{riviere2008discontinuous, hesthaven2008nodal}, the finite difference method \cite{smith1985numerical}, and the wavelet method \cite{cohen2000wavelet}.
Among these methods, the spectral element method (SEM) \cite{patera1984spectral, dubiner1991spectral, hesthaven2000stable, pasquetti2004spectral, huang2019spectral, weiss2023spectral, jia2022sparse, zhan2023novel} provides a robust compromise by combining spectral convergence with the local flexibility of finite elements.
In particular, compared to hexahedral spectral element methods \cite{huang2019spectral, weiss2023spectral}, tetrahedral spectral elements (TSEM) \cite{hesthaven2000stable, jia2022sparse, zhan2023novel} offer notable flexibility toward complex geometries and adaptivity \cite{li2010spectral, schneider2022large}.
Moreover, the presence of multi-scale and singular features complicates mesh generation for high-order methods like TSEM, necessitating adaptivity to ensure computational robustness.


Adaptive methods \cite{
melenk2001residual, leitner1995three, dorfler1996convergent, bonito2024adaptive} address the inefficiency of discretizations by concentrating computational effort in regions where the solution varies rapidly or exhibits local singularities. Three principal types of adaptivity are commonly employed: $h$-adaptivity \cite{zienkiewicz1991adaptivity, hu1998h, benito2003h}, which locally subdivides elements; $p$-adaptivity \cite{demkowicz1989toward, salagame1997simple, mitchell2011survey}, which varies the polynomial degree of basis functions; and $r$-adaptivity \cite{budd2009adaptivity, browne2014fast}, which redistributes the mesh. The adaptive SEM mainly focuses on $p$-adaptive \cite{moxey2017towards} methods, $h$-adaptive methods \cite{mavriplis1994adaptive, hsu1997adaptive, henderson1999adaptive}, and their combination \cite{valenciano2000h,galvao2008hp}. Specifically, $h$-adaptive TSEM provides the requisite resolution for complex topologies and fine-scale structures,
enabling the accurate representation of complex phenomena in computational science and engineering.


As a result, adaptive SEMs have been extensively utilized in fields ranging from fluid dynamics \cite{mavriplis1994adaptive, henderson1999adaptive, moxey2017towards} and electromagnetics \cite{mahariq2017application, kopriva2002computation} to electronic structure theory \cite{motamarri2013higher, kanungo2019real, zhan2023novel, UNKNOWN-CITE}. Within this context, all-electron first-principles calculations, specifically Kohn-Sham Density Functional Theory (KSDFT) \cite{hohenberg1964inhomogeneous, kohn1965self, parr1995density}, present formidable numerical challenges. These simulations demand the simultaneous resolution of Coulomb singularities near atomic nuclei, highly accurate evaluation of long-range Hartree potentials, and the management of multiscale coupling in large systems. Since the solutions are typically smooth in bulk regions but exhibit cusp-like behavior near the nuclei, adaptive high-order discretizations are well suited to capture such localized features \cite{motamarri2013higher, kanungo2019real, zhan2023novel, UNKNOWN-CITE}. 
The growing demand for related applications has driven the development of computational software and numerical libraries.

Beyond the theoretical formulation and numerical algorithms, the robust software implementation of adaptive high-order methods remains a non-trivial challenge. Current state-of-the-art libraries exhibit distinct architectural priorities: deal.II \cite{2025:arndt.bangerth.ea:deal} is primarily optimized for quadrilateral and hexahedral elements, lacking the specialized infrastructure required for fully adaptive tetrahedral meshes. Although MFEM \cite{mfem-2024} offers significant flexibility across various geometries and high-order spaces, it typically manages non-conformal interfaces via algebraic constraints on hanging nodes. This adversely affects the condition number of the global system and complicates the implementation of efficient solvers. Similarly, despite its strengths in spectral element formulations for complex flows, Nektar++ \cite{CANTWELL2015205, MOXEY2020107110} often encounters significant implementation complexity when integrating $h$-adaptivity with tetrahedral elements. 
Hence, owing to challenges such as geometric data management and parallel implementation, the development of robust computational frameworks for $h$-adaptive high-order methods remains an active area of research.

An issue that merits attention in such frameworks is the efficient and accurate mapping of solutions across disparate mesh levels. Several strategies are typically employed for solution interpolation between coarse and fine meshes. For nested grids with nodal basis functions, prolongation matrices are standard \cite{zhang2011high, cerveny2019nonconforming}. However, this approach is often inapplicable to modal basis functions that depend on element orientation. For more general, non-nested cases, one common approach involves evaluating values at quadrature points for interpolation \cite{ortiz1991adaptive, bussetta2011comparison}, though this incurs significant computational overhead due to cell-searching algorithms. Alternatively, the solution can be represented as a point cloud \cite{costin2013numerical, garon2020mesh} to reconstruct the field on the new mesh, albeit at the potential cost of numerical accuracy.



In this work, we design an efficient $h$-adaptive tetrahedral spectral element framework utilizing the generalized Koornwinder polynomials. To address the $C^0$ continuity challenges inherent in high-order modal bases, a constructive approach is proposed to ensure a consistent orientation between adjacent elements. Building upon this framework, a geometric red-green refinement strategy is adopted to maintain mesh conformity. By resolving non-conformal interfaces geometrically rather than through algebraic hanging-node constraints, our approach decouples the global assembly from the refinement logic, thereby preserving the native sparsity and symmetry of the system matrix. Subsequently, a fast point-location algorithm leveraging prior local topological information is integrated into the interpolation procedure to achieve efficient solution transfer between meshes. This integration allows the framework to inherit the high-order precision of the TSEM while gaining the spatial flexibility required for multi-scale problems.

We first verify the proposed method using the Poisson problem and the Laplacian eigenvalue problem, where the expected spectral accuracy is successfully observed. Subsequently, our framework is applied to the all-electron Kohn-Sham equations, where the nuclear singularities are well resolved. It is worth mentioning that a specialized high-order interpolation operator provides a ten-fold acceleration in convergence compared to standard initializations. 
Moreover, parallel performance tests on 64 cores illustrate that speedups for assembly and adaptive steps consistently exceed 10-fold, frequently reaching a more typical range of 15 times to 20 times. Meanwhile, a speedup greater than 25-fold is observed for interpolation. These speedups are already close to the limit of Amdahl's law imposed by the serial portions, highlighting the potential of our framework for large-scale, high-resolution physical simulations.

The paper is structured as follows. In Section \ref{chapter::Methodology}, the fundamentals of the spectral element method and $h$-adaptive strategies are reviewed. Section \ref{chapter::Implementation} details the implementation and presents numerical results for the Poisson equation and Laplacian eigenvalue problem. The framework is further applied to KSDFT in Section \ref{chapter::KohnShamEquations}. Section \ref{chapter::Conclusion} concludes the paper.

\section{\zw{An Adaptive Tetrahedral Spectral Element Method}}
\label{chapter::Methodology}

In this section, we employ the classical Poisson problem as a benchmark to introduce the spectral element basis functions, the discretization procedure, and the $h$-adaptive framework. We begin by introducing the modal basis functions defined on tetrahedral elements, as well as the standard Galerkin framework to formulate the TSEM discretization. Subsequently, in the context of adaptive refinement, the red-green refinement algorithm and a residual-based error indicator are discussed.



The Poisson problem serves as the first benchmark to validate the accuracy and efficiency of the proposed framework. It is a fundamental equation widely used in fluid mechanics \cite{mohammad2025advanced}, heat conduction \cite{frkackowiak2011solution}, and electronic structure calculations \cite{garcia2014survey}. In particular, let $\Omega = (-1, 1)^{3}$ be the physical domain and $H^{1}(\Omega)$ be the associated Sobolev space. The Poisson equation in $\Omega$ subject to Dirichlet boundary conditions is given by
\begin{equation}
   \begin{cases}
       \displaystyle - \Delta u(\mathbf{x}) = f(\mathbf{x}), & \mathbf{x} \in \Omega,          \\
       \displaystyle u(\mathbf{x}) = u_{b} (\mathbf{x}),     & \mathbf{x} \in \partial \Omega,
   \end{cases}
\end{equation}
where $u_{b} (\mathbf{x})$ is a prescribed function on the boundary. The corresponding weak formulation seeks $u \in H^1(\Omega)$ with $u|_{\partial \Omega} = u_b$ such that
\begin{equation}
   \int_{\Omega} \nabla u \cdot \nabla v \ \mathrm{d} x = \int_{\Omega} f \cdot v \ \mathrm{d} x,\ \forall v \in H_{0}^{1} (\Omega).
\end{equation}

\subsection{\zw{Tetrahedral Spectral Element Discretization}}
\label{sec::Methodology::TetrahedralSpectralElementMethod}

Following previous studies \cite{jia2022sparse, zhan2023novel}, we utilize basis functions constructed from generalized Jacobi polynomials \cite{guo2006optimal, guo2009generalized} and generalized Koornwinder polynomials \cite{li2010optimal, shan2017triangular}. Specifically, let $\boldsymbol{\alpha} = (\alpha_{0}, \alpha_{1}, \alpha_{2}, \alpha_{3}) \in [-1, +\infty)^{4}$ be a parameter vector, and let $(i, j, k) \in \mathbb{N}^{3}$ denote the multi-index. Defined on the reference tetrahedron

\begin{equation}
    \mathcal{T} := \{ (x, y, z) \in [0, 1]^{3}: x + y + z \leq 1 \},
\end{equation}

\noindent the generalized Koornwinder polynomials are formulated as

\begin{equation}
    \begin{aligned}
        \mathcal{J}_{i, j, k}^{\boldsymbol{\alpha}}(x, y, z) = & (1 - y - z)^{i} J^{\alpha_{0}, \alpha_{1}}_{i}\left(\frac{2 x}{1 - y - z} - 1\right) (1 - z)^{j}                                                                                  \\
                                                                  & \times J^{2 i + \alpha_{0} + \alpha_{1} + 1, \alpha_{2}}_{j}\left(\frac{2 y}{1 - z} - 1\right) J^{2 i + 2 j + \alpha_{0} + \alpha_{1} + \alpha_{2} + 2, \alpha_{3}}_{k}(2 z - 1),
    \end{aligned}
\end{equation}

\noindent where $J^{\alpha, \beta}_{n}(x)$ represents the $n$-th order generalized Jacobi polynomial.

\begin{figure}[H]
    \centering
    \begin{tikzpicture}[scale=4, line join=round, line cap=round]
        \coordinate (V0) at (0.0, 0.0);
        \coordinate (V1) at (1.0, 0.0);
        \coordinate (V2) at (0.9, 0.6);
        \coordinate (V3) at (0.0, 1.0);
            
        \draw[red, ultra thick] (V0) -- (V1);
        \draw[ultra thick, dashed] (V0) -- (V2);
        \draw[orange, ultra thick] (V0) -- (V3);
        \draw[ultra thick] (V1) -- (V2);
        \draw[blue, ultra thick] (V1) -- (V3);
        \draw[ultra thick] (V2) -- (V3);

        \fill[cyan, opacity=0.1] (V0) -- (V1) -- (V3) -- cycle;
    
        \filldraw (V0) circle (0.5pt) node[anchor=north east] {$\mathcal{V}_{0} (0,0,0)$};
        \filldraw (V1) circle (0.5pt) node[anchor=north west] {$\mathcal{V}_{1} (1,0,0)$};
        \filldraw (V2) circle (0.5pt) node[anchor=south west] {$\mathcal{V}_{2} (0,1,0)$};
        \filldraw (V3) circle (0.5pt) node[anchor=south]      {$\mathcal{V}_{3} (0,0,1)$};
    
        \node[cyan] at (0.3, 0.3) {$\mathcal{F}_2$};

        \path[red] (V0) -- (V1) node[midway, below] {$\mathcal{E}_{01}$};
        \path[orange] (V0) -- (V3) node[midway, left] {$\mathcal{E}_{03}$};
        \path[blue] (V1) -- (V3) node[midway, left] {$\mathcal{E}_{13}$};
    \end{tikzpicture}
    \caption{The entities of a reference tetrahedron.}
    \label{figure::reference-tetrahedron}
\end{figure}
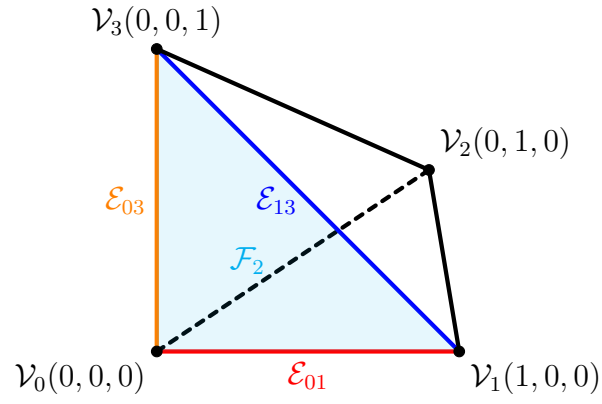

Setting $\boldsymbol{\alpha} = -\mathds{1}$, we obtain the modal basis functions summarized in Table \ref{table::BasisPolynomialFunction}. The hierarchical structure of these basis functions corresponds to specific geometric entities of the tetrahedron, as illustrated in Figure \ref{figure::reference-tetrahedron}:

\begin{itemize}
    \item \textbf{Vertex modes}: Coincide with standard linear finite element shape functions, each associated with a specific vertex;
    \item \textbf{Edge modes}: Defined on edges that connect two vertices, such as $\mathcal{E}_{12}$ associated with the edge between vertices $\mathcal{V}_{1}$ and $\mathcal{V}_{2}$, with polynomial degrees $p \geq 2$;
    \item \textbf{Face modes}: Defined on the four triangular faces. We denote $\mathcal{F}_{i}$ as the face opposite to vertex $\mathcal{V}_{i}$, with degrees $p \geq 3$;
    \item \textbf{Interior modes}: Associated with the tetrahedral volume itself, with degrees $p \geq 4$.
\end{itemize}

This construction ensures that each class of basis functions vanishes on lower-dimensional geometric entities. Specifically, the interior modes vanish on the entire tetrahedral boundary, face modes vanish on all edges and vertices, while edge modes take the value 0 on all vertices. Consequently, the numerical solution on any given face is uniquely determined by the modes associated with that face and its constituent edges and vertices. 
Therefore, the establishment of $C^0$ continuity critically depends on the correspondence between degrees of freedom (DoFs) from adjacent elements.

\begin{table}[H]
    \centering
    \caption{Polynomial basis functions for tetrahedron spectral element method.}
    \label{table::BasisPolynomialFunction}
    \renewcommand{\arraystretch}{2}
    \begin{tabular}{|l|lr|}
        \hline
        \multicolumn{3}{|c|}{Vertex modes}                                                                                                                                                                                                                                                                                                                          \\
        \hline
        $\mathcal{V}_{0}$  & $\displaystyle \phi_{0, 0, 0}(\mathbf{x}) =  \frac{1}{8} \mathcal{J}_{0, 0, 0}^{-\mathds{1}}(\mathbf{x}) - \frac{1}{2} \mathcal{J}_{1, 0, 0}^{-\mathds{1}}(\mathbf{x}) - \frac{1}{4} \mathcal{J}_{0, 1, 0}^{-\mathds{1}}(\mathbf{x}) - \frac{1}{8} \mathcal{J}_{0, 0, 1}^{-\mathds{1}}(\mathbf{x})$ &                                  \\
        $\mathcal{V}_{1}$  & $\displaystyle \phi_{1, 0, 0}(\mathbf{x}) =  \frac{1}{8} \mathcal{J}_{0, 0, 0}^{-\mathds{1}}(\mathbf{x}) + \frac{1}{2} \mathcal{J}_{1, 0, 0}^{-\mathds{1}}(\mathbf{x}) - \frac{1}{4} \mathcal{J}_{0, 1, 0}^{-\mathds{1}}(\mathbf{x}) - \frac{1}{8} \mathcal{J}_{0, 0, 1}^{-\mathds{1}}(\mathbf{x})$ &                                  \\
        $\mathcal{V}_{2}$  & $\displaystyle \phi_{0, 1, 0}(\mathbf{x}) =  \frac{1}{4} \mathcal{J}_{0, 0, 0}^{-\mathds{1}}(\mathbf{x}) + \frac{1}{2} \mathcal{J}_{0, 1, 0}^{-\mathds{1}}(\mathbf{x}) - \frac{1}{4} \mathcal{J}_{0, 0, 1}^{-\mathds{1}}(\mathbf{x})$                                                               &                                  \\
        $\mathcal{V}_{3}$  & $\displaystyle \phi_{0, 0, 1}(\mathbf{x}) =  \frac{1}{2} \mathcal{J}_{0, 0, 0}^{-\mathds{1}}(\mathbf{x}) + \frac{1}{2} \mathcal{J}_{0, 0, 1}^{-\mathds{1}}(\mathbf{x})$                                                                                                                             &                                  \\
        \hline
        \multicolumn{3}{|c|}{Edge modes}                                                                                                                                                                                                                                                                                                                            \\
        \hline
        $\mathcal{E}_{01}$ & $\displaystyle \phi_{i, 0, 0}(\mathbf{x}) = 2 \mathcal{J}_{i, 0, 0}^{-\mathds{1}}(\mathbf{x})$,                                                                                                                                                                                                     & $(i \geq 2)$                     \\
        $\mathcal{E}_{02}$ & $\displaystyle \phi_{0, j, 0}(\mathbf{x}) = \mathcal{J}_{0, j, 0}^{-\mathds{1}}(\mathbf{x}) + \frac{j - 1}{j} \mathcal{J}_{1, j - 1, 0}^{-\mathds{1}}(\mathbf{x})$,                                                                                                                                 & $(j \geq 2)$                     \\
        $\mathcal{E}_{03}$ & $\displaystyle \phi_{0, 0, k}(\mathbf{x}) = \frac{1}{2} \mathcal{J}_{0, 0, k}^{-\mathds{1}}(\mathbf{x}) + \frac{k - 1}{2 k} \mathcal{J}_{0, 1, k - 1}^{-\mathds{1}}(\mathbf{x}) + \frac{k - 1}{k} \mathcal{J}_{1, 0, k - 1}^{-\mathds{1}}(\mathbf{x})$,                                             & $(k \geq 2)$                     \\
        $\mathcal{E}_{13}$ & $\displaystyle \phi_{1, 0, k - 1}(\mathbf{x}) = \frac{1}{2} \mathcal{J}_{0, 0, k}^{-\mathds{1}}(\mathbf{x}) + \frac{k - 1}{2 k} \mathcal{J}_{0, 1, k - 1}^{-\mathds{1}}(\mathbf{x}) - \frac{k - 1}{k} \mathcal{J}_{1, 0, k - 1}^{-\mathds{1}}(\mathbf{x})$,                                         & $(k \geq 2)$                     \\
        $\mathcal{E}_{12}$ & $\displaystyle \phi_{1, j - 1, 0}(\mathbf{x}) = \mathcal{J}_{0, j, 0}^{-\mathds{1}}(\mathbf{x}) - \frac{j - 1}{j} \mathcal{J}_{1, j - 1, 0}^{-\mathds{1}}(\mathbf{x})$,                                                                                                                             & $(j \geq 2)$                     \\
        $\mathcal{E}_{23}$ & $\displaystyle \phi_{0, 1, k - 1}(\mathbf{x}) = \mathcal{J}_{0, 0, k}^{-\mathds{1}}(\mathbf{x}) - \frac{k - 1}{k} \mathcal{J}_{0, 1, k - 1}^{-\mathds{1}}(\mathbf{x})$,                                                                                                                             & $(k \geq 2)$                     \\
        \hline
        \multicolumn{3}{|c|}{Face modes}                                                                                                                                                                                                                                                                                                                            \\
        \hline
        $\mathcal{F}_{0}$  & $\displaystyle \phi_{1, j - 1, k}(\mathbf{x}) = \mathcal{J}_{0, j, k}^{-\mathds{1}}(\mathbf{x}) - \frac{j - 1}{j} \mathcal{J}_{1, j - 1, k}^{-\mathds{1}}(\mathbf{x})$,                                                                                                                             & $(j \geq 2, k \geq 1)$           \\
        $\mathcal{F}_{1}$  & $\displaystyle \phi_{0, j, k}(\mathbf{x}) = \mathcal{J}_{0, j, k}^{-\mathds{1}}(\mathbf{x}) + \frac{j - 1}{j} \mathcal{J}_{1, j - 1, k}^{-\mathds{1}}(\mathbf{x})$,                                                                                                                                 & $(j \geq 2, k \geq 1)$           \\
        $\mathcal{F}_{2}$  & $\displaystyle \phi_{i, 0, k}(\mathbf{x}) = 2 \mathcal{J}_{i, 0, k}^{-\mathds{1}}(\mathbf{x})$,                                                                                                                                                                                                     & $(i \geq 2, k \geq 1)$           \\
        $\mathcal{F}_{3}$  & $\displaystyle \phi_{i, j, 0}(\mathbf{x}) = 2 \mathcal{J}_{i, j, 0}^{-\mathds{1}}(\mathbf{x})$,                                                                                                                                                                                                     & $(i \geq 2, j \geq 1)$           \\
        \hline
        \multicolumn{3}{|c|}{Interior modes}                                                                                                                                                                                                                                                                                                                        \\
        \hline
        $\mathcal{I}$      & $\displaystyle \phi_{i, j, k}(\mathbf{x}) = \mathcal{J}_{i, j, k}^{-\mathds{1}}(\mathbf{x})$,                                                                                                                                                                                                       & $(i \geq 2, j \geq 1, k \geq 1)$ \\
        \hline
    \end{tabular}
\end{table}

To provide a comprehensive demonstration of the TSEM discretization, we first give the discretization for each tetrahedron. Let $\{ \mathcal{T}_{i} \}_{i=1}^{N_{\text{ele}}}$ denote tetrahedra, and let $\left\{ \phi^{[i]}_{j} \right\}_{j=1}^{N_{\text{dof}}}$ denote the basis functions on the tetrahedron $\mathcal{T}_{i}$, ordered as vertex, edge, face, and interior modes. Substituting the discrete solution
\begin{equation}
    u_{h}(\mathbf{x}) = \sum_{i=1}^{N_{\text{ele}}} \sum_{j=1}^{N_{\text{dof}}} u^{[i]}_{j} \phi^{[i]}_{j}(\mathbf{x}),
\end{equation}
whose coefficients $u^{[i]}_{j}$ satisfy the constraint that $u_{h}(\mathbf{x}) \in C^{0}(\Omega)$, into the weak formulation, we obtain the following linear system for each tetrahedron:

\begin{equation}
    \sum_{j=1}^{N_{\text{dof}}} u^{[i]}_{j} \int_{\mathcal{T}_{i}} \nabla \phi^{[i]}_{j}(\mathbf{x}) \cdot \nabla \phi^{[i]}_{k} \ \mathrm{d} x = \int_{\mathcal{T}_{i}} f \cdot \phi^{[i]}_{k} \ \mathrm{d} x.
    \label{eqn:weak-form}
\end{equation}

Next, we establish the assembly procedure by imposing the necessary constraints to ensure the $C^0$ continuity of $u_h$ across the global domain. Consider an interface $\mathcal{F}_{ij}$ shared by two adjacent tetrahedra $\mathcal{T}_{i}$ and $\mathcal{T}_{j}$. Let $\mathbf{u}^{[i]}_{\mathcal{F}_{ij}}$ be the DoFs corresponding to the face $\mathcal{F}_{ij}$ in the $i$-th tetrahedron.
There exists a linear transform $T_{ij}$ such that
\begin{equation}
        \mathbf{u}^{[i]}_{\mathcal{F}_{ij}} \\
    = T_{ij} 
        \mathbf{u}^{[j]}_{\mathcal{F}_{ij}} \\
\end{equation}
This transformation is uniquely determined by the requirement that $u_h \in C^0(\Omega)$. A similar mapping can be applied to the edge modes. By incorporating these interface and edge constraints into the local weak form (Equation (\ref{eqn:weak-form})), we arrive at the global linear system for the TSEM.

\begin{remark}
    The matrix $T_{ij}$ has a strong impact on the sparsity pattern of the linear system, which directly influences computational efficiency. Therefore, it requires careful treatment. Moreover, its structure depends on element orientations, which will be discussed in Section \ref{sec::Orientation}.
\end{remark}

\subsection{\zw{$h$-adaptive Strategy}}
\label{sec::Methodology::ErrorIndicatorAdaptiveStrategy}

While spectral element methods offer exponential convergence for smooth solutions, this advantage is often compromised by the localized singularities or sharp gradients typical of practical Poisson problems. Consequently, an $h$-adaptive strategy is indispensable to maintain its intrinsic high accuracy.

The adaptive mesh refinement (AMR) procedure in this work is based on the red-green refinement scheme \cite{leitner1995three}, driven by a residual-type error indicator \cite{melenk2001residual}. The process follows the standard adaptive loop:

$$
    \text{Solve} \quad\to\quad \text{Estimate} \quad\to\quad \text{Mark} \quad\to\quad \text{Refine}.
$$

In the estimation step, we define the residual-based error indicator $\eta_{\mathcal{T}}$ for a tetrahedron $\mathcal{T}$ as:
\begin{equation}
    \eta_{\mathcal{T}} = \sqrt{\left(\frac{h_{\mathcal{T}}}{p}\right)^{2} \int_{\mathcal{T}} \vert \mathcal{L}u \vert^{2} \mathrm{d} \mathbf{x} + \frac{1}{2} \sum_{\mathcal{F} \in \partial \mathcal{T}} \frac{h_{\mathcal{F}}}{p} \int_{\mathcal{F}} \left\vert \frac{\partial u_{\mathcal{T}}}{\partial \mathbf{n}_{\mathcal{F}}} - \frac{\partial u_{\mathcal{F}}}{\partial \mathbf{n}_{\mathcal{F}}} \right\vert^{2} \mathrm{d} \mathbf{x}},
\end{equation}
where $\mathcal{L}u=-\Delta u-f$ is the interior residual of the Poisson equation, $u_{\mathcal{T}}$ denotes the numerical solution on $\mathcal{T}$, and $u_{\mathcal{F}}$ represents the solution on the adjacent tetrahedron to face $\mathcal{F}$. 
The parameters $h_{\mathcal{T}}$ and $h_{\mathcal{F}}$ are the diameters of the corresponding tetrahedron and triangle, respectively. Following the maximum marking strategy, an element $\mathcal{T}$ is marked for refinement if $\eta_{\mathcal{T}} > \theta \max_{\tilde{\mathcal{T}}} \eta_{\tilde{\mathcal{T}}}$, with $\theta = 0.5$ in this study.

To preserve mesh conformity, the refinement is performed using a red-green refinement that comprises two rules:
\begin{itemize}
    \item Red Refinement: A tetrahedron is subdivided into eight smaller, similar tetrahedra by bisecting all six edges;
    \item Green Refinement: Transitional templates are applied to non-conforming elements to eliminate hanging nodes locally, ensuring the refinement does not propagate unnecessarily throughout the domain.    
\end{itemize}

The implementation utilizes two priority queues: a red queue for elements marked by the indicator and a green queue for those requiring transitional refinement. The algorithm iterates until both queues are empty. Notably, if a ``green'' element requires further subdivision that violates its transitional template, it is converted to a ``red'' refinement by reverting its temporary subdivision and re-entering the red queue. Detailed refinement templates and algorithmic specifics are provided in Appendix \ref{appendix:red-green-refinement}.



\begin{remark}
    The practical implementation of the red-green refinement necessitates robust handling of complex mesh topologies and efficient data transfer between nested grids. Specifically, maintaining full adjacency information is crucial for the localized green refinement templates. Furthermore, the transfer of solutions and coefficients during the refinement-coarsening cycle requires high-order interpolation operators to preserve accuracy. We provide a comprehensive discussion on the management of tetrahedral topology data in Section \ref{sec::Tetrahedral-Mesh-with-Full-Topology-Data} and detail the interpolation schemes in Section \ref{sec::Interpolation}.
\end{remark}

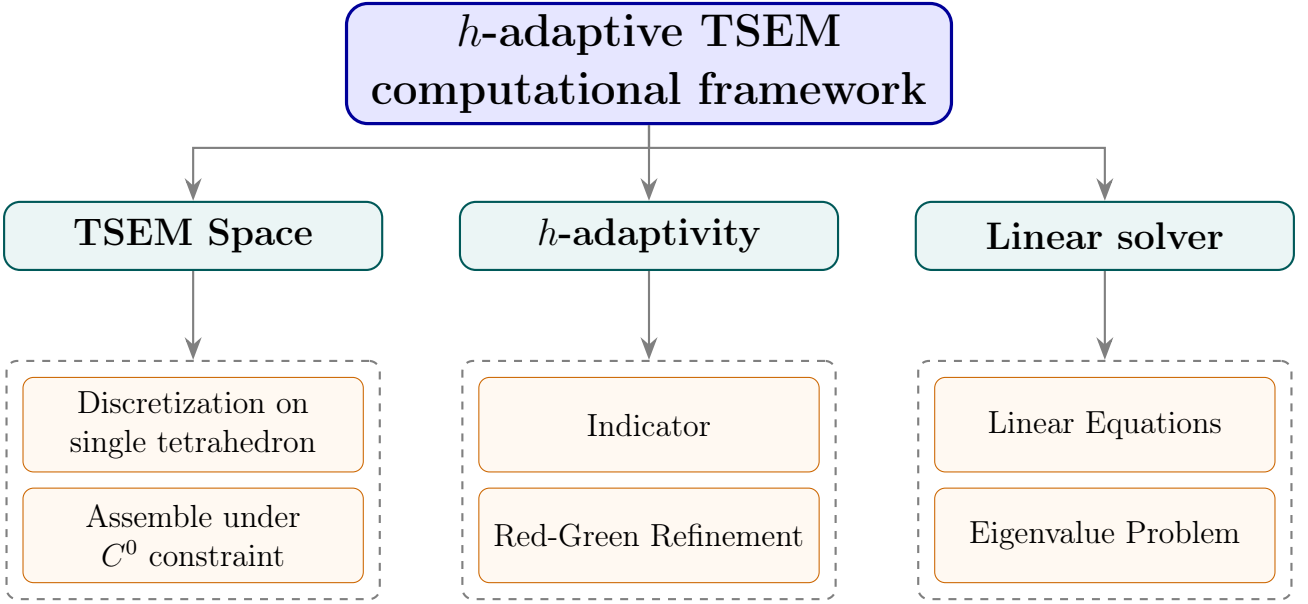
\begin{figure}[H]
    \centering
    \begin{tikzpicture}[
        root/.style = {
            rectangle, rounded corners=8pt, draw=blue!60!black, fill=blue!10, 
            very thick, minimum width=8cm, minimum height=1.2cm, align=center, 
            font=\Large\bfseries
        },
        module/.style = {
            rectangle, rounded corners=6pt, draw=teal!70!black, fill=teal!8, 
            thick, minimum width=5cm, minimum height=0.9cm, align=center, 
            font=\large\bfseries
        },
        sub/.style = {
            rectangle, rounded corners=4pt, draw=orange!80!black, fill=orange!5, 
            thin, minimum width=4.5cm, minimum height=1.25cm, align=center, 
            font=\normalsize
        },
        line/.style = {draw=gray, thick, -{Stealth[length=3mm]}},
    ]
        \node[root] (root) {$h$-adaptive TSEM \\ computational framework};
        
        \node[module, below=1cm of root] (adapt) {$h$-adaptivity};
        \node[module, left=1cm of adapt] (space) {TSEM Space};
        \node[module, right=1cm of adapt] (solver) {Linear solver};
        
        \node[sub, below=1.4cm of space] (space1) {Discretization on \\ single tetrahedron};
        \node[sub, below=0.2cm of space1] (space2) {Assemble under \\ $C^{0}$ constraint};
        
        \node[sub, below=1.4cm of adapt] (adapt1) {Indicator};
        \node[sub, below=0.2cm of adapt1] (adapt2) {Red-Green Refinement};
        
        \node[sub, below=1.4cm of solver] (sol1) {Linear Equations};
        \node[sub, below=0.2cm of sol1] (sol2) {Eigenvalue Problem};
        
        \begin{scope}[on background layer]
            \node[draw=gray, dashed, thick, inner sep=6pt, rounded corners=4pt,
                  fit=(space1)(space2)] (m1) {};
            \node[draw=gray, dashed, thick, inner sep=6pt, rounded corners=4pt,
                  fit=(adapt1)(adapt2)] (m2) {};
            \node[draw=gray, dashed, thick, inner sep=6pt, rounded corners=4pt,
                  fit=(sol1)(sol2)] (m3) {};
        \end{scope}
        
        \draw[line] (root.south) -- ++(0,-0.3) -| (space.north);
        \draw[line] (root.south) -- ++(0,-0.3) -| (adapt.north);
        \draw[line] (root.south) -- ++(0,-0.3) -| (solver.north);
        
        \draw[line] (space.south) -- (m1.north);
        \draw[line] (adapt.south) -- (m2.north);
        \draw[line] (solver.south) -- (m3.north);
    \end{tikzpicture}
    \caption{The architecture for our framework.}
    \label{figure::architecture}
\end{figure}

After establishing the basis of TSEM and $h$-adaptivity, we now give an overview of our framework. The Figure \ref{figure::architecture} shows the architecture of our $h$-adaptive TSEM computational framework. The framework is decomposed into three principal modules:
\begin{enumerate}
    \item \textrm{TSEM Space}: performs the discretization on each single tetrahedron and assembles the global system subject to the global $C^{0}$ continuity constraint;
    \item \textrm{$h$-adaptivity}: contains the adaptivity logic with corresponding residual-type indicator and red-green refinement;
    \item \textrm{Linear Solver}: handles the algebraic problems arising from the discretization, including linear equations and eigenvalue problems.
\end{enumerate}
Next, in Section \ref{chapter::Implementation}, we elaborate on the implementation details of the relevant modules.

\section{\zw{Implementations and Numerical Tests}}
\label{chapter::Implementation}

Building upon the theoretical framework of the tetrahedral spectral element method and the $h$-adaptive refinement strategies established in the preceding sections, this section focuses on the practical implementation. 
In the following, we detail the software architecture and the specific algorithmic solutions, including the element orientation in Section \ref{sec::Orientation}, the data structure for the full topology data in Section \ref{sec::Tetrahedral-Mesh-with-Full-Topology-Data}, and the interpolation in Section \ref{sec::Interpolation}.

\subsection{\zw{Element Orientation}}
\label{sec::Orientation}

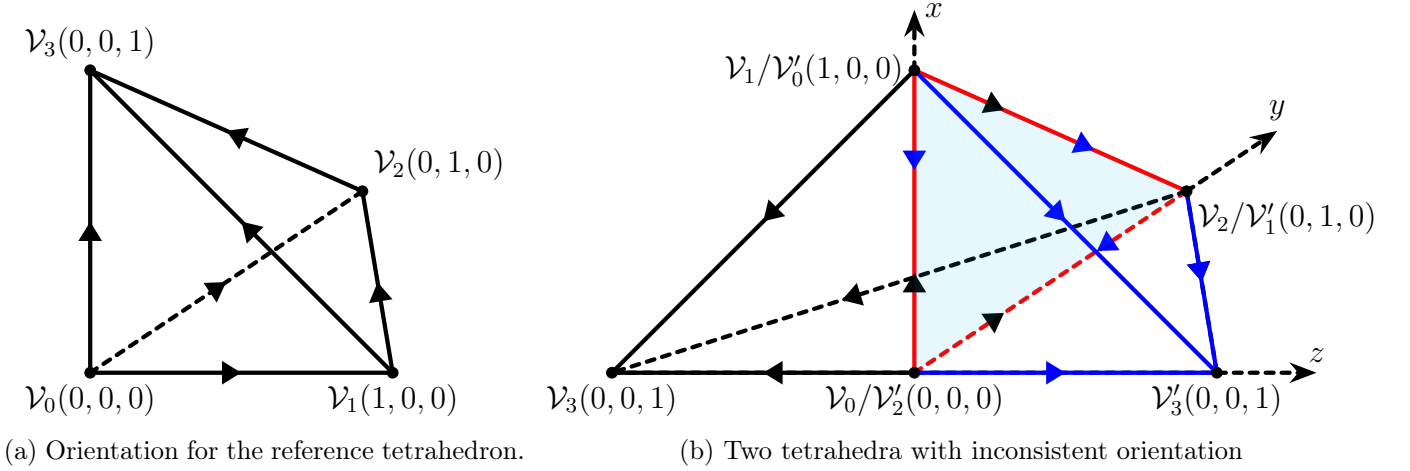
\begin{figure}[H]
    \centering
    \subfloat[Orientation for the reference tetrahedron.]
    {
    \begin{tikzpicture}[scale=4, line join=round, line cap=round]
        \coordinate (V0) at (0.0, 0.0);
        \coordinate (V1) at (1.0, 0.0);
        \coordinate (V2) at (0.9, 0.6);
        \coordinate (V3) at (0.0, 1.0);
            
        \draw[ultra thick, postaction={decorate, decoration={markings, mark=at position 0.5 with {\arrow{Triangle}}}}] (V0) -- (V1);
        \draw[ultra thick, dashed, postaction={decorate, decoration={markings, mark=at position 0.5 with {\arrow{Triangle}}}}] (V0) -- (V2);
        \draw[ultra thick, postaction={decorate, decoration={markings, mark=at position 0.5 with {\arrow{Triangle}}}}] (V0) -- (V3);
        \draw[ultra thick, postaction={decorate, decoration={markings, mark=at position 0.5 with {\arrow{Triangle}}}}] (V1) -- (V2);
        \draw[ultra thick, postaction={decorate, decoration={markings, mark=at position 0.5 with {\arrow{Triangle}}}}] (V1) -- (V3);
        \draw[ultra thick, postaction={decorate, decoration={markings, mark=at position 0.5 with {\arrow{Triangle}}}}] (V2) -- (V3);

        \filldraw (V0) circle (0.5pt) node[anchor=north] {$\mathcal{V}_{0} (0,0,0)$};
        \filldraw (V1) circle (0.5pt) node[anchor=north] {$\mathcal{V}_{1} (1,0,0)$};
        \filldraw (V2) circle (0.5pt) node[anchor=south west] {$\mathcal{V}_{2} (0,1,0)$};
        \filldraw (V3) circle (0.5pt) node[anchor=south]      {$\mathcal{V}_{3} (0,0,1)$};
    \end{tikzpicture}
    }
    \subfloat[Two tetrahedra with inconsistent orientation]
    {
    \begin{tikzpicture}[scale=4, line join=round, line cap=round]
        \coordinate (V0) at ( 0.0,  0.0);
        \coordinate (V1) at ( 0.0,  1.0);
        \coordinate (V2) at ( 0.9,  0.6);
        \coordinate (V3) at (-1.0,  0.0);
        \coordinate (V4) at ( 1.0,  0.0);
        
        \coordinate (X) at (0.0,  1.2);
        \coordinate (Y) at (1.2,  0.8);
        \coordinate (Z) at (1.33,  0.0);
            
        \draw[dashed, ultra thick, postaction={decorate, decoration={markings, mark=at position 0.5 with {\arrow{Triangle}}}}] (V2) -- (V4);
        \draw[dashed, ultra thick, -Stealth] (V1) -- (X);
        \draw[dashed, ultra thick, -Stealth] (V2) -- (Y);
        \draw[dashed, ultra thick, -Stealth] (V3) -- (Z);
        
        \draw[red, ultra thick, postaction={decorate, decoration={markings, mark=at position 0.33 with {\arrow[black]{Triangle}}}}] (V0) -- (V1);
        \draw[red, ultra thick, dashed, postaction={decorate, decoration={markings, mark=at position 0.33 with {\arrow[black]{Triangle}}}}] (V0) -- (V2);
        \draw[red, ultra thick, postaction={decorate, decoration={markings, mark=at position 0.33 with {\arrow[black]{Triangle}}}}] (V1) -- (V2);
        \path[red, ultra thick, postaction={decorate, decoration={markings, mark=at position 0.33 with {\arrow[blue]{Triangle}}}}] (V1) -- (V0);
        \path[red, ultra thick, dashed, postaction={decorate, decoration={markings, mark=at position 0.33 with {\arrow[blue]{Triangle}}}}] (V2) -- (V0);
        \path[red, ultra thick, postaction={decorate, decoration={markings, mark=at position 0.66 with {\arrow[blue]{Triangle}}}}] (V1) -- (V2);
        
        \draw[ultra thick, postaction={decorate, decoration={markings, mark=at position 0.5 with {\arrow{Triangle}}}}] (V0) -- (V3);
        \draw[ultra thick, postaction={decorate, decoration={markings, mark=at position 0.5 with {\arrow{Triangle}}}}] (V1) -- (V3);
        \draw[ultra thick, dashed, postaction={decorate, decoration={markings, mark=at position 0.6 with {\arrow{Triangle}}}}] (V2) -- (V3);
        
        \draw[blue, ultra thick, postaction={decorate, decoration={markings, mark=at position 0.5 with {\arrow{Triangle}}}}] (V0) -- (V4);
        \draw[blue, ultra thick, postaction={decorate, decoration={markings, mark=at position 0.5 with {\arrow{Triangle}}}}] (V1) -- (V4);
        \draw[blue, ultra thick, postaction={decorate, decoration={markings, mark=at position 0.5 with {\arrow{Triangle}}}}] (V2) -- (V4);

        \fill[cyan, opacity=0.1] (V0) -- (V1) -- (V2) -- cycle;
    
        \filldraw (V0) circle (0.5pt) node[anchor=north] {$\mathcal{V}_{0}/\mathcal{V}^{\prime}_{2} (0,0,0)$};
        \filldraw (V1) circle (0.5pt) node[anchor=east] {$\mathcal{V}_{1}/\mathcal{V}^{\prime}_{0} (1,0,0)$};
        \filldraw (V2) circle (0.5pt) node[anchor=north west] {$\mathcal{V}_{2}/\mathcal{V}^{\prime}_{1} (0,1,0)$};
        \filldraw (V3) circle (0.5pt) node[anchor=north]      {$\mathcal{V}_{3} (0,0,1)$};
        \filldraw (V4) circle (0.5pt) node[anchor=north]      {$\mathcal{V}^{\prime}_{3} (0,0,1)$};

        \filldraw (X) circle (0.0pt) node[anchor=west]      {$x$};
        \filldraw (Y) circle (0.0pt) node[anchor=south]      {$y$};
        \filldraw (Z) circle (0.0pt) node[anchor=south]      {$z$};
    \end{tikzpicture} \label{figure::twin-reference-tetrahedron::2}
    }
    \caption{The orientation for the reference tetrahedron and two tetrahedra with inconsistent orientation on the shared interface.}
    \label{figure::twin-reference-tetrahedron}
\end{figure}

In the tetrahedral spectral element method employed in this study, the construction of basis functions is explicitly dependent on the vertex ordering (i.e., the orientation) of each element. This dependency stems from the fact that the projections of edge and face modes onto their respective manifolds are defined by generalized Koornwinder polynomials. Specifically, Figure \ref{figure::twin-reference-tetrahedron} illustrates two tetrahedra with inconsistent orientations on face $\mathcal{V}_{0}\mathcal{V}_{1}\mathcal{V}_{2}$. The projection of the modal basis functions onto $\mathcal{F}_{2}$ yields:
\begin{equation}
    \begin{aligned}
        \phi_{i, j, 0}(x, y, 0) &= 2 (1 - y)^{i} J^{-1, -1}_{i}\left(\frac{2 x + y - 1}{1 - y}\right) J^{2 i - 1, -1}_{j}(2 y - 1), \\
        \phi^{\prime}_{i, j, 0}(x, y, 0) &= 2 (x + y)^{i} J^{-1, -1}_{i}\left(\frac{2 y - x}{x + y}\right) J^{2 i - 1, -1}_{j}(1 - 2 x - 2 y),
    \end{aligned}
\end{equation}
where $i \geq 2$ and $j \geq 1$. However, a change in triangular orientation induces two distinct sets of basis functions, which differ by a linear transformation \cite{UNKNOWN-CITE}.

Previous studies have developed two distinct strategies to ensure continuity. The first approach relies on constructing a feasible mesh based on a predefined template element \cite{zhan2023novel}. While more general tetrahedral meshes can be generated through the stacking of these template elements and homeomorphic mappings, the inherent topological constraints of such meshes preclude the use of $h$-adaptivity. An alternative approach \cite{UNKNOWN-CITE} treats the degrees of freedom on edges and faces independently across different elements, employing auxiliary linear mappings to enforce inter-element continuity. However, this method inevitably leads to a denser linear system, thereby incurring significant computational overhead in terms of both time and memory.

This work develops an algorithm for determining a consistent orientation for every tetrahedron within a general mesh. This consistency is formally established in Definitions \ref{def::orientation} and \ref{def::consistent-orientation}. Specifically, the intersection parameter $\vert \mathcal{S} \vert$ in Definition \ref{def::consistent-orientation} takes four discrete values $\{0, 1, 2, 3\}$, corresponding to intersections at the empty set, a vertex, an edge, or a face, respectively. The conditions in Definition \ref{def::consistent-orientation} ensure that all shared edges and faces maintain a uniform orientation across adjacent tetrahedra, thereby naturally satisfying the continuity requirements of the spectral element basis.

\begin{definition}
    \label{def::orientation}
    Let $\mathcal{T} = \{ \mathcal{V}_{0}, \mathcal{V}_{1}, \mathcal{V}_{2}, \mathcal{V}_{3} \}$ be a set of four non-coplanar vertices in $\mathbb{R}^3$ defining a tetrahedron. An orientation $\mathcal{O}$ of $\mathcal{T}$ is defined as an ordered sequence:
    \begin{equation}
        \mathcal{O} = (v_0, v_1, v_2, v_3)
    \end{equation}

    \noindent where the sequence is a specific permutation of the vertices in $\mathcal{T}$. This ordering establishes a unique bijective mapping from the physical tetrahedron to the vertices of a reference tetrahedron.
    


\end{definition}

\begin{definition}
    \label{def::consistent-orientation}
    Given a tetrahedral mesh $\mathcal{M}$ with tetrahedron $\{\mathcal{T}_j\}_{j=1}^m$ and associated orientations $O = \{ \mathcal{O}_{j} \}_{j=1}^{m}$, the orientation $O$ is said to be consistent if, for every pair of distinct tetrahedra $\mathcal{T}_{i}$ and $\mathcal{T}_{j}$ sharing a common set of vertices $\mathcal{S} = \mathcal{T}_{i} \cap \mathcal{T}_{j}$:
    
    \begin{enumerate}
        \item if $\vert \mathcal{S} \vert \leq 1$, the intersection is empty or a single vertex, then no additional constraint is imposed;
        \item if $\vert \mathcal{S} \vert > 1$, the intersection is an edge or a face, then the relative ordering of the vertices in $\mathcal{S}$ as they appear in sequence $\mathcal{O}_i$ must be identical to that in sequence $\mathcal{O}_j$.
    \end{enumerate}
    
\end{definition}

\begin{figure}[H]
    \centering
    \subfloat[Initialization at $\mathcal{V}_{0} \mathcal{V}_{3} \mathcal{V}_{5}$.]
    {
    \begin{tikzpicture}[scale=4, line join=round, line cap=round]
        \coordinate (V0) at ( 0.0,  0.0);
        \coordinate (V1) at ( 1.0,  0.0);
        \coordinate (V2) at ( 1.4,  0.2);
        \coordinate (V3) at ( 0.4,  0.2);
        \coordinate (V4) at ( 0.0,  1.0);
        \coordinate (V5) at ( 1.0,  1.0);
        \coordinate (V6) at ( 1.4,  1.2);
        \coordinate (V7) at ( 0.4,  1.2);

        \draw[ultra thick] (V0) -- (V1);
        \draw[ultra thick] (V1) -- (V2);
        \draw[ultra thick, dashed] (V2) -- (V3);
        \draw[brown, ultra thick, dashed, postaction={decorate, decoration={markings, mark=at position 0.5 with {\arrow{Triangle}}}}] (V0) -- (V3);
        
        \draw[ultra thick] (V4) -- (V5);
        \draw[ultra thick] (V5) -- (V6);
        \draw[ultra thick] (V6) -- (V7);
        \draw[ultra thick] (V4) -- (V7);

        \draw[ultra thick] (V0) -- (V4);
        \draw[ultra thick] (V1) -- (V5);
        \draw[ultra thick] (V2) -- (V6);
        \draw[ultra thick, dashed] (V3) -- (V7);

        \draw[dashed, ultra thick] (V1) -- (V3);
        \draw[dashed, ultra thick] (V3) -- (V4);
        \draw[brown, ultra thick, postaction={decorate, decoration={markings, mark=at position 0.5 with {\arrow{Triangle}}}}] (V0) -- (V5);
        \draw[ultra thick] (V2) -- (V5);
        \draw[ultra thick] (V5) -- (V7);
        \draw[brown, dashed, ultra thick, postaction={decorate, decoration={markings, mark=at position 0.5 with {\arrow{Triangle}}}}] (V3) -- (V5);
        \draw[dashed, ultra thick] (V3) -- (V6);

        \filldraw[brown] (V0) circle (0.5pt) node[anchor=north] {$\mathcal{V}_{0}$};
        \filldraw (V1) circle (0.5pt) node[anchor=north] {$\mathcal{V}_{1}$};
        \filldraw (V2) circle (0.5pt) node[anchor=north] {$\mathcal{V}_{2}$};
        \filldraw[brown] (V3) circle (0.5pt) node[anchor=north] {$\mathcal{V}_{3}$};
        \filldraw (V4) circle (0.5pt) node[anchor=south] {$\mathcal{V}_{4}$};
        \filldraw[brown] (V5) circle (0.5pt) node[anchor=south] {$\mathcal{V}_{5}$};
        \filldraw (V6) circle (0.5pt) node[anchor=south] {$\mathcal{V}_{6}$};
        \filldraw (V7) circle (0.5pt) node[anchor=south] {$\mathcal{V}_{7}$};
    \end{tikzpicture} \label{figure::ErrorOrientation::init}
    }
    \subfloat[Propagation to $\mathcal{V}_{1}$ and $\mathcal{V}_{4}$.]
    {
    \begin{tikzpicture}[scale=4, line join=round, line cap=round]
        \coordinate (V0) at ( 0.0,  0.0);
        \coordinate (V1) at ( 1.0,  0.0);
        \coordinate (V2) at ( 1.4,  0.2);
        \coordinate (V3) at ( 0.4,  0.2);
        \coordinate (V4) at ( 0.0,  1.0);
        \coordinate (V5) at ( 1.0,  1.0);
        \coordinate (V6) at ( 1.4,  1.2);
        \coordinate (V7) at ( 0.4,  1.2);

        \draw[brown, ultra thick, postaction={decorate, decoration={markings, mark=at position 0.5 with {\arrow{Triangle}}}}] (V0) -- (V1);
        \draw[ultra thick] (V1) -- (V2);
        \draw[ultra thick, dashed] (V3) -- (V2);
        \draw[ultra thick, dashed, postaction={decorate, decoration={markings, mark=at position 0.5 with {\arrow{Triangle}}}}] (V0) -- (V3);
        
        \draw[brown, ultra thick, postaction={decorate, decoration={markings, mark=at position 0.5 with {\arrow{Triangle}}}}] (V4) -- (V5);
        \draw[ultra thick] (V5) -- (V6);
        \draw[ultra thick] (V6) -- (V7);
        \draw[ultra thick] (V7) -- (V4);

        \draw[brown, ultra thick, postaction={decorate, decoration={markings, mark=at position 0.5 with {\arrow{Triangle}}}}] (V4) -- (V0);
        \draw[brown, ultra thick, postaction={decorate, decoration={markings, mark=at position 0.5 with {\arrow{Triangle}}}}] (V5) -- (V1);
        \draw[ultra thick] (V2) -- (V6);
        \draw[ultra thick] (V7) -- (V3);

        \draw[brown, dashed, ultra thick, postaction={decorate, decoration={markings, mark=at position 0.5 with {\arrow{Triangle}}}}] (V3) -- (V1);
        \draw[brown, dashed, ultra thick, postaction={decorate, decoration={markings, mark=at position 0.5 with {\arrow{Triangle}}}}] (V4) -- (V3);
        \draw[ultra thick, postaction={decorate, decoration={markings, mark=at position 0.5 with {\arrow{Triangle}}}}] (V0) -- (V5);
        \draw[ultra thick] (V5) -- (V2);
        \draw[ultra thick] (V7) -- (V5);
        \draw[dashed, ultra thick, postaction={decorate, decoration={markings, mark=at position 0.5 with {\arrow{Triangle}}}}] (V3) -- (V5);
        \draw[dashed, ultra thick] (V3) -- (V6);

        \path[ultra thick] (V2) -- (V6);
        \path[ultra thick] (V3) -- (V6);
        \path[ultra thick] (V5) -- (V6);

        \path[ultra thick] (V6) -- (V3);
        \path[ultra thick] (V6) -- (V5);
        \path[ultra thick] (V6) -- (V7);

        \filldraw (V0) circle (0.5pt) node[anchor=north] {$\mathcal{V}_{0}$};
        \filldraw[brown] (V1) circle (0.5pt) node[anchor=north] {$\mathcal{V}_{1}$};
        \filldraw (V2) circle (0.5pt) node[anchor=north] {$\mathcal{V}_{2}$};
        \filldraw (V3) circle (0.5pt) node[anchor=north] {$\mathcal{V}_{3}$};
        \filldraw[brown] (V4) circle (0.5pt) node[anchor=south] {$\mathcal{V}_{4}$};
        \filldraw (V5) circle (0.5pt) node[anchor=south] {$\mathcal{V}_{5}$};
        \filldraw (V6) circle (0.5pt) node[anchor=south] {$\mathcal{V}_{6}$};
        \filldraw (V7) circle (0.5pt) node[anchor=south] {$\mathcal{V}_{7}$};
    \end{tikzpicture} \label{figure::ErrorOrientation::1}
    } \\
    \subfloat[Propagation to $\mathcal{V}_{2}$ and $\mathcal{V}_{7}$.]
    {
    \begin{tikzpicture}[scale=4, line join=round, line cap=round]
        \coordinate (V0) at ( 0.0,  0.0);
        \coordinate (V1) at ( 1.0,  0.0);
        \coordinate (V2) at ( 1.4,  0.2);
        \coordinate (V3) at ( 0.4,  0.2);
        \coordinate (V4) at ( 0.0,  1.0);
        \coordinate (V5) at ( 1.0,  1.0);
        \coordinate (V6) at ( 1.4,  1.2);
        \coordinate (V7) at ( 0.4,  1.2);

        \draw[ultra thick, postaction={decorate, decoration={markings, mark=at position 0.5 with {\arrow{Triangle}}}}] (V0) -- (V1);
        \draw[brown, ultra thick, postaction={decorate, decoration={markings, mark=at position 0.5 with {\arrow{Triangle}}}}] (V1) -- (V2);
        \draw[brown, ultra thick, dashed, postaction={decorate, decoration={markings, mark=at position 0.5 with {\arrow{Triangle}}}}] (V3) -- (V2);
        \draw[ultra thick, dashed, postaction={decorate, decoration={markings, mark=at position 0.5 with {\arrow{Triangle}}}}] (V0) -- (V3);
        
        \draw[ultra thick, postaction={decorate, decoration={markings, mark=at position 0.5 with {\arrow{Triangle}}}}] (V4) -- (V5);
        \draw[ultra thick] (V5) -- (V6);
        \draw[ultra thick] (V6) -- (V7);
        \draw[brown, ultra thick, postaction={decorate, decoration={markings, mark=at position 0.5 with {\arrow{Triangle}}}}] (V7) -- (V4);

        \draw[ultra thick, postaction={decorate, decoration={markings, mark=at position 0.5 with {\arrow{Triangle}}}}] (V4) -- (V0);
        \draw[ultra thick, postaction={decorate, decoration={markings, mark=at position 0.5 with {\arrow{Triangle}}}}] (V5) -- (V1);
        \draw[ultra thick] (V2) -- (V6);
        \draw[brown, ultra thick, dashed, postaction={decorate, decoration={markings, mark=at position 0.5 with {\arrow{Triangle}}}}] (V7) -- (V3);

        \draw[dashed, ultra thick, postaction={decorate, decoration={markings, mark=at position 0.5 with {\arrow{Triangle}}}}] (V3) -- (V1);
        \draw[dashed, ultra thick, postaction={decorate, decoration={markings, mark=at position 0.5 with {\arrow{Triangle}}}}] (V4) -- (V3);
        \draw[ultra thick, postaction={decorate, decoration={markings, mark=at position 0.5 with {\arrow{Triangle}}}}] (V0) -- (V5);
        \draw[brown, ultra thick, postaction={decorate, decoration={markings, mark=at position 0.5 with {\arrow{Triangle}}}}] (V5) -- (V2);
        \draw[brown, ultra thick, postaction={decorate, decoration={markings, mark=at position 0.5 with {\arrow{Triangle}}}}] (V7) -- (V5);
        \draw[dashed, ultra thick, postaction={decorate, decoration={markings, mark=at position 0.5 with {\arrow{Triangle}}}}] (V3) -- (V5);
        \draw[dashed, ultra thick] (V3) -- (V6);

        \path[ultra thick] (V2) -- (V6);
        \path[ultra thick] (V3) -- (V6);
        \path[ultra thick] (V5) -- (V6);

        \path[ultra thick] (V6) -- (V3);
        \path[ultra thick] (V6) -- (V5);
        \path[ultra thick] (V6) -- (V7);

        \filldraw (V0) circle (0.5pt) node[anchor=north] {$\mathcal{V}_{0}$};
        \filldraw (V1) circle (0.5pt) node[anchor=north] {$\mathcal{V}_{1}$};
        \filldraw (V2)[brown] circle (0.5pt) node[anchor=north] {$\mathcal{V}_{2}$};
        \filldraw (V3) circle (0.5pt) node[anchor=north] {$\mathcal{V}_{3}$};
        \filldraw (V4) circle (0.5pt) node[anchor=south] {$\mathcal{V}_{4}$};
        \filldraw (V5) circle (0.5pt) node[anchor=south] {$\mathcal{V}_{5}$};
        \filldraw (V6) circle (0.5pt) node[anchor=south] {$\mathcal{V}_{6}$};
        \filldraw (V7)[brown] circle (0.5pt) node[anchor=south] {$\mathcal{V}_{7}$};
    \end{tikzpicture} \label{figure::ErrorOrientation::2}
    }
    \subfloat[Contradiction at $\mathcal{V}_{6}$.]
    {
    \begin{tikzpicture}[scale=4, line join=round, line cap=round]
        \coordinate (V0) at ( 0.0,  0.0);
        \coordinate (V1) at ( 1.0,  0.0);
        \coordinate (V2) at ( 1.4,  0.2);
        \coordinate (V3) at ( 0.4,  0.2);
        \coordinate (V4) at ( 0.0,  1.0);
        \coordinate (V5) at ( 1.0,  1.0);
        \coordinate (V6) at ( 1.4,  1.2);
        \coordinate (V7) at ( 0.4,  1.2);

        \draw[ultra thick, postaction={decorate, decoration={markings, mark=at position 0.5 with {\arrow{Triangle}}}}] (V0) -- (V1);
        \draw[ultra thick, postaction={decorate, decoration={markings, mark=at position 0.5 with {\arrow{Triangle}}}}] (V1) -- (V2);
        \draw[ultra thick, dashed, postaction={decorate, decoration={markings, mark=at position 0.5 with {\arrow{Triangle}}}}] (V3) -- (V2);
        \draw[ultra thick, dashed, postaction={decorate, decoration={markings, mark=at position 0.5 with {\arrow{Triangle}}}}] (V0) -- (V3);
        
        \draw[ultra thick, postaction={decorate, decoration={markings, mark=at position 0.5 with {\arrow{Triangle}}}}] (V4) -- (V5);
        \draw[brown, ultra thick] (V5) -- (V6);
        \draw[brown, ultra thick] (V6) -- (V7);
        \draw[ultra thick, postaction={decorate, decoration={markings, mark=at position 0.5 with {\arrow{Triangle}}}}] (V7) -- (V4);

        \draw[ultra thick, postaction={decorate, decoration={markings, mark=at position 0.5 with {\arrow{Triangle}}}}] (V4) -- (V0);
        \draw[ultra thick, postaction={decorate, decoration={markings, mark=at position 0.5 with {\arrow{Triangle}}}}] (V5) -- (V1);
        \draw[brown, ultra thick] (V2) -- (V6);
        \draw[ultra thick, dashed, postaction={decorate, decoration={markings, mark=at position 0.5 with {\arrow{Triangle}}}}] (V7) -- (V3);

        \draw[dashed, ultra thick, postaction={decorate, decoration={markings, mark=at position 0.5 with {\arrow{Triangle}}}}] (V3) -- (V1);
        \draw[dashed, ultra thick, postaction={decorate, decoration={markings, mark=at position 0.5 with {\arrow{Triangle}}}}] (V4) -- (V3);
        \draw[ultra thick, postaction={decorate, decoration={markings, mark=at position 0.5 with {\arrow{Triangle}}}}] (V0) -- (V5);
        \draw[ultra thick, postaction={decorate, decoration={markings, mark=at position 0.5 with {\arrow{Triangle}}}}] (V5) -- (V2);
        \draw[ultra thick, postaction={decorate, decoration={markings, mark=at position 0.5 with {\arrow{Triangle}}}}] (V7) -- (V5);
        \draw[dashed, ultra thick, postaction={decorate, decoration={markings, mark=at position 0.5 with {\arrow{Triangle}}}}] (V3) -- (V5);
        \draw[brown, dashed, ultra thick] (V3) -- (V6);

        \path[ultra thick, postaction={decorate, decoration={markings, mark=at position 0.5 with {\arrow[blue]{Triangle}}}}] (V2) -- (V6);
        \path[ultra thick, postaction={decorate, decoration={markings, mark=at position 0.33 with {\arrow[blue]{Triangle}}}}] (V3) -- (V6);
        \path[ultra thick, postaction={decorate, decoration={markings, mark=at position 0.33 with {\arrow[blue]{Triangle}}}}] (V5) -- (V6);

        \path[ultra thick, postaction={decorate, decoration={markings, mark=at position 0.33 with {\arrow[red]{Triangle}}}}] (V6) -- (V3);
        \path[ultra thick, postaction={decorate, decoration={markings, mark=at position 0.33 with {\arrow[red]{Triangle}}}}] (V6) -- (V5);
        \path[ultra thick, postaction={decorate, decoration={markings, mark=at position 0.5 with {\arrow[red]{Triangle}}}}] (V6) -- (V7);

        \filldraw (V0) circle (0.5pt) node[anchor=north] {$\mathcal{V}_{0}$};
        \filldraw (V1) circle (0.5pt) node[anchor=north] {$\mathcal{V}_{1}$};
        \filldraw (V2) circle (0.5pt) node[anchor=north] {$\mathcal{V}_{2}$};
        \filldraw (V3) circle (0.5pt) node[anchor=north] {$\mathcal{V}_{3}$};
        \filldraw (V4) circle (0.5pt) node[anchor=south] {$\mathcal{V}_{4}$};
        \filldraw (V5) circle (0.5pt) node[anchor=south] {$\mathcal{V}_{5}$};
        \filldraw (V6)[brown] circle (0.5pt) node[anchor=south] {$\mathcal{V}_{6}$};
        \filldraw (V7) circle (0.5pt) node[anchor=south] {$\mathcal{V}_{7}$};
    \end{tikzpicture} \label{figure::ErrorOrientation::error}
    }
    \caption{The propagation based on the right-hand rule, where (a) shows the initial coloring, (b) and (c) show the propagation to new vertices, and (d) shows the contradiction at $\mathcal{V}_{6}$.}
    \label{figure::ErrorOrientation}
\end{figure}
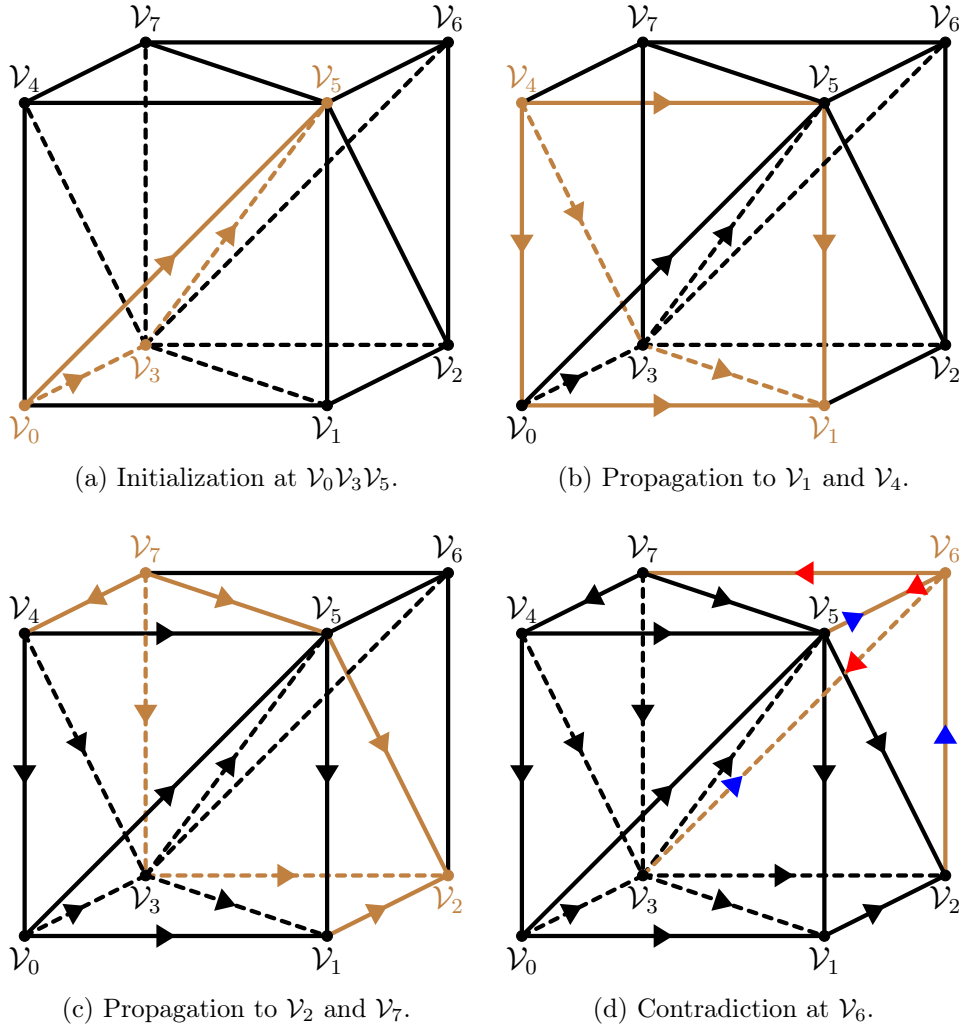

However, implementing this consistent orientation at the vertex level presents significant challenges. A natural strategy is to start with one face and then propagate the orientation according to the right-hand rule. However, although this propagation ensures a positive Jacobian determinant for each individual tetrahedron, such a scheme is not universally applicable and often fails to achieve global consistency. For the example illustrated in Figure \ref{figure::ErrorOrientation}, we apply the strategy by:
\begin{enumerate}
    \item Starting with an initial tetrahedron (e.g., $\mathcal{V}_{0} \mathcal{V}_{3} \mathcal{V}_{5}$), with an arbitrary orientation;
    \item For face $\mathcal{V}_{0} \mathcal{V}_{3} \mathcal{V}_{5}$, the normal vector gives the direction $\overrightarrow{\mathcal{V}_{4} \mathcal{V}_{1}}$, then the orientation of its neighbor is given accordingly (e.g., $\mathcal{V}_{4} \to \mathcal{V}_{0} \to \mathcal{V}_{1}$);
    \item For faces $\mathcal{V}_{4} \mathcal{V}_{3} \mathcal{V}_{5}$ and $\mathcal{V}_{3} \mathcal{V}_{5} \mathcal{V}_{1}$, the orientations of two faces and the newly encountered vertices $\mathcal{V}_{2}$ and $\mathcal{V}_{7}$ are determined similarly;
    \item With $\mathcal{V}_{6}$ left as the final vertex, a contradiction arises regarding the edges $\mathcal{V}_{5}\mathcal{V}_{6}$ and $\mathcal{V}_{3}\mathcal{V}_{6}$ since the orientations propagated from different adjacent paths assign a conflicting direction.
\end{enumerate}
Evidently, this strategy may fail even in simple cube cases. Consequently, we propose a generalized orientation scheme that circumvents the rigid reliance on predefined Jacobian positivity, rendering it applicable to arbitrary and complex mesh topologies.


We develop a systematic strategy to address the challenges of establishing a globally consistent orientation. In particular, we first investigate feasible ordering relations on vertices and then derive the global element orientations based on these relations.
As illustrated in Definition \ref{def::orientation}, the canonical orientation of a reference tetrahedron inherently induces a total ordering of its local vertices. Therefore, if a global total ordering is established for the entire mesh, any sub-simplex, such as an edge or a face, can be uniquely and consistently oriented according to the relative ranks of its constituent vertices.

\begin{figure}[H]
    \centering
    \subfloat[The orientation based on the indices.]
    {
    \begin{tikzpicture}[scale=4, line join=round, line cap=round]
        \coordinate (V0) at ( 0.0,  0.0);
        \coordinate (V1) at ( 1.0,  0.0);
        \coordinate (V2) at ( 1.4,  0.2);
        \coordinate (V3) at ( 0.4,  0.2);
        \coordinate (V4) at ( 0.0,  1.0);
        \coordinate (V5) at ( 1.0,  1.0);
        \coordinate (V6) at ( 1.4,  1.2);
        \coordinate (V7) at ( 0.4,  1.2);

        \draw[red, ultra thick, postaction={decorate, decoration={markings, mark=at position 0.5 with {\arrow{Triangle}}}}] (V0) -- (V1);
        \draw[red, ultra thick, postaction={decorate, decoration={markings, mark=at position 0.5 with {\arrow{Triangle}}}}] (V1) -- (V2);
        \draw[red, ultra thick, dashed, postaction={decorate, decoration={markings, mark=at position 0.5 with {\arrow{Triangle}}}}] (V2) -- (V3);
        \draw[red, ultra thick, dashed, postaction={decorate, decoration={markings, mark=at position 0.5 with {\arrow{Triangle}}}}] (V0) -- (V3);
        
        \draw[blue, ultra thick, postaction={decorate, decoration={markings, mark=at position 0.5 with {\arrow{Triangle}}}}] (V4) -- (V5);
        \draw[blue, ultra thick, postaction={decorate, decoration={markings, mark=at position 0.5 with {\arrow{Triangle}}}}] (V5) -- (V6);
        \draw[blue, ultra thick, postaction={decorate, decoration={markings, mark=at position 0.5 with {\arrow{Triangle}}}}] (V7) -- (V6);
        \draw[blue, ultra thick, postaction={decorate, decoration={markings, mark=at position 0.5 with {\arrow{Triangle}}}}] (V7) -- (V4);

        \draw[brown, ultra thick, postaction={decorate, decoration={markings, mark=at position 0.5 with {\arrow{Triangle}}}}] (V0) -- (V4);
        \draw[brown, ultra thick, postaction={decorate, decoration={markings, mark=at position 0.5 with {\arrow{Triangle}}}}] (V1) -- (V5);
        \draw[brown, ultra thick, postaction={decorate, decoration={markings, mark=at position 0.5 with {\arrow{Triangle}}}}] (V2) -- (V6);
        \draw[brown, ultra thick, dashed, postaction={decorate, decoration={markings, mark=at position 0.5 with {\arrow{Triangle}}}}] (V3) -- (V7);

        \draw[red, dashed, ultra thick, postaction={decorate, decoration={markings, mark=at position 0.5 with {\arrow{Triangle}}}}] (V1) -- (V3);
        \draw[brown, dashed, ultra thick, postaction={decorate, decoration={markings, mark=at position 0.5 with {\arrow{Triangle}}}}] (V3) -- (V4);
        \draw[brown, ultra thick, postaction={decorate, decoration={markings, mark=at position 0.5 with {\arrow{Triangle}}}}] (V0) -- (V5);
        \draw[brown, ultra thick, postaction={decorate, decoration={markings, mark=at position 0.5 with {\arrow{Triangle}}}}] (V2) -- (V5);
        \draw[blue, ultra thick, postaction={decorate, decoration={markings, mark=at position 0.5 with {\arrow{Triangle}}}}] (V7) -- (V5);
        \draw[brown, dashed, ultra thick, postaction={decorate, decoration={markings, mark=at position 0.5 with {\arrow{Triangle}}}}] (V3) -- (V5);
        \draw[brown, dashed, ultra thick, postaction={decorate, decoration={markings, mark=at position 0.5 with {\arrow{Triangle}}}}] (V3) -- (V6);

        \filldraw (V0) circle (0.5pt) node[anchor=north] {$\mathcal{V}_{0}$};
        \filldraw (V1) circle (0.5pt) node[anchor=north] {$\mathcal{V}_{1}$};
        \filldraw (V2) circle (0.5pt) node[anchor=north] {$\mathcal{V}_{2}$};
        \filldraw (V3) circle (0.5pt) node[anchor=north] {$\mathcal{V}_{3}$};
        \filldraw (V4) circle (0.5pt) node[anchor=south] {$\mathcal{V}_{5}$};
        \filldraw (V5) circle (0.5pt) node[anchor=south] {$\mathcal{V}_{6}$};
        \filldraw (V6) circle (0.5pt) node[anchor=south] {$\mathcal{V}_{7}$};
        \filldraw (V7) circle (0.5pt) node[anchor=south east] {$\mathcal{V}_{4}$};
        
        \coordinate (X)  at (-0.20, -0.10);
        \coordinate (Y)  at ( 1.60,  0.20);
        \coordinate (Z)  at ( 0.40,  1.50);
        
        \draw[dashed, ultra thick, -Stealth] (V0) -- (X);
        \draw[dashed, ultra thick, -Stealth] (V2) -- (Y);
        \draw[dashed, ultra thick, -Stealth] (V7) -- (Z);
        \filldraw (X) circle (0.0pt) node[anchor=north]      {$x$};
        \filldraw (Y) circle (0.0pt) node[anchor=south]      {$y$};
        \filldraw (Z) circle (0.0pt) node[anchor=east]      {$z$};

        \fill[cyan, opacity=0.1] (V0) -- (V1) -- (V2) -- (V3) -- cycle;
        \fill[cyan, opacity=0.1] (V4) -- (V5) -- (V6) -- (V7) -- cycle;
    \end{tikzpicture} \label{figure::Orientation::index}
    }
    \subfloat[The orientation based on the vertex positions.]
    {
    \begin{tikzpicture}[scale=4, line join=round, line cap=round]
        \coordinate (V0) at ( 0.0,  0.0);
        \coordinate (V1) at ( 1.0,  0.0);
        \coordinate (V2) at ( 1.4,  0.2);
        \coordinate (V3) at ( 0.4,  0.2);
        \coordinate (V4) at ( 0.0,  1.0);
        \coordinate (V5) at ( 1.0,  1.0);
        \coordinate (V6) at ( 1.4,  1.2);
        \coordinate (V7) at ( 0.4,  1.2);

        \draw[brown, ultra thick, postaction={decorate, decoration={markings, mark=at position 0.5 with {\arrow{Triangle}}}}] (V0) -- (V1);
        \draw[red, ultra thick, postaction={decorate, decoration={markings, mark=at position 0.5 with {\arrow{Triangle}}}}] (V2) -- (V1);
        \draw[brown, ultra thick, dashed, postaction={decorate, decoration={markings, mark=at position 0.5 with {\arrow{Triangle}}}}] (V3) -- (V2);
        \draw[red, ultra thick, dashed, postaction={decorate, decoration={markings, mark=at position 0.5 with {\arrow{Triangle}}}}] (V3) -- (V0);
        
        \draw[brown, ultra thick, postaction={decorate, decoration={markings, mark=at position 0.5 with {\arrow{Triangle}}}}] (V4) -- (V5);
        \draw[red, ultra thick, postaction={decorate, decoration={markings, mark=at position 0.5 with {\arrow{Triangle}}}}] (V6) -- (V5);
        \draw[brown, ultra thick, postaction={decorate, decoration={markings, mark=at position 0.5 with {\arrow{Triangle}}}}] (V7) -- (V6);
        \draw[red, ultra thick, postaction={decorate, decoration={markings, mark=at position 0.5 with {\arrow{Triangle}}}}] (V7) -- (V4);

        \draw[blue, ultra thick, postaction={decorate, decoration={markings, mark=at position 0.5 with {\arrow{Triangle}}}}] (V0) -- (V4);
        \draw[blue, ultra thick, postaction={decorate, decoration={markings, mark=at position 0.5 with {\arrow{Triangle}}}}] (V1) -- (V5);
        \draw[blue, ultra thick, postaction={decorate, decoration={markings, mark=at position 0.5 with {\arrow{Triangle}}}}] (V2) -- (V6);
        \draw[blue, ultra thick, dashed, postaction={decorate, decoration={markings, mark=at position 0.5 with {\arrow{Triangle}}}}] (V3) -- (V7);

        \draw[red, dashed, ultra thick, postaction={decorate, decoration={markings, mark=at position 0.5 with {\arrow{Triangle}}}}] (V3) -- (V1);
        \draw[red, dashed, ultra thick, postaction={decorate, decoration={markings, mark=at position 0.5 with {\arrow{Triangle}}}}] (V3) -- (V4);
        \draw[brown, ultra thick, postaction={decorate, decoration={markings, mark=at position 0.5 with {\arrow{Triangle}}}}] (V0) -- (V5);
        \draw[red, ultra thick, postaction={decorate, decoration={markings, mark=at position 0.5 with {\arrow{Triangle}}}}] (V2) -- (V5);
        \draw[red, ultra thick, postaction={decorate, decoration={markings, mark=at position 0.5 with {\arrow{Triangle}}}}] (V7) -- (V5);
        \draw[red, dashed, ultra thick, postaction={decorate, decoration={markings, mark=at position 0.5 with {\arrow{Triangle}}}}] (V3) -- (V5);
        \draw[brown, dashed, ultra thick, postaction={decorate, decoration={markings, mark=at position 0.5 with {\arrow{Triangle}}}}] (V3) -- (V6);

        \filldraw (V0) circle (0.5pt) node[anchor=north] {$\mathcal{V}_{0}$};
        \filldraw (V1) circle (0.5pt) node[anchor=north] {$\mathcal{V}_{1}$};
        \filldraw (V2) circle (0.5pt) node[anchor=north] {$\mathcal{V}_{2}$};
        \filldraw (V3) circle (0.5pt) node[anchor=north] {$\mathcal{V}_{3}$};
        \filldraw (V4) circle (0.5pt) node[anchor=south] {$\mathcal{V}_{5}$};
        \filldraw (V5) circle (0.5pt) node[anchor=south] {$\mathcal{V}_{6}$};
        \filldraw (V6) circle (0.5pt) node[anchor=south] {$\mathcal{V}_{7}$};
        \filldraw (V7) circle (0.5pt) node[anchor=south east] {$\mathcal{V}_{4}$};
        
        \coordinate (X)  at (-0.20, -0.10);
        \coordinate (Y)  at ( 1.60,  0.20);
        \coordinate (Z)  at ( 0.40,  1.50);
        
        \draw[dashed, ultra thick, -Stealth] (V0) -- (X);
        \draw[dashed, ultra thick, -Stealth] (V2) -- (Y);
        \draw[dashed, ultra thick, -Stealth] (V7) -- (Z);
        \filldraw (X) circle (0.0pt) node[anchor=north]      {$x$};
        \filldraw (Y) circle (0.0pt) node[anchor=south]      {$y$};
        \filldraw (Z) circle (0.0pt) node[anchor=east]      {$z$};
        
        \fill[cyan, opacity=0.1] (V0) -- (V1) -- (V2) -- (V3) -- cycle;
        \fill[cyan, opacity=0.1] (V4) -- (V5) -- (V6) -- (V7) -- cycle;
    \end{tikzpicture} \label{figure::Orientation::position}
    }
    \caption{The example of orientations, where (a) shows the orientation based on the vertex indices, and (b) shows the orientation based on the vertex position. For (b), the red, brown and blue edges are determined by $x$, $y$ and $z$ respectively, and the orientation of top and bottom surface shows the translation invariance.}
    \label{figure::Orientation}
\end{figure}
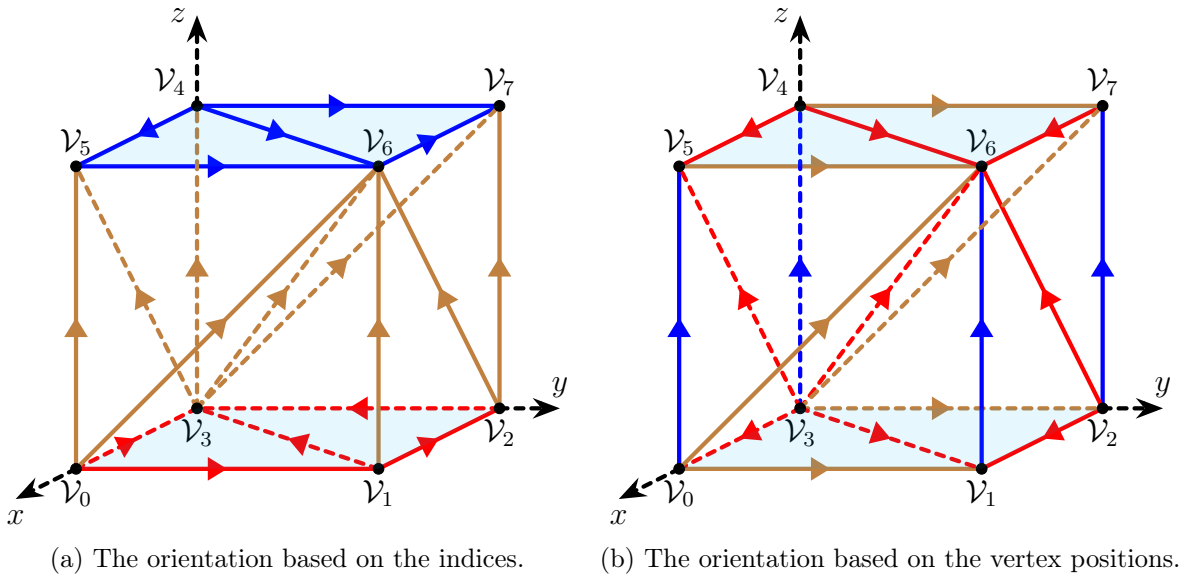

We explore two distinct strategies for implementing this total ordering framework, as compared in Figure \ref{figure::Orientation}. The first approach, shown in Figure \ref{figure::Orientation}\subref{figure::Orientation::index}, utilizes global vertex indices directly, which can also be found in previous work \cite{ainsworth2003hierarchic}. While straightforward, this method lacks robustness in multi-block meshes or domain decomposition scenarios, where local index renumbering can lead to orientation mismatches across partitions. Furthermore, in $h$-adaptive refinement, the removal or addition of vertices often triggers a re-indexing process. Such volatility flips the tetrahedral orientations and necessitates complex ``re-orientation'' operators to correct the affected basis functions.

The second strategy, depicted in Figure \ref{figure::Orientation}\subref{figure::Orientation::position}, employs a lexicographical ordering based on spatial coordinates $(x, y, z)$. This criterion rigorously satisfies the mathematical requirements of a total ordering, i.e., antisymmetric, transitive, and total, ensuring that any two distinct vertices are uniquely comparable. The primary advantage of this coordinate-based approach is its geometric invariance: unlike volatile indices, spatial positions remain unchanged during mesh migration and refinement. For example, for the top and bottom surfaces, the position-based method shows translation invariance, which is not held for the index-based method. This stability ensures that the relative ordering of vertices, and thus the orientation of the tetrahedra, is kept invariant throughout the simulation. Consequently, by anchoring the orientation to an invariant total ordering, continuity is enforced by construction when degrees of freedom are assigned to their respective geometric entities, eliminating the overhead of auxiliary linear transformations.

\subsection{\zw{Tetrahedral Mesh Data Structure with Complete Topology}}
\label{sec::Tetrahedral-Mesh-with-Full-Topology-Data}

A tetrahedral mesh with full topological information is indispensable for our $h$-adaptive spectral element method. Specifically, three fundamental challenges must be addressed: (1) matrix assembly requires ``downward'' adjacency, i.e., the mapping from a tetrahedron to its constituent faces, edges, and vertices; (2) red-green refinement necessitates ``upward'' adjacency, i.e., the connectivity from edges or faces back to their incident tetrahedra; and (3) inter-grid interpolation requires the maintenance of the entire refinement hierarchy across multiple levels. Prominent finite-element libraries such as deal.II \cite{2025:arndt.bangerth.ea:deal} and MFEM \cite{mfem-2024} typically employ algebraic constraints to manage hanging nodes, which prioritize memory efficiency by storing only fixed-size ``downward'' relations. However, this approach is ill-suited for unstructured tetrahedral SEM. Given the high order of basis polynomials and the intricate orientation issues involved, constructing interpolation matrices via constraints becomes prohibitively complex. Consequently, we adopt the red-green refinement strategy to circumvent these constraints, necessitating the design and implementation of a custom mesh data structure that supports full topological queries.

Our data structure encapsulates comprehensive geometric and topological information by assigning a unique identifier to every entity, including vertices, edges, faces, and tetrahedra. Beyond standard geometric data, such as nodal coordinates and precomputed Jacobian matrices, each entity maintains multiple hash-based lookup tables. These tables define mappings to both constituent lower-dimensional entities (e.g., the vertices and edges defining a face) and adjacent entities of equal or higher dimension (e.g., neighboring faces or incident cells). This architectural design ensures that any adjacency query between two entities is resolved via direct indexing or constant-time lookup, achieving a time complexity of $O(1)$.

The refinement tree constitutes another pivotal component of our $h$-adaptive framework. In the context of red-green refinement, ``green'' elements are conventionally treated as temporary transitions to ensure mesh conformity and are restricted from being further bisected. When employing algebraic constraints to manage hanging nodes, these green elements are entirely bypassed. Consequently, traditional hierarchical representations, such as octrees or quadtrees, are optimized for ``red'' refinement and typically exclude green elements from the formal tree hierarchy. However, such an exclusion introduces topological inconsistencies during element traversal and necessitates specialized handling that incurs additional computational overhead. To address this, our refinement tree explicitly incorporates green elements, ensuring a bijective mapping where every physical element corresponds to a unique leaf node. This design is primarily motivated by the requirements of efficient interpolation: by maintaining a comprehensive parent-child lineage, we enable uniform and rapid traversal across all entities—a critical prerequisite for the prolongation operators in multi-level adaptive cycles.

Furthermore, the availability of full topological data enables a robust integrity verification mechanism. To ensure the correctness of the adaptive process, we implement a topological consistency check based on the Euler characteristic. Specifically, the algorithm computes the number of connected components of the mesh ($N_{\text{mesh}}$), the boundary components ($N_{\text{boundary}}$), and the Euler characteristic $\chi = N_{\text{vert}} - N_{\text{edge}} + N_{\text{face}} - N_{\text{cell}}$. From these, we derive the number of handles (genus) $N_{\text{handle}} = N_{\text{boundary}} - \chi$ and cavities $N_{\text{cavity}} = N_{\text{boundary}} - N_{\text{mesh}}$. By validating these invariants against the analytical topology of the domain, we guarantee that the $h$-adaptive refinement preserves the original manifold topology without introducing spurious connectivity errors.

\subsection{Interpolation}
\label{sec::Interpolation}


In $h$-adaptive simulations, transferring solution fields from a coarse mesh to a refined mesh is a fundamental operation. Leveraging our fully topological mesh data structure, we further optimize the transfer method to achieve higher efficiency while preserving accuracy. Our approach adopts the idea of evaluating the source solution at the quadrature points of the refined grid. However, high-order quadrature rules can lead to a large number of evaluation points, making the associated location step rather expensive. To overcome this, we exploit a key topological property of red-green refinement: for any given point, the coarse-mesh cell and the fine-mesh cell that contain it are guaranteed to share a common ancestor in the refinement tree. Moreover, due to the rule of red-green refinement, this common ancestor is either the fine-mesh cell itself or the direct parent of the fine-mesh cell. Consequently, the distance in the refinement tree to this ancestor is at most one, which drastically accelerates the location of quadrature points.


Based on this property, we precompute a mapping from each refined cell to its common coarse ancestor. This reduces the location cost from one search per evaluation point to a single cell-level lookup shared by all quadrature points within that cell. Then, as outlined in Algorithms \ref{algo::interpolation} and \ref{algo::construct-function}, the interpolation can be carried out locally within each element since entities of the same dimension remain mutually independent. The overall computational complexity of this approach scales linearly with the total number of tetrahedra and the order of the basis polynomials, ensuring high efficiency for large-scale adaptive cycles.


\begin{algorithm}[H]
    \caption{Interpolation from $M_{c}$ to $M_{f}$ using precomputed ancestor mapping}
    \label{algo::interpolation}
    \KwIn{Coarse mesh: $M_{c}$; Fine mesh: $M_{f}$; Coefficients of basis functions on $M_{c}$: $\mathbf{u}_{c}$; Refinement tree of coarse mesh: $\mathbf{T}_{c}$; Precomputed mapping $\pi: M_{f} \to \mathbf{T}_{c}$ associating each fine cell with its containing coarse ancestor.
    }
    \KwOut{Coefficients of basis functions on $M_{f}$: $\mathbf{u}_{f}$.}
    \BlankLine
    \For {tetrahedra $\mathcal{T} \in M_{f}$}
    {
        \eIf{$\mathcal{T} \in M_{c}$}
        {
            \tcc{Not refined: directly copy coefficients}
            $\mathcal{T}$: $\left. \mathbf{u}_{f} \right|_\mathcal{T} \leftarrow \left. \mathbf{u}_{c} \right|_\mathcal{T}$\;
        }
        {
            \tcc{Refined: interpolate from its coarse ancestor}
            Retrieve the coarse ancestor $\mathcal{T}_{p} \leftarrow \pi(\mathcal{T}) \in \mathbf{T}_{c}$\;
            Construct $u_{\mathcal{T}_{p}}(\mathbf{x})$ that returns the coarse solution inside $\mathcal{T}_{p}$ through Algorithm \ref{algo::construct-function}\;
            Interpolate DoFs on all entities of $\mathcal{T}$: $\left. \mathbf{u}_{f} \right|_\mathcal{T} \leftarrow \text{interpolate}_{\mathcal{T}}(u_{\mathcal{T}_{p}}(\mathbf{x}))$\;
        }
    }
\end{algorithm}

\begin{algorithm}[H]
    \caption{Evaluate the coarse solution at a point inside a coarse tetrahedron.}
    \label{algo::construct-function}
    \KwIn{Coarse ancestor tetrahedron: $\mathcal{T}_{p}$; Refinement subtree rooted at $\mathcal{T}_{p}$: $\mathbf{T}_{\mathcal{T}_{p}}$; Coefficients of basis polynomials on the coarse mesh: $\mathbf{u}_{c}$.}
    \KwOut{A function on $\mathcal{T}_{p}$: $u_{\mathcal{T}_{p}}(\mathbf{x})$.}
    \BlankLine
    For a point $\mathbf{x} \in \mathcal{T}_{p}$: \\
    Locate the unique leaf cell $\mathcal{T}_{c} \in \mathbf{T}_{\mathcal{T}_{p}}$ that contains $\mathbf{x}$\;
    \tcc{The height of the subtree is at most $1$, thus this locating step is $O(1)$}
    Extract the coefficients $\mathbf{u}_{\mathcal{T}_{c}} \leftarrow \left. \mathbf{u}_{c} \right|_{\mathcal{T}_{c}}$\;
    Compute the solution $u_{\mathcal{T}_{p}}(\mathbf{x}) \leftarrow \sum_{\mathbf{i}} \mathbf{u}_{\mathcal{T}_{c}}^{\mathbf{i}} \phi_{\mathbf{i}}(\mathbf{x})$ where $\phi_{\mathbf{i}}(\mathbf{x})$ are modal basis functions on $\mathcal{T}_{c}$;
\end{algorithm}

\subsection{\zw{Numerical Verification}}

\zw{The performance of the proposed adaptive tetrahedral spectral element method is assessed through a set of numerical experiments. In particular, the Poisson equation is first considered on a cubic domain and an L-shaped prismatic domain to verify the accuracy and robustness. Adaptive refinement results are then presented to validate the proposed adaptive strategies. Moreover, parallel performance is illustrated to show the efficiency of our algorithm. Finally, the method is applied to the Laplace eigenvalue problem on a cubic domain. The resulting eigenvalue and eigenfunction approximations confirm the accuracy of the method and its applicability to elliptic eigenvalue problems.
All computations were carried out on a workstation with an AMD EPYC 9534 64-core processor and 1 TB of RAM using double-precision arithmetic, and the source code is publicly available at \url{https://github.com/zeyu-math/h-ATSEM}.}

\subsubsection{\zw{Poisson Equation}}

\paragraph{Cubic Domain:}

We consider a Poisson equation with Dirichlet boundary on $\Omega = (-1, 1)^{3}$, which is defined as

\begin{equation}
    \begin{cases}
        \displaystyle - \Delta u(\mathbf{x}) = f(\mathbf{x}), & \mathbf{x} \in \Omega,          \\
        \displaystyle u(\mathbf{x}) = u_{0} (\mathbf{x}),     & \mathbf{x} \in \partial \Omega.
    \end{cases}
\end{equation}

\begin{figure}[H]
    \centering
    \subfloat[Error with uniform refinement for fixed $p$.]{\includegraphics[width=0.45\textwidth]{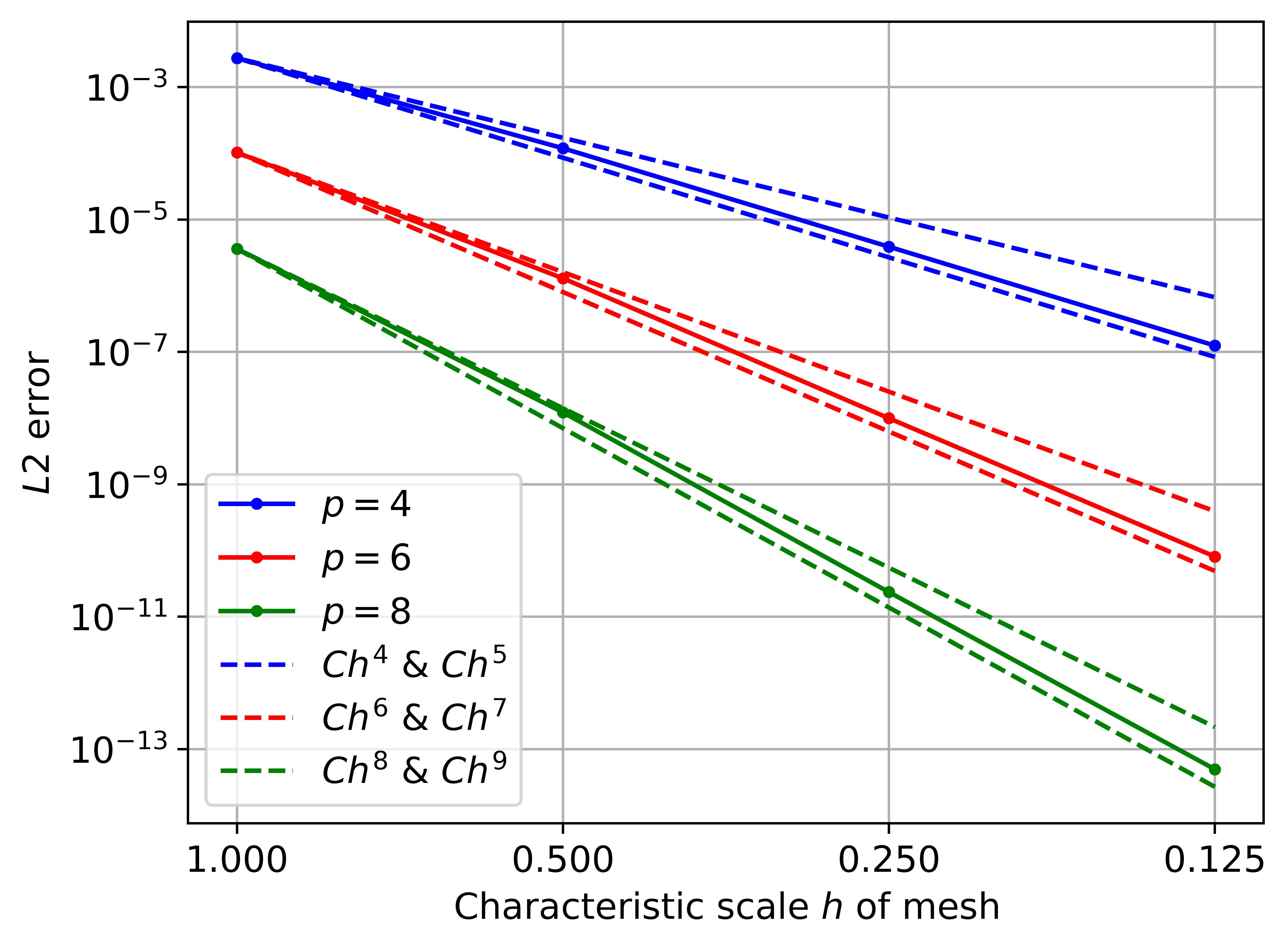}}
    \subfloat[Error with different $p$ on the same mesh.]{\includegraphics[width=0.45\textwidth]{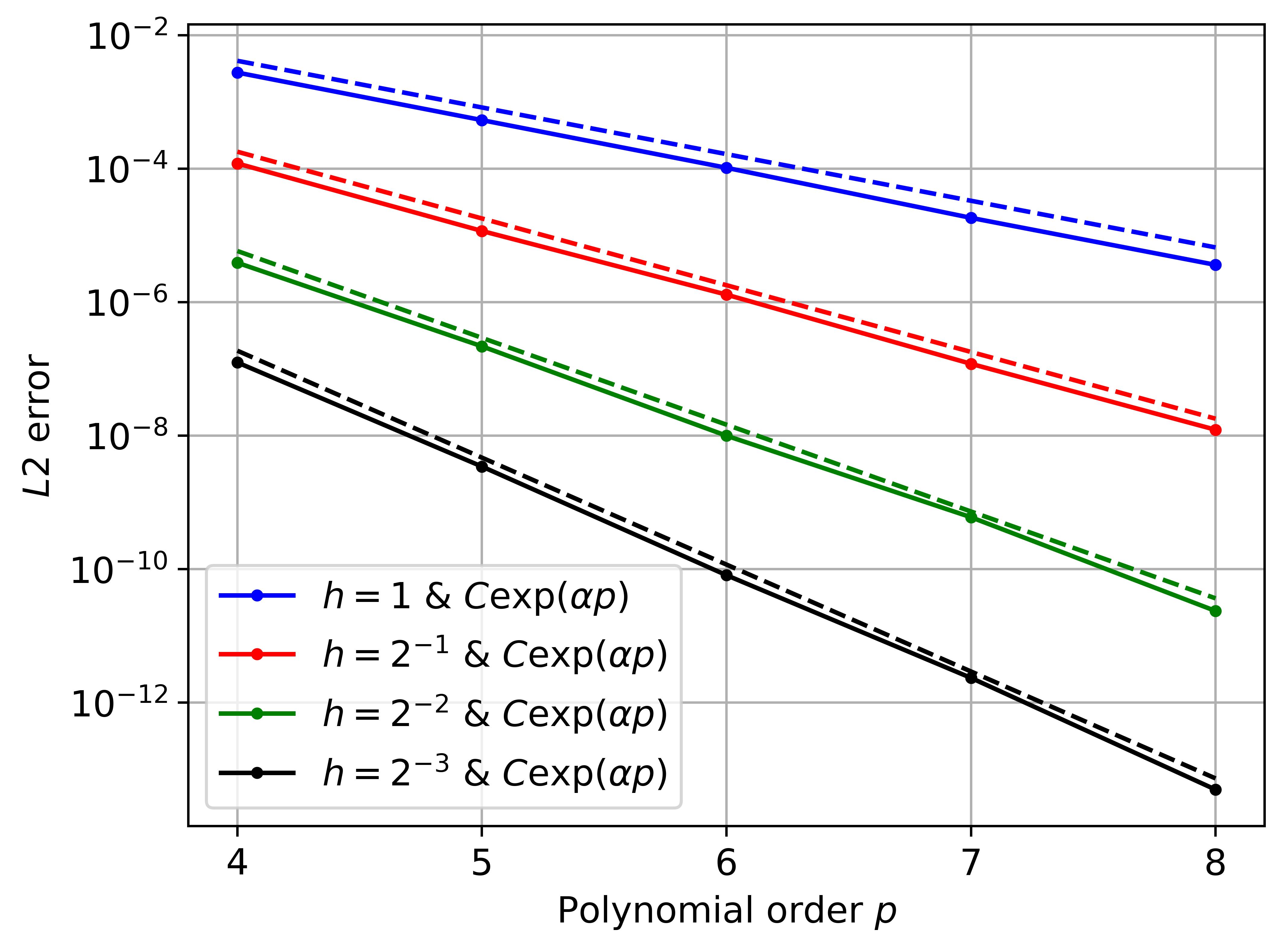}}
    \caption{The $L_{2}$ error for $u_{1} = \exp(-(x^{2} + y^{2} + z^{2}))$, where (a) shows the dependence of $L_{2}$ error on $h$, and (b) shows the relation between $L_{2}$ error and order $p$ with same mesh.}
    \label{figure::ConvergenceOrder}
\end{figure}

We first examine the convergence order with a smooth solution. Figure \ref{figure::ConvergenceOrder} shows the anticipated spectral accuracy of $L_{2}$ error for different polynomial degrees $p$ and edge lengths $h$, with solution

\begin{equation}
    u_{1} = \exp(-(x^{2} + y^{2} + z^{2})),
\end{equation}

\noindent and uniform refinements. The solid lines give the numerical results, and the dashed lines give the references. When fixing $p$ and applying the uniform refinement, $h$ is divided by $2$ each time, and the $L_{2}$ error is bounded by $C h^{p}$ and $C h^{p + 1}$, where $C$ is a constant related to $p$, as expected. Then, for different $p$ on the same mesh, the $L_{2}$ error follows the curve $C \exp(p)$, where $C$ is a constant related to the characteristic scale $h$.

\begin{figure}[H]
    \centering
    \subfloat[The DoF-error plot for solving $u_{2}$.]{\includegraphics[width=0.45\textwidth]{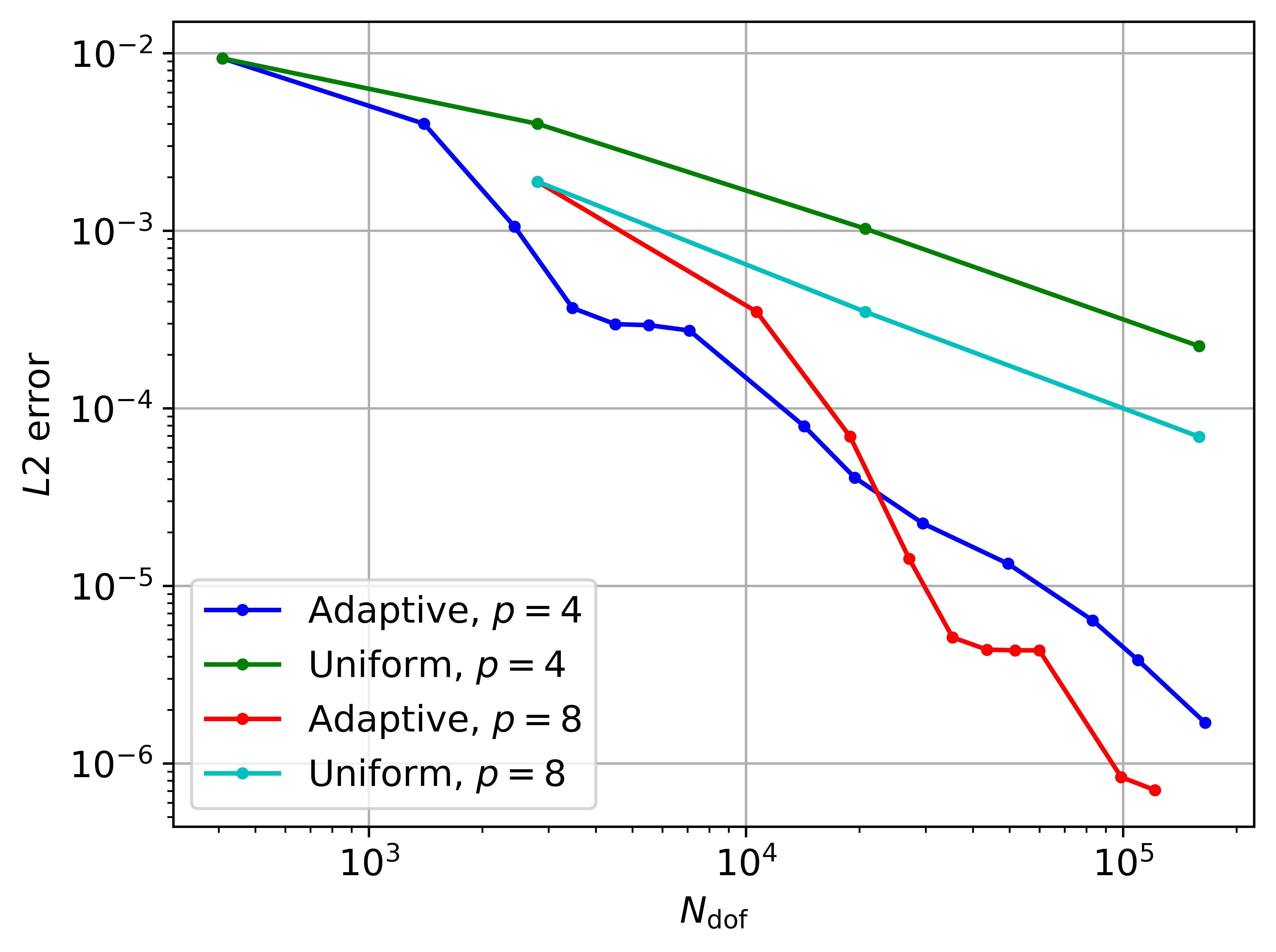}}
    \subfloat[The DoF-error plot for solving $u_{3}$.]{\includegraphics[width=0.45\textwidth]{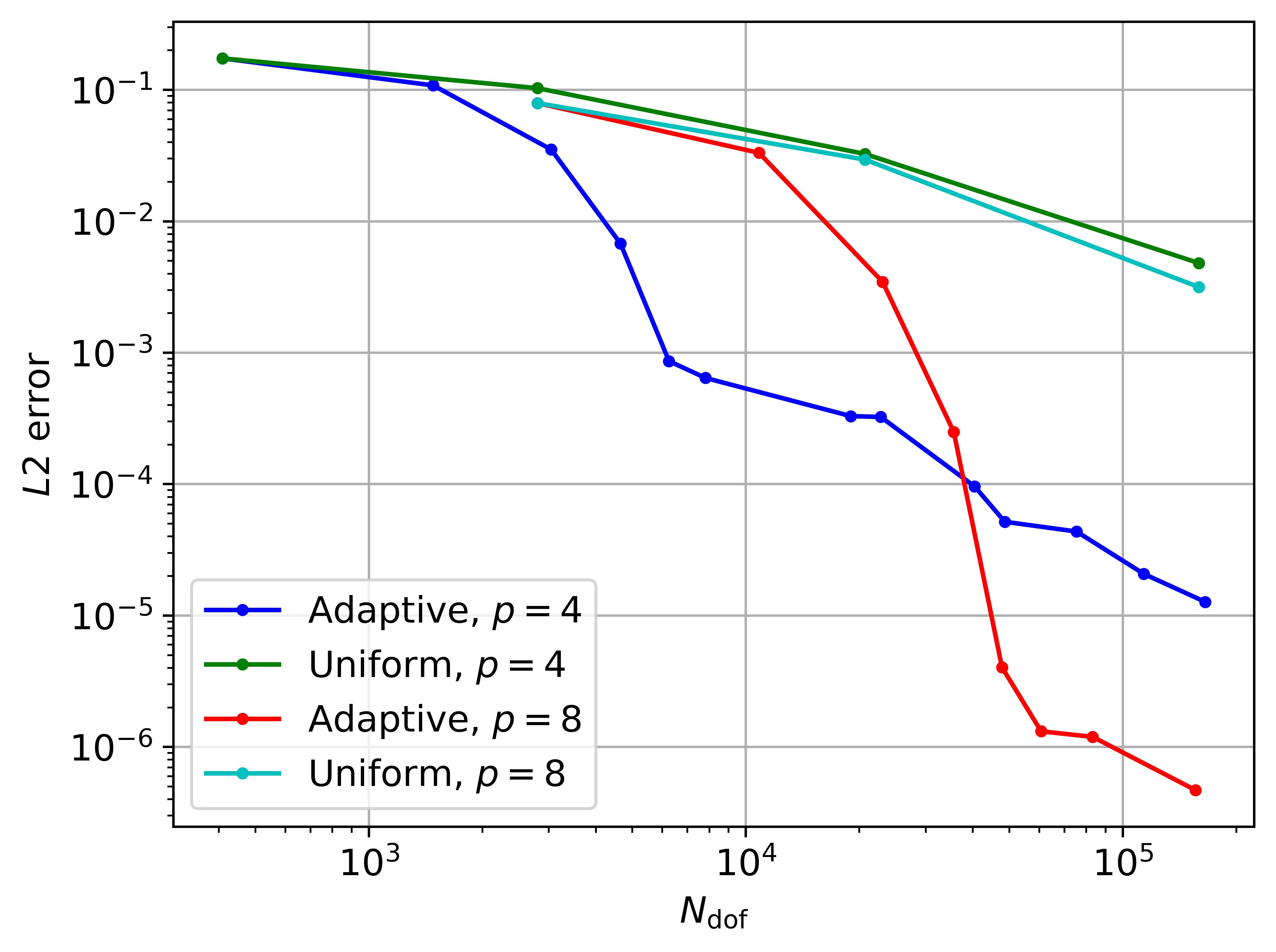}}
    \caption{The DoF-error plots for solving $u_{2}$ and $u_{3}$.}
    \label{figure::PoissonDofError}
\end{figure}

\begin{figure}[H]
    \centering
    \subfloat[The final mesh for solving $u_{2}$ with $p = 4$]{\includegraphics[width=0.225\textwidth]{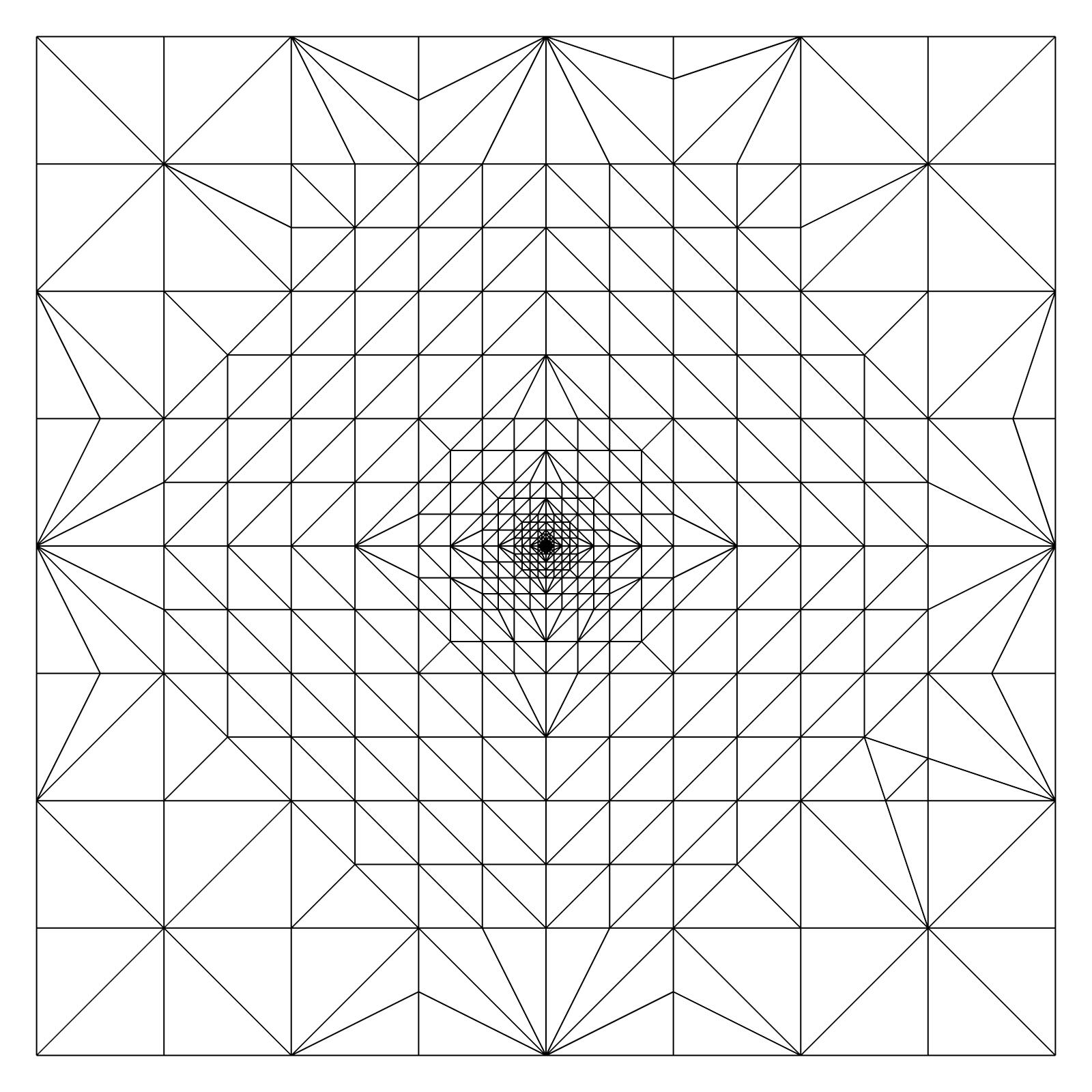}}
    \subfloat[The final mesh for solving $u_{2}$ with $p = 8$]{\includegraphics[width=0.225\textwidth]{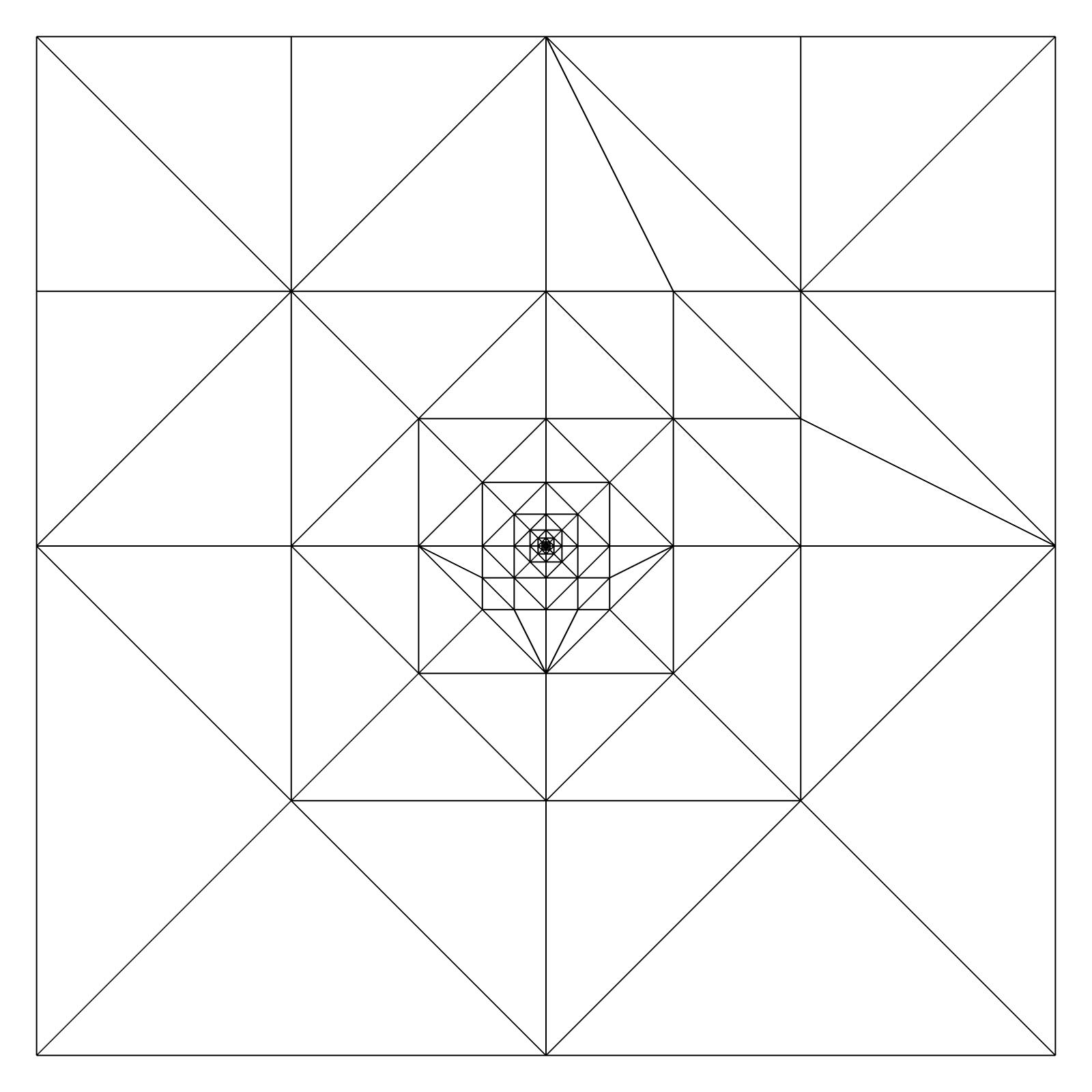}}
    \subfloat[The final mesh for solving $u_{3}$ with $p = 4$]{\includegraphics[width=0.225\textwidth]{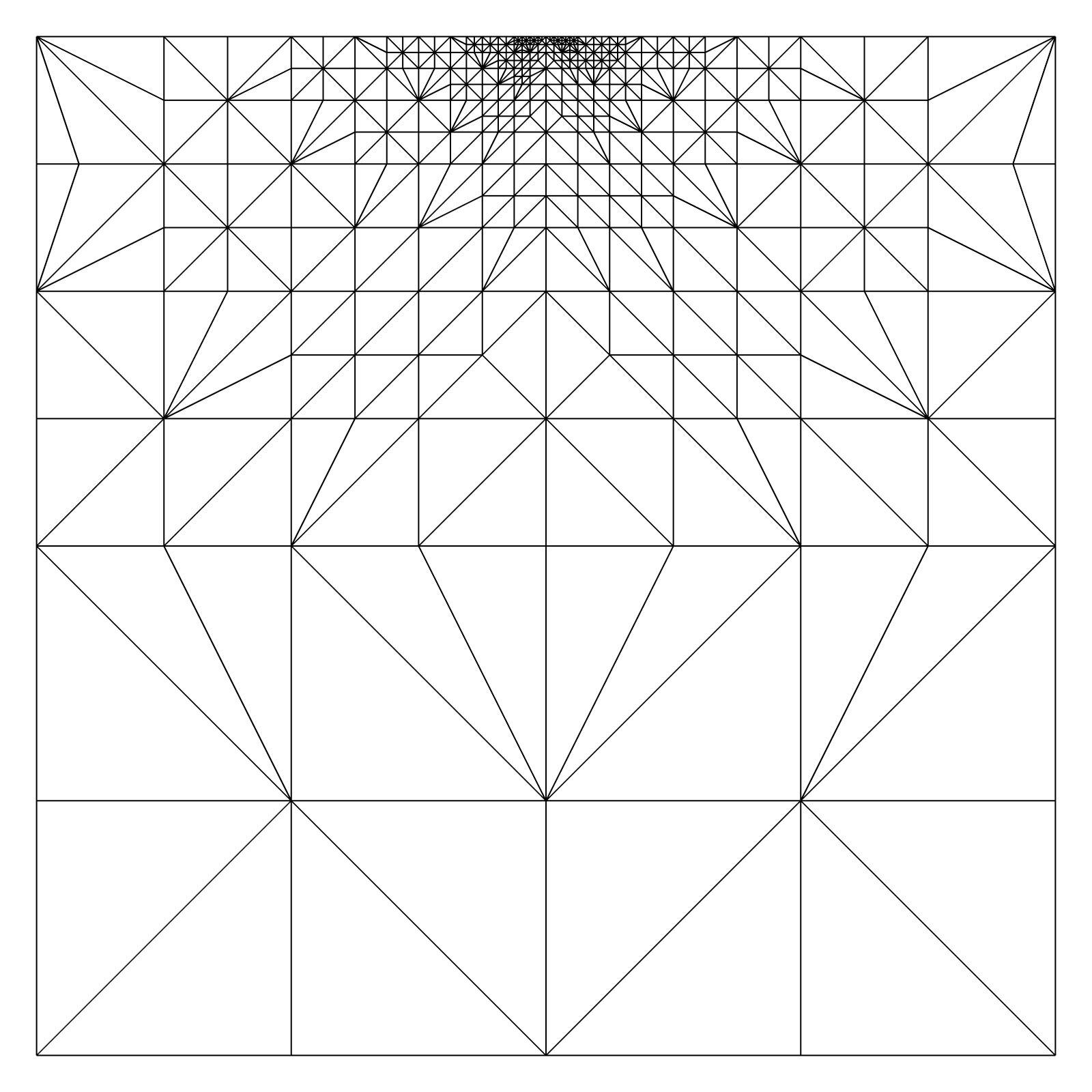}}
    \subfloat[The final mesh for solving $u_{3}$ with $p = 8$]{\includegraphics[width=0.225\textwidth]{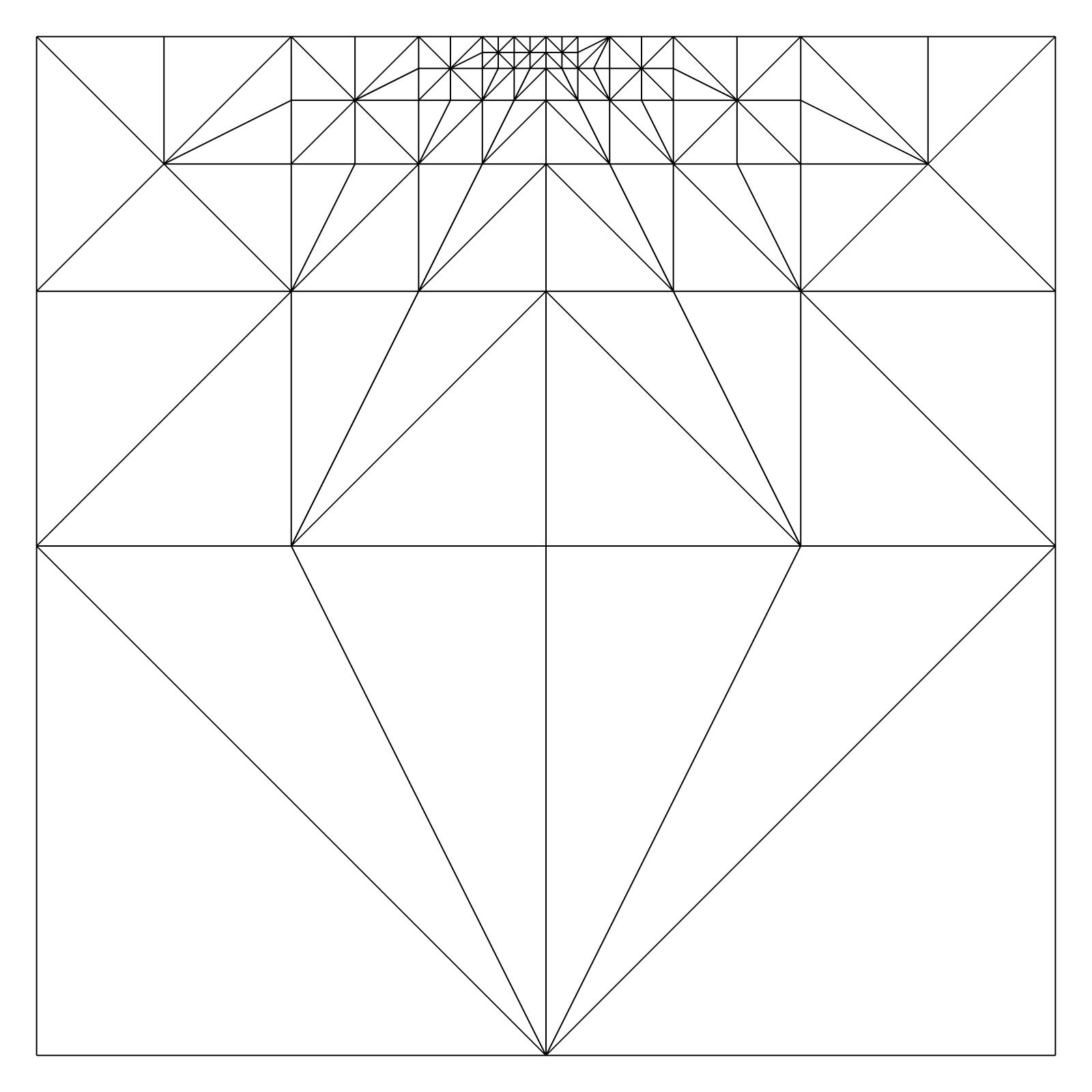}}
    \caption{The final meshes for solving $u_{2}$ and $u_{3}$, with $p = 4$ and $p = 8$.}
    \label{figure::PoissonMesh}
\end{figure}

Then, two irregular solutions were tested, one with a singularity and another with an oscillatory region, defined as

\begin{equation}
    u_{2} = \exp(-10 \sqrt{x^{2} + y^{2} + z^{2}}), \quad
    u_{3} = \sin\left(\frac{1}{\sqrt{x^{2} + y^{2} + (z - 1.1)^{2}}}\right),
\end{equation}

\noindent where $u_{2}$ includes a singularity at the origin point, and $u_{3}$ oscillates around $(0, 0, 1.1)$. As shown in Figure \ref{figure::PoissonDofError} and Figure \ref{figure::PoissonMesh}, the $h$-adaptive method refines the mesh around the singularity and the oscillatory region, achieving a higher accuracy with the same number of degrees of freedom, compared to uniform refinement. On the other hand, the $h$-adaptive method also benefits from higher-order basis polynomial functions. After sufficient iterations, the higher-order cases perform better, as a relatively coarse mesh is enough for higher-order basis functions to achieve the same accuracy.

The program was tested for core counts from $1$ to $64$, doubling at each step. As shown in Table \ref{table::TimeCostPoisson}, $p=4$ and $p=8$ were tested with the same initial mesh and the same quadrature formula. When reaching the same level of $L_{2}$ error ($8.27 \times 10^{-8}$ for $p = 4$, and $4.3 \times 10^{-8}$ for $p = 8$), the lower-order case takes $1209245$ DoFs, while high-order case takes only $296285$ DoFs due to the higher convergence rate.

\begin{table}[H]
    \centering
    \caption{Computational time (second) for solving Poisson equation with different numbers of CPU cores. The programs stops as the error is less than $10^{-7}$. For $p = 4$ and $p = 8$, the final numbers of degrees of freedom are $1209245$ and $296285$, respectively.}
    \label{table::TimeCostPoisson}
    \begin{tabular}{|c|c|c|c|c|c|c|c|c|}
        \hline
        \multirow{2}{*}{Core} & \multicolumn{4}{c|}{$p = 4, \text{final } N_{\text{dof}} = 1209245$} & \multicolumn{4}{c|}{$p = 8, \text{final } N_{\text{dof}} = 296285$} \\
        \cline{2-9}
        & PCG & Adaptivity & Assembly & Total & PCG & Adaptivity & Assembly & Total \\
        \hline
        $1$  & $83.45$ & $617.09$ & $131.02$ & $832.30$ & $219.63$ & $148.75$ & $116.27$ & $485.27$ \\
        $2$  & $50.97$ & $377.75$ & $73.56$  & $503.23$ & $132.60$ & $81.77$  & $62.20$  & $277.26$ \\
        $4$  & $38.05$ & $215.79$ & $39.85$  & $294.41$ & $118.01$ & $45.83$  & $33.93$  & $198.57$ \\
        $8$  & $38.15$ & $131.70$ & $22.04$  & $192.24$ & $117.75$ & $27.27$  & $19.49$  & $165.20$ \\
        $16$ & $29.81$ & $96.93$  & $13.80$  & $140.98$ & $66.57$  & $17.17$  & $11.86$  & $96.24$ \\
        $32$ & $22.86$ & $72.41$  & $9.45$   & $105.05$ & $45.20$  & $12.59$  & $8.48$   & $66.96$ \\
        $64$ & $20.38$ & $64.10$  & $8.05$   & $92.86$  & $39.06$  & $10.05$  & $7.05$   & $56.84$ \\
        \hline
    \end{tabular}
\end{table}

It is noticed that when stopping at an error less than $10^{-7}$, the case of $p=8$ takes considerably less time than the case of $p=4$. More specifically, when using a single CPU core, the overall time cost of the higher-order case is about $58\%$ of the lower-order case, and for multiple cores, it is on average $67\%$. This result is mainly because the higher-order case converges faster and thus requires fewer adaptive iterations. In particular, the computational cost of the adaptivity and assembly components is given by 

\begin{equation}
    \begin{aligned}
        & T^{p=4}_{\text{adapt}} = 35 C_{\text{adapt}} N_{\text{tet}}, & & T^{p=4}_{\text{assem}} = 35^{2} C_{\text{assem}} N_{\text{tet}}, \\
        & T^{p=8}_{\text{adapt}} = 165 C_{\text{adapt}} N_{\text{tet}}, & & T^{p=8}_{\text{assem}} = 165^{2} C_{\text{assem}} N_{\text{tet}}, \\
    \end{aligned}
\end{equation}

\noindent where $N_{\text{tet}}$ is the number of tetrahedra, $C_{\text{adapt}}$ and $C_{\text{assem}}$ are constants, and $35$ and $165$ are the numbers of basis functions for $p=4$ and $p=8$, respectively. Therefore, with the same mesh, the computational costs for $p=8$ are about $4.7$ and $22.2$ times those for $p = 4$ for adaptivity and assembly, respectively. However, the final mesh contains $113008$ tetrahedra for $p = 4$, approximately $33$ times as many as the $3448$ tetrahedra for $p = 8$. Consequently, it can be observed that the computational time for $p=8$ is significantly lower than that for $p=4$ in the adaptive refinement and assembly stages. In contrast, the opposite trend is observed in the PCG stage, as a much denser linear system is generated for higher-order cases.

We further examined the performance for different core counts. According to Amdahl's law, if $\alpha_{p}$ denotes the parallelizable proportion of the program, the maximum speedup on $N$ cores is
\begin{equation}
    S(N, \alpha_{p}) = \frac{1}{1 - \alpha_{p} + \frac{\alpha_{p}}{N}}.
\end{equation}
As shown in Figure \ref{table::SpeedupPoisson}, we set $\alpha_{p} = 0.9$, the speedups for total time align well with theoretical predictions for $32$ and $64$ cores. However, a performance plateau is found from $4$ to $16$ cores, pointing to certain constraints that cause this behavior.

\begin{table}[H]
    \centering
    \caption{The speedup for different cores, where for theoretical predictions from Amdahl's Law, we set $\alpha_{p} = 0.9$.}
    \label{table::SpeedupPoisson}
    \begin{tabular}{|c|c|c|c|c|c|c|c|c|}
        \hline
        $N$ & $1$ & $2$ & $4$ & $8$ & $16$ & $32$ & $64$ \\
        \hline
        $p = 4$ & $1.00$ & $1.65$ & $2.83$ & $4.33$ & $5.90$ & $7.92$ & $8.96$ \\
        $p = 8$ & $1.00$ & $1.75$ & $2.44$ & $2.94$ & $5.04$ & $7.25$ & $8.54$ \\
        $S(N, 0.9)$ & $1.00$ & $1.82$ & $3.08$ & $4.71$ & $6.40$ & $7.80$ & $8.77$ \\
        \hline
    \end{tabular}
\end{table}

\begin{figure}[H]
    \centering
    \subfloat[Speedup for $p = 4$.]{\includegraphics[width=0.45\textwidth]{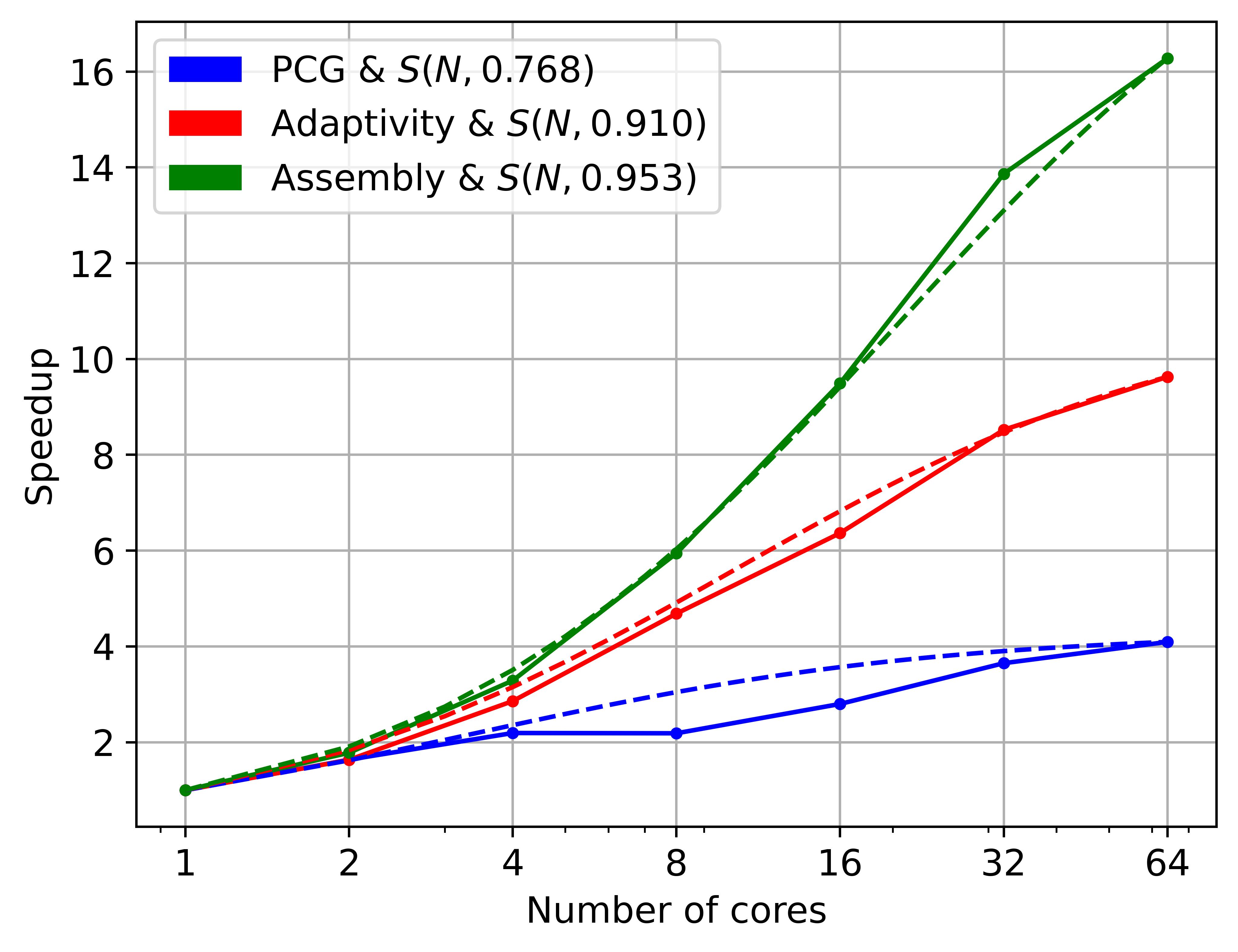}}
    \subfloat[Speedup for $p = 8$.]{\includegraphics[width=0.45\textwidth]{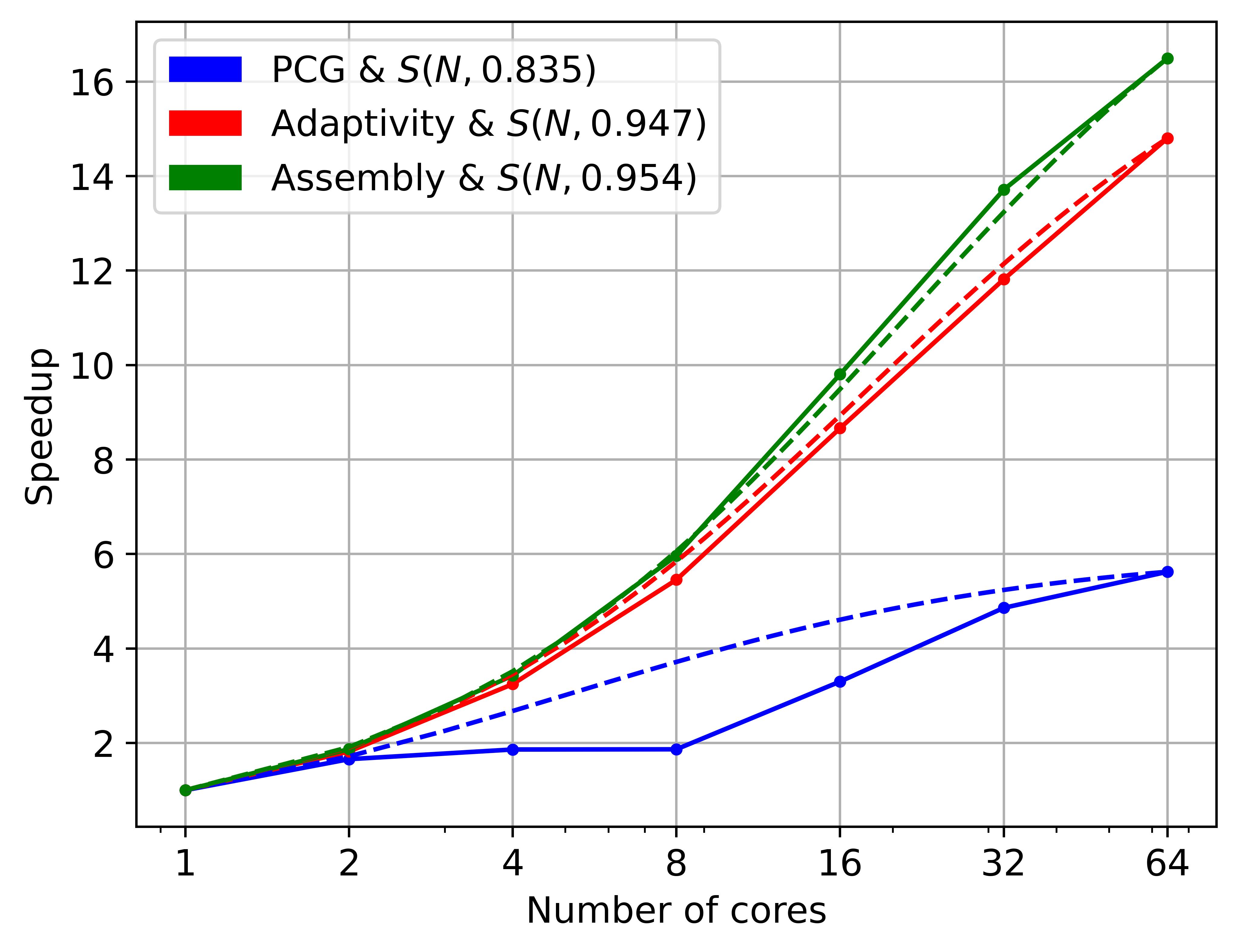}}
    \caption{The speedup of each module for solving Poisson equation with $p = 4$ and $p = 8$, where $S(N, \alpha)$ is the theoretical predictions from Amdahl's law.}
    \label{figure::TimeCostPoisson}
\end{figure}

Therefore, we further examined the performance of each module. As shown in Figure \ref{figure::TimeCostPoisson}, we fit the curve according to Amdahl's Law, and the speedups for adaptivity and assembly align well with theoretical predictions. However, the PCG solver exhibits a performance plateau between $4$ and $16$ cores, consistent with the behavior observed above. 
The observed performance is characteristic of the architecture of the AMD EPYC 9534 64-core processor \cite{amd_4th_gen_epyc_2024}. Specifically, this CPU consists of $8$ Core Complex Dies (CCDs), each housing $8$ cores and a shared $32$ MB L3 cache.
Therefore, for thread counts up to $8$, the workload is likely confined within a single CCD. In these cases, the observed plateau indicates a memory-bound regime where the memory bandwidth of a single CCD reaches saturation.
Then, as the thread count exceeds $16$, the speedup recovers as the workload spreads across multiple CCDs, tapping into higher aggregate memory bandwidth.
Moreover, as the per-core working set shrinks under strong scaling, it eventually becomes cache-resident within the expanded aggregate L3 cache, further accelerating the solver by transitioning from memory-bound to cache-bound execution.



Based on the above analysis, when solving the Poisson equation, the higher-order method outperforms the lower-order one in both total execution time and parallel performance. An analysis based on Amdahl's law indicates that up to 95\% of the adaptation and assembly steps can be parallelized, which is remarkably high given the inherently irregular data structures encountered in 3D tetrahedral $h$-adaptivity. Based on this parallelization profile, the computational speedup is projected to scale effectively with increasing core counts, gradually approaching the theoretical asymptotic limit of a 20-fold improvement. We also note that the PCG solver shows a relatively poor parallel speedup in this problem, which will be a focus of future optimization.

\paragraph{L-Shaped Domain:}

\begin{figure}[H]
    \centering
    \subfloat[Error with uniform refinement for fixed $p$.]{\includegraphics[width=0.45\textwidth]{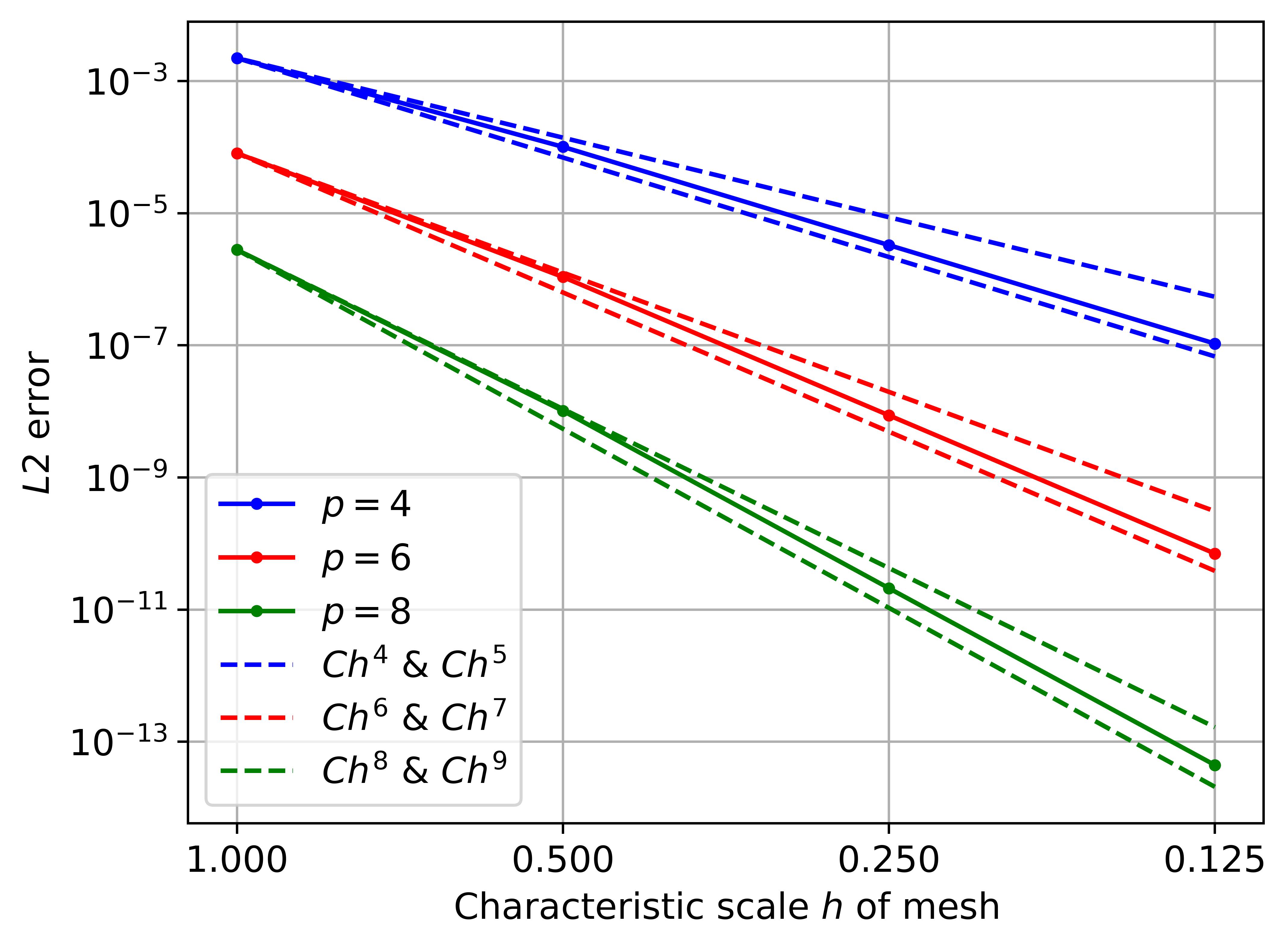}}
    \subfloat[Error with different $p$ on the same mesh.]{\includegraphics[width=0.45\textwidth]{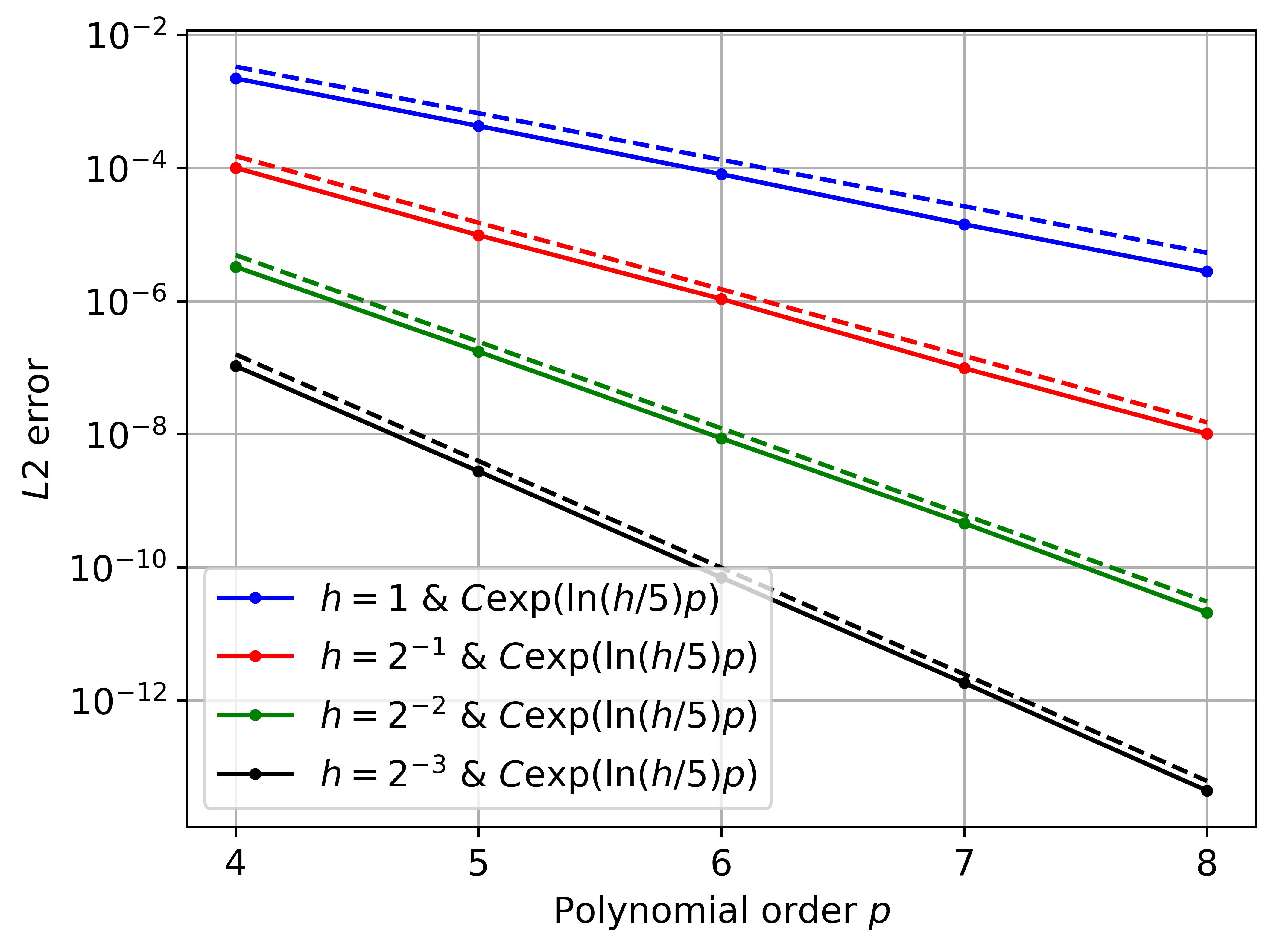}}
    \caption{The $L_{2}$ error for $u_{1} = \exp(-(x^{2} + y^{2} + z^{2}))$, where (a) shows the dependence of $L_{2}$ error on $h$, and (b) shows the relation between $L_{2}$ error and order $p$ with same mesh.}
    \label{figure::ConvergenceOrderL}
\end{figure}

\zw{The Poisson solver is then tested on an L-shaped prismatic domain $\Omega = ((-1, 1)^{2} \setminus (0, 1)^{2}) \times (-1, 1)$. We first consider the solution $u_{1} = \exp(-(x^{2} + y^{2} + z^{2}))$. The convergence behaviors are presented in Figure \ref{figure::ConvergenceOrderL}, where the solid lines represent the numerical results and the dashed lines indicate the reference convergence rates. It can be observed that, despite the non-convex geometry of the L-shaped domain, the proposed solver demonstrates the same spectral convergence behavior in the $L_{2}$ norm as that obtained on the cubic domain.}

\begin{figure}[H]
    \centering
    \subfloat[The DoF-error plot for solving $u_{2}$.]{\includegraphics[width=0.45\textwidth]{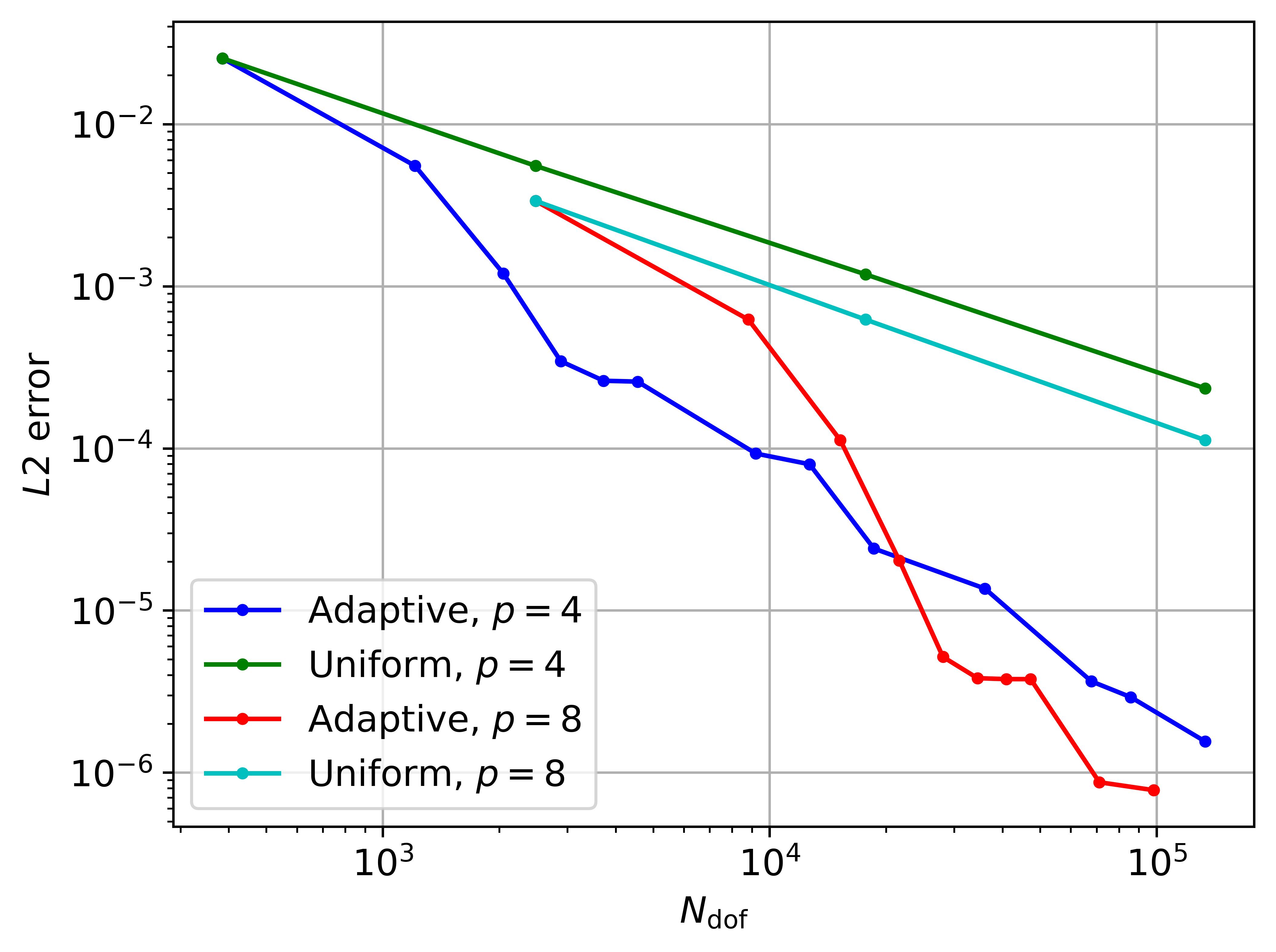}}
    \subfloat[The DoF-error plot for solving $u_{3}$.]{\includegraphics[width=0.45\textwidth]{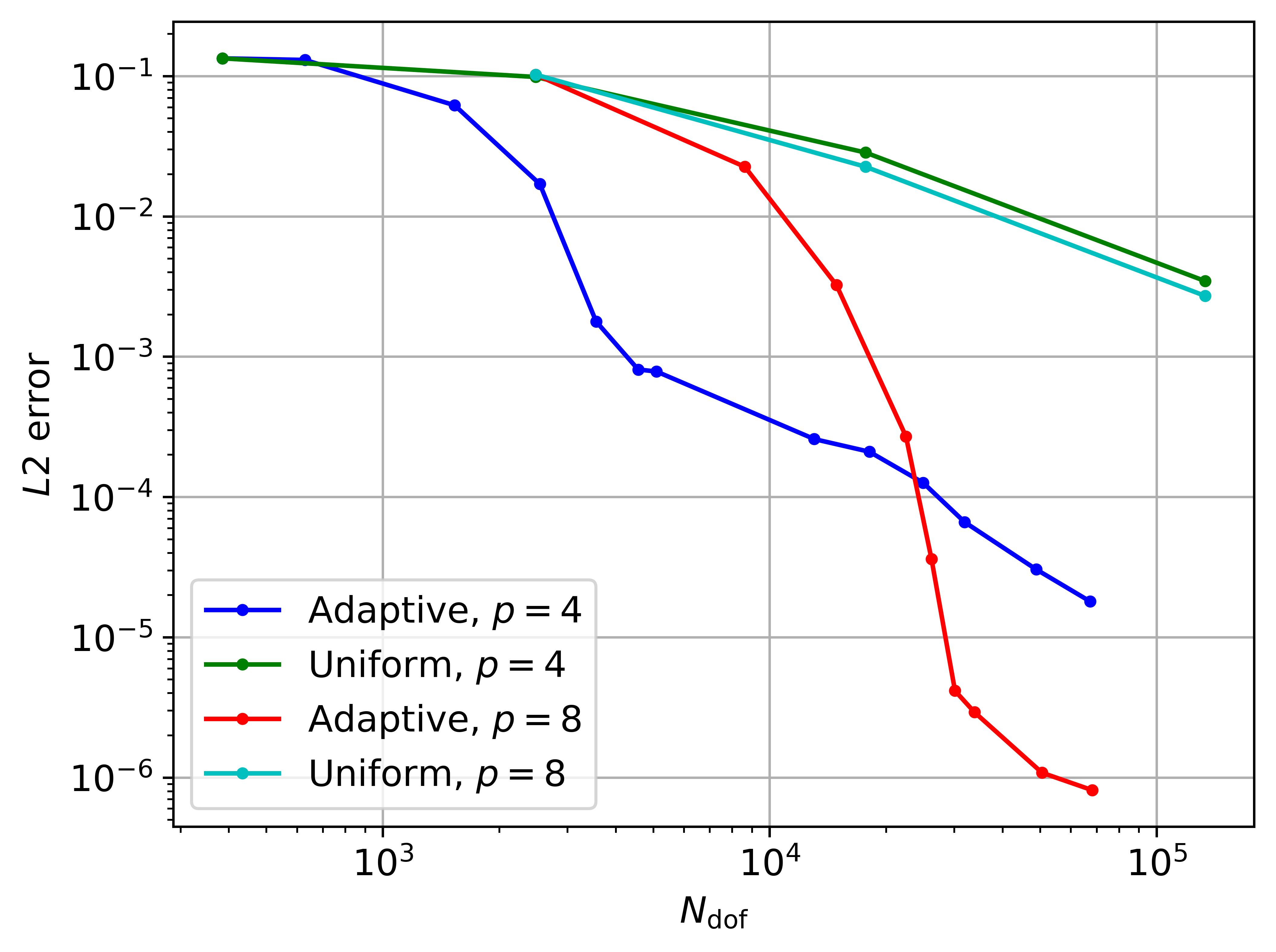}}
    \caption{The DoF-error plots for solving $u_{2}$ and $u_{3}$.}
    \label{figure::PoissonDofErrorL}
\end{figure}

\begin{figure}[H]
    \centering
    \subfloat[The final mesh for solving $u_{2}$ with $p = 4$]{\includegraphics[width=0.225\textwidth]{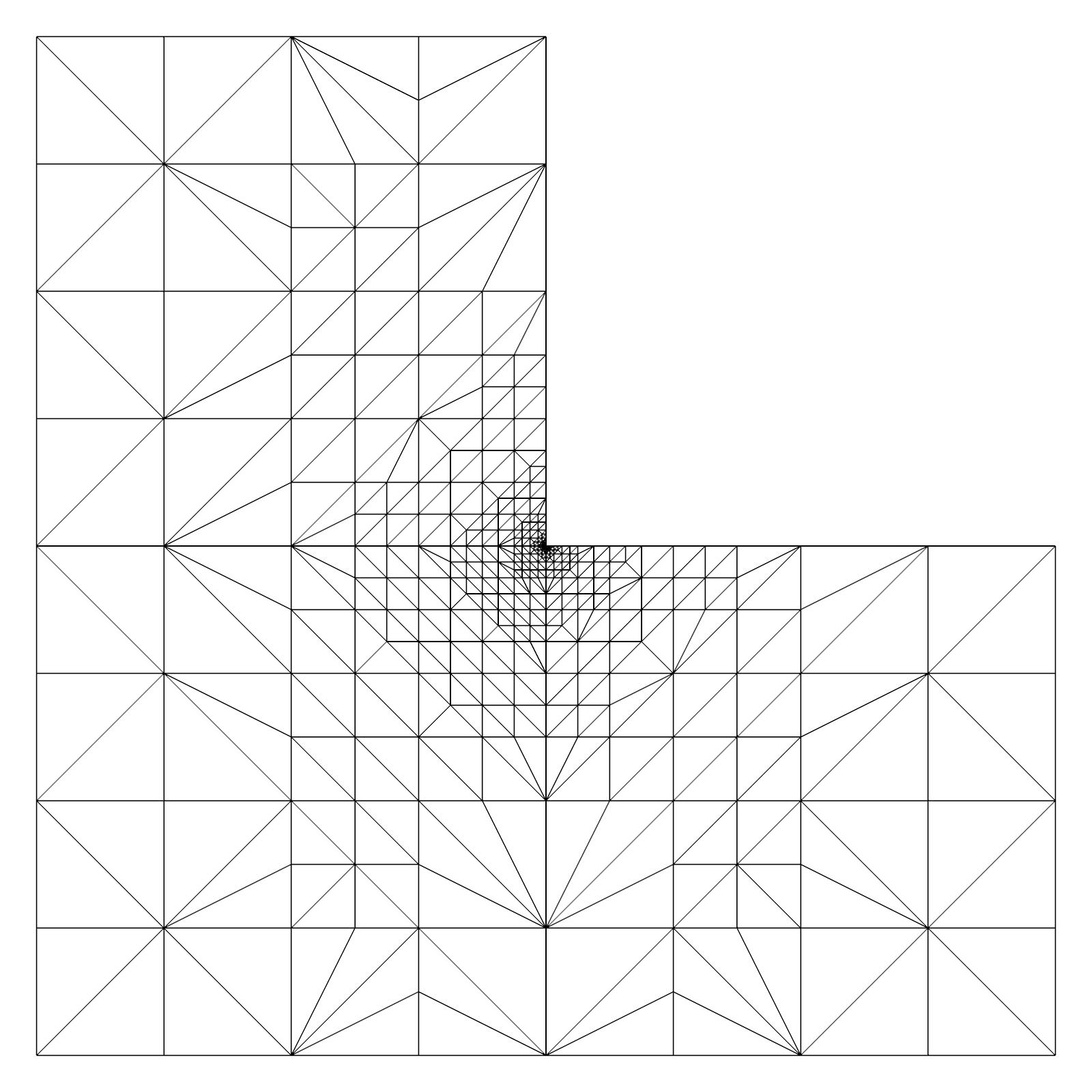}}
    \subfloat[The final mesh for solving $u_{2}$ with $p = 8$]{\includegraphics[width=0.225\textwidth]{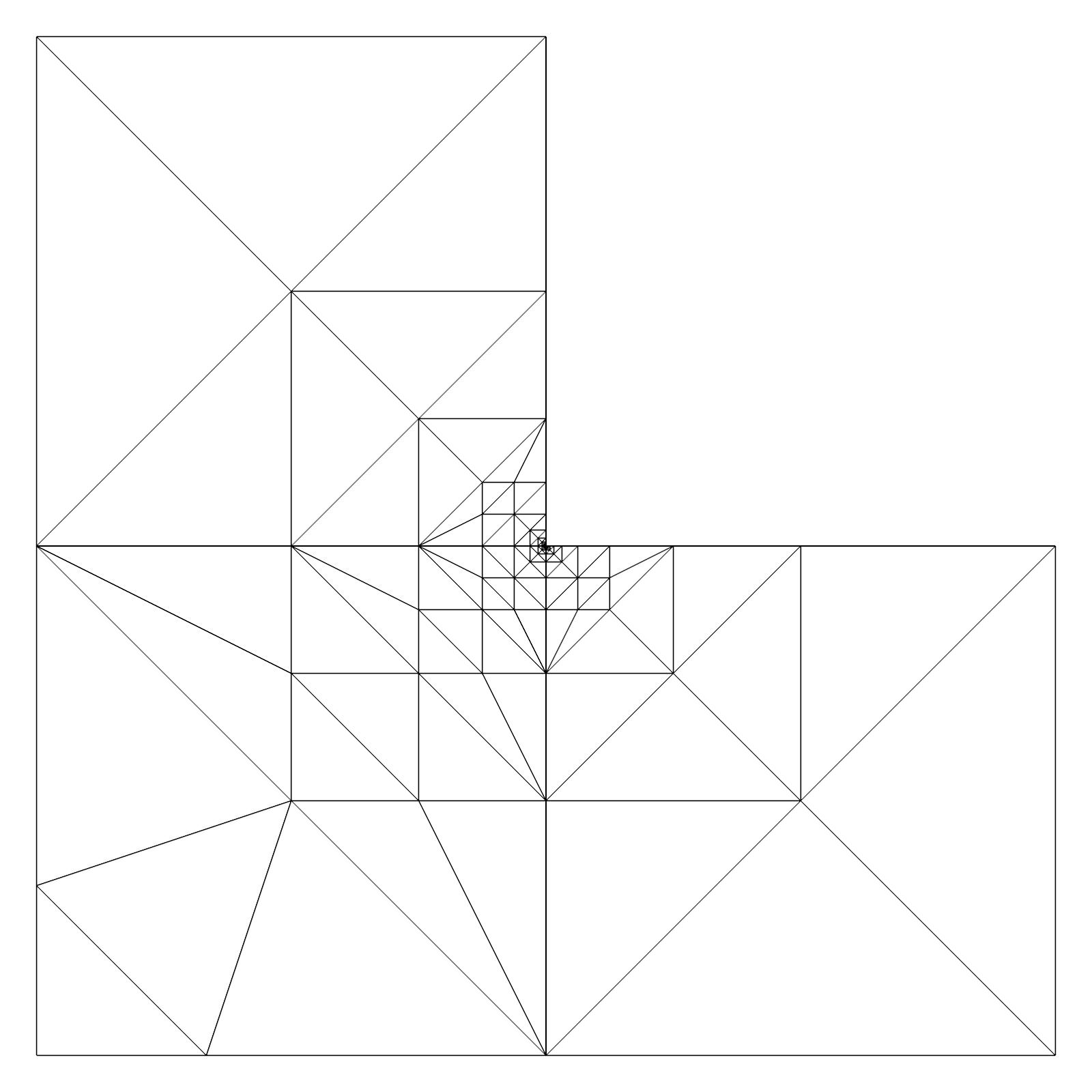}}
    \subfloat[The final mesh for solving $u_{3}$ with $p = 4$]{\includegraphics[width=0.225\textwidth]{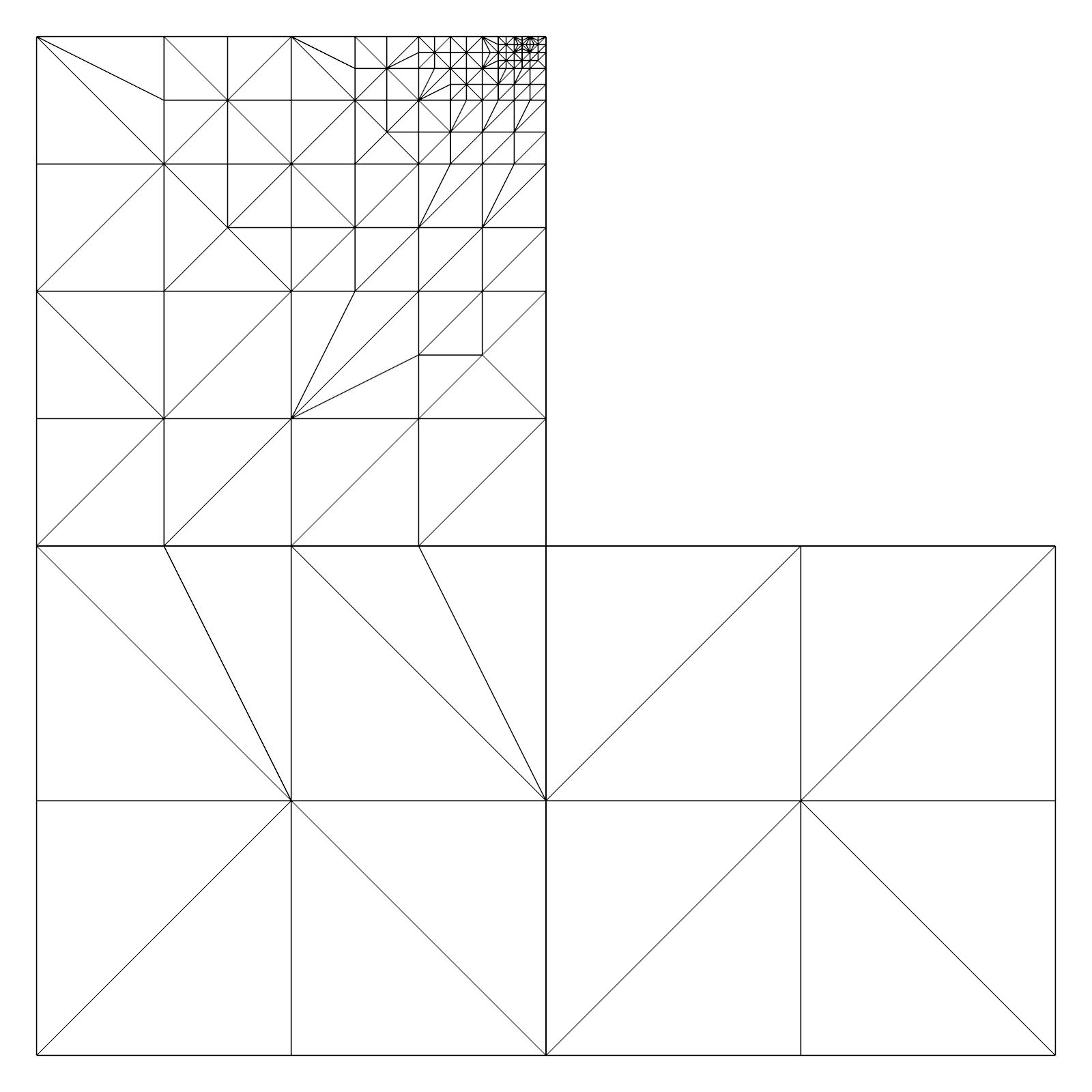}}
    \subfloat[The final mesh for solving $u_{3}$ with $p = 8$]{\includegraphics[width=0.225\textwidth]{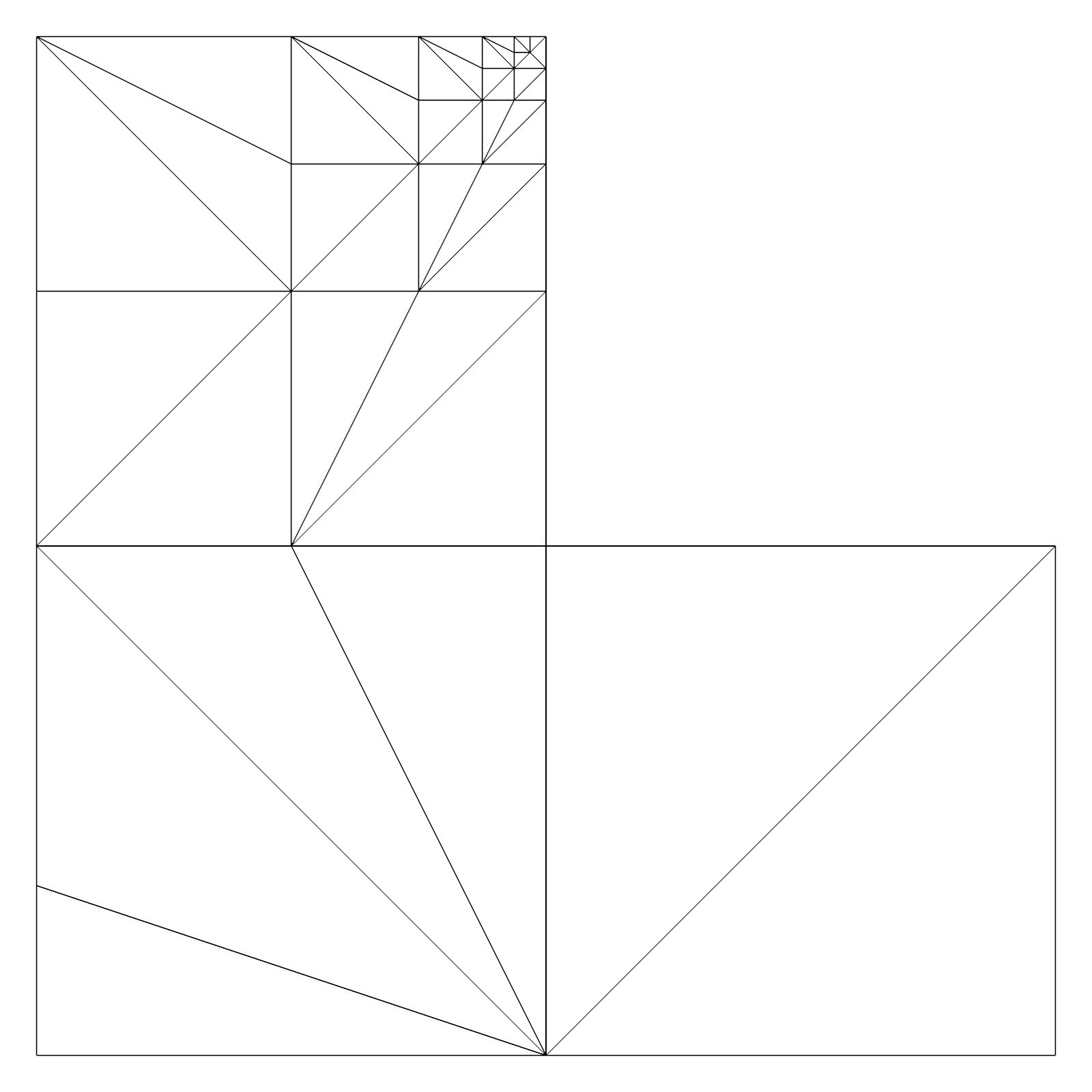}}
    \caption{The final meshes for solving $u_{2}$ and $u_{3}$, with $p = 4$ and $p = 8$.}
    \label{figure::PoissonMeshL}
\end{figure}

\zw{We also consider $u_{2} = \exp(-10 \sqrt{x^{2} + y^{2} + z^{2}})$ and $u_{3} = \sin\left(\frac{1}{\sqrt{x^{2} + y^{2} + (z - 1.1)^{2}}}\right)$, which contain localized variations and oscillatory features, respectively. The dof-error plots and refined meshes are shown in Figure \ref{figure::PoissonDofErrorL} and Figure \ref{figure::PoissonMeshL}, respectively. 
Similarly, the proposed method achieves satisfactory accuracy with relatively few degrees of freedom. The adaptively refined meshes effectively resolve the singularity in each solution. It is also worth noting that higher-order basis functions lead to faster convergence and require fewer mesh elements to attain the same accuracy.
}

\subsubsection{\zw{Laplacian Eigenvalue Problem}}

The Laplacian eigenvalue problem in $\Omega = (-1, 1)^{3}$ takes the form

\begin{equation}
    \begin{cases}
        \displaystyle - \Delta u(\mathbf{x}) = \lambda u(\mathbf{x}), & \mathbf{x} \in \Omega,          \\
        \displaystyle u(\mathbf{x}) = 0,                              & \mathbf{x} \in \partial \Omega,
    \end{cases}
\end{equation}

\noindent and the eigenvalues are 

\begin{equation}
    {\lambda_{i, j, k} = \frac{(i^{2} + j^{2} + k^{2}) \pi^{2}}{4}, \forall i, j, k \in \mathbb{N}^*}
\end{equation}

\begin{figure}[H]
    \centering
    \subfloat[$N_{\text{dof}} = 2031$]{\includegraphics[width=0.45\textwidth]{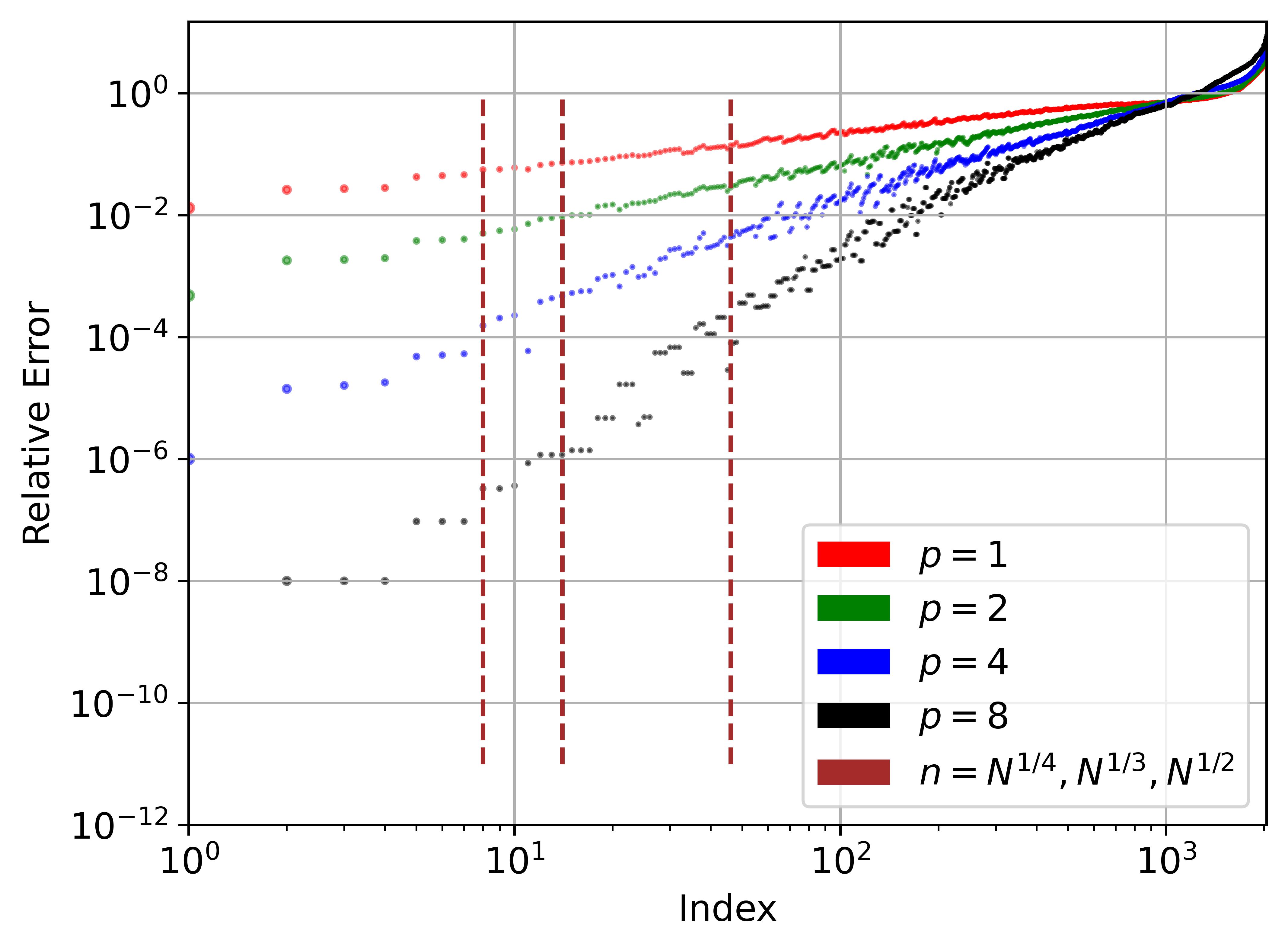}}
    \subfloat[$N_{\text{dof}} = 17631$]{\includegraphics[width=0.45\textwidth]{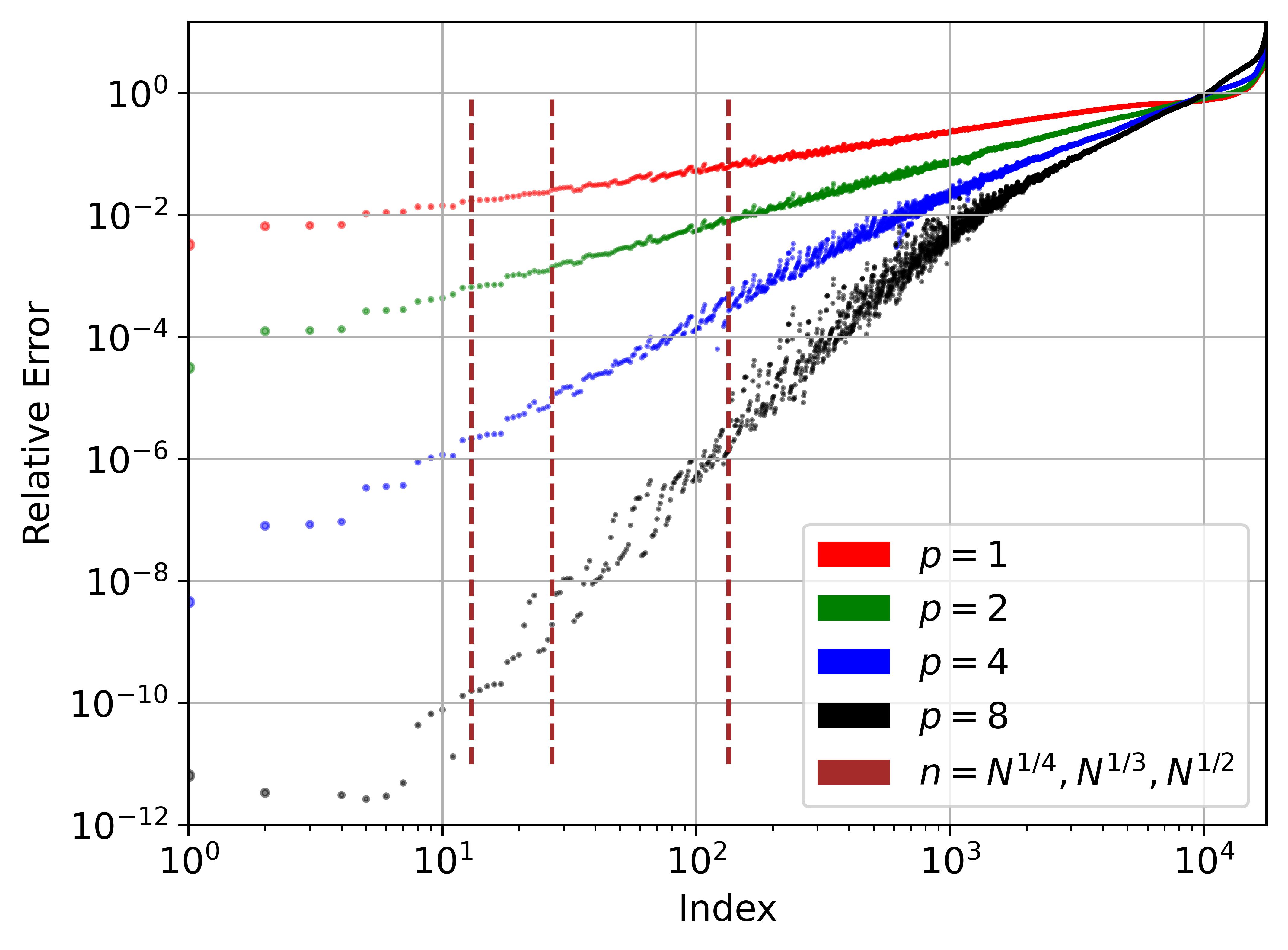}}
    \caption{Relative errors of the eigenvalues with $p = 1, 2, 4, 8$ and $N_{\text{dof}} = 2031, 17631$. For $N_{\text{dof}} = 2031$ and $p = 8$, the relative error of the first eigenvalue is small than $10^{-12}$, which is considered a coincidence.}
    \label{figure::Laplace}
\end{figure}

\begin{figure}[H]
    \centering
    \subfloat[$N_{\text{dof}} = 2031$]{\includegraphics[width=0.45\textwidth]{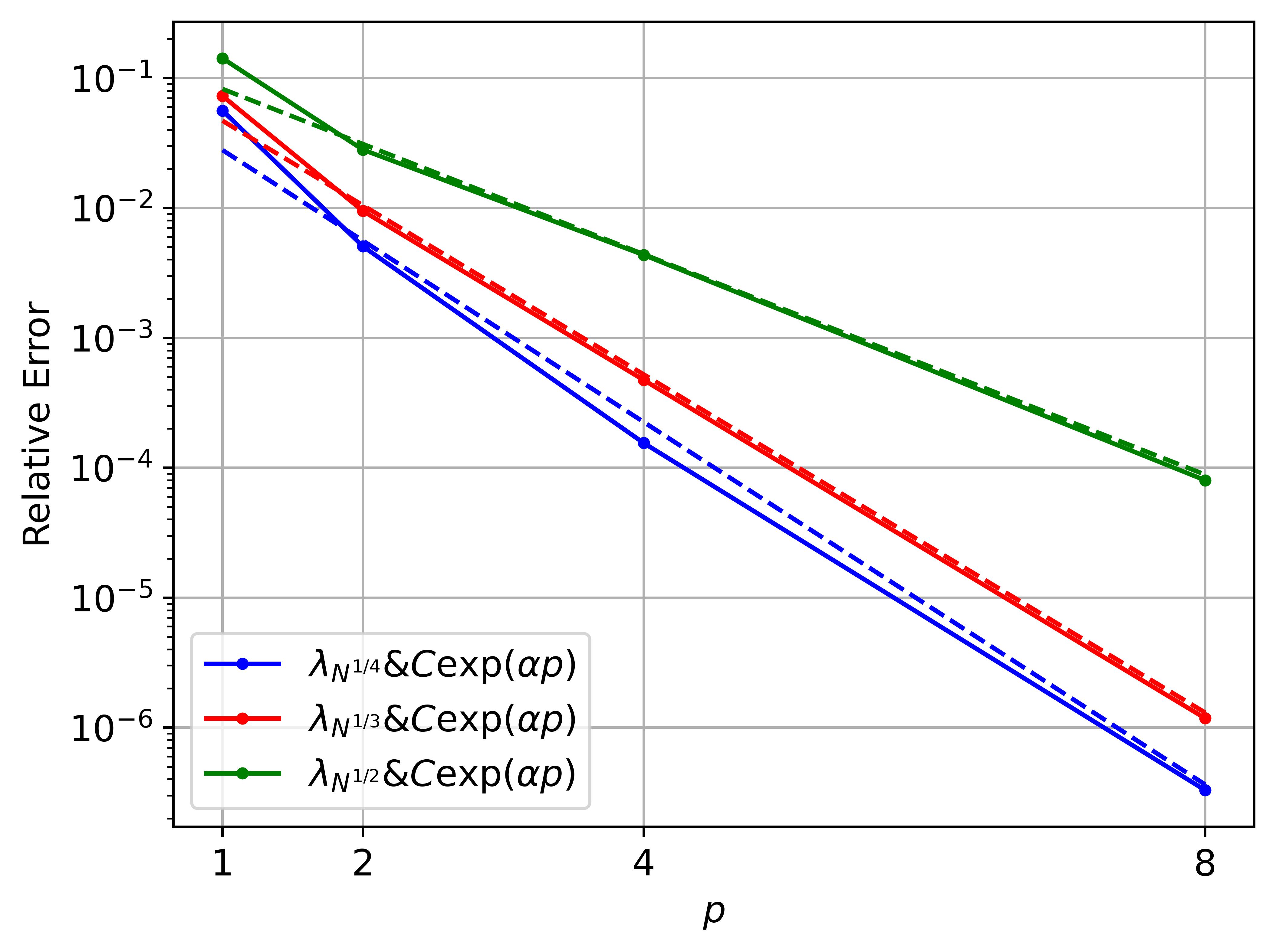}}
    \subfloat[$N_{\text{dof}} = 17631$]{\includegraphics[width=0.45\textwidth]{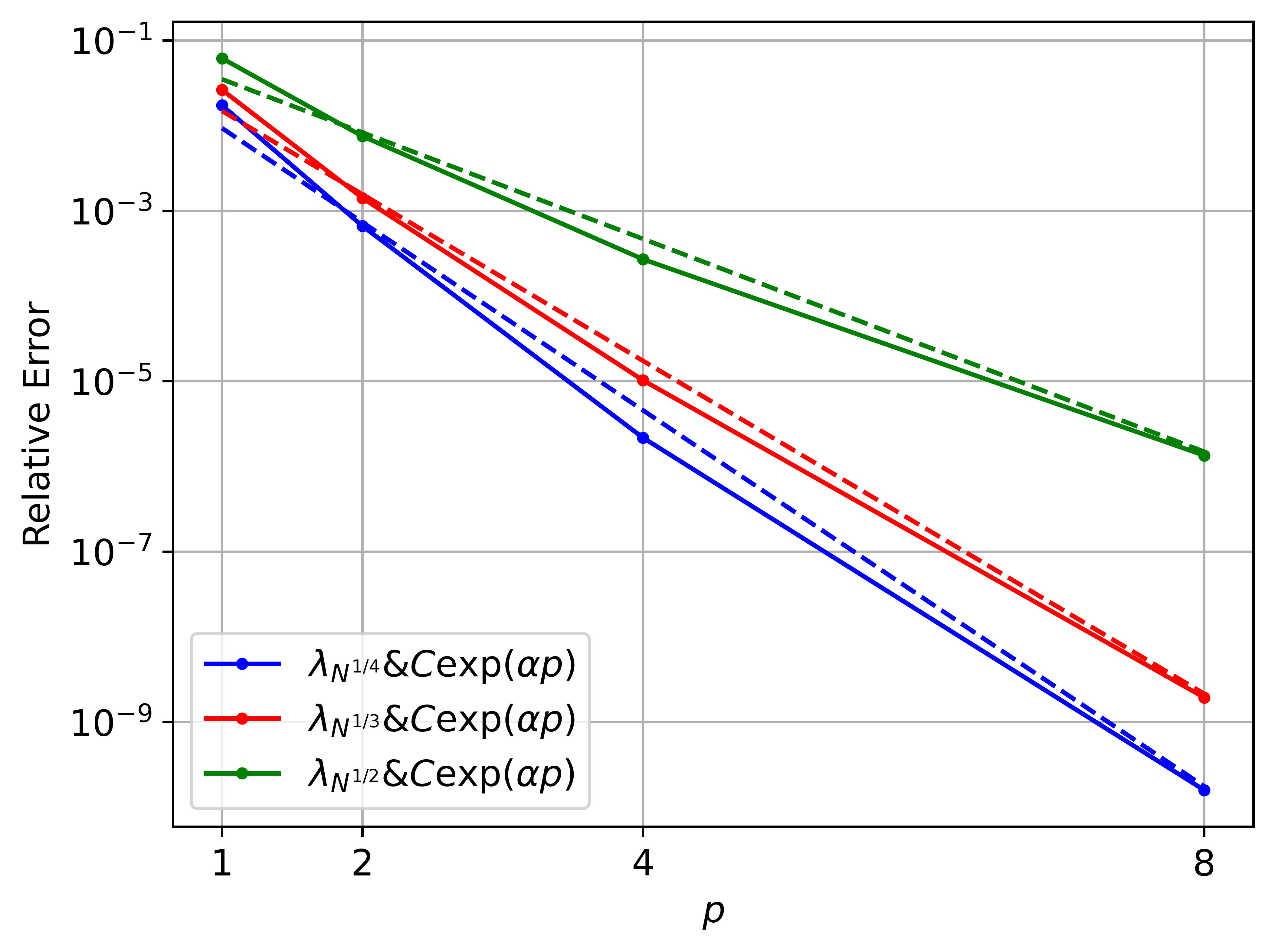}}
    \caption{Relative errors of the eigenvalues at positions $N^{\frac{1}{4}}$, $N^{\frac{1}{3}}$ and $N^{\frac{1}{2}}$ with different $p$.}
    \label{figure::LaplaceRate}
\end{figure}

Figure \ref{figure::Laplace} demonstrates the relative error for basis functions with $p=1,2,4,8$ and different uniformly refined meshes. Different orders deliver the same trend, consistent with the results in \cite{zhang2015many}.
Figure \ref{figure::LaplaceRate} illustrates the relative errors in the eigenvalue at positions $N^{\frac{1}{4}}$, $N^{\frac{1}{3}}$ and $N^{\frac{1}{2}}$ with different $p$, where the higher-order discretization demonstrates better accuracy.




\section{\zw{Application to Kohn-Sham Density Functional Theory}}
\label{chapter::KohnShamEquations}

In this section, we apply our framework to KSDFT calculations. We first give a brief introduction to the Kohn-Sham equations, including their mathematical form, discretization, and solving strategy. Subsequently, we present and discuss numerical results for various atoms and molecules.


\subsection{Kohn-Sham Equations}

Within the framework of Kohn-Sham density functional theory, the ground state is determined by solving the all-electron, time-independent Kohn-Sham equations for the orbitals $\psi_i$.
For a system with $N_{\text{ele}}$ electrons in $N_{\text{orb}}$ occupied orbitals having occupation numbers $\{f_{i}\}_{i=1}^{N_{\text{orb}}}$, and $N_{\text{nu}}$ nuclei, the spin-polarized Kohn-Sham equations take the form
\begin{equation}
    \begin{cases}
        \displaystyle \left(- \frac{1}{2} \nabla^{2} + V^{\sigma}_{\text{eff}}([\rho^{\uparrow}, \rho^{\downarrow}], \mathbf{x})\right) \psi^{\sigma}_{i}(\mathbf{x}) = \varepsilon^{\sigma}_{i} \psi^{\sigma}_{i}(\mathbf{x}), & i = 1, \dots, N^{\sigma}_{\text{orb}},    \\
        \displaystyle \int_{\Omega} \psi^{\sigma}_{i}(\mathbf{x}) \psi^{\sigma}_{j}(\mathbf{x}) \mathrm{d} \mathbf{x} = \delta_{i j},                                                                                           & i, j = 1, \dots, N^{\sigma}_{\text{orb}},
    \end{cases}
\end{equation}
where $\sigma \in \{ \uparrow, \downarrow \}$ is the spin state, $\Omega$ is the truncated domain, $\rho^{\sigma}(\mathbf{x}) = \sum_{i=1}^{N^{\sigma}_{\text{orb}}} f^{\sigma}_{i} \vert \psi_{i}(\mathbf{x}) \vert^{2}$ is the density for the corresponding spin state.



In the equations, $V_{\text{eff}}([\rho^{\uparrow}, \rho^{\downarrow}], \mathbf{x})$ denotes the Kohn-Sham potential, 
\begin{equation}\label{eqn:eff pot}
    V_{\text{eff}}([\rho^{\uparrow}, \rho^{\downarrow}], \mathbf{x}) = V_{\text{ext}}(\mathbf{x}) + V_{\text{Har}}([\rho]; \mathbf{x}) + V_{\text{xc}}([\rho^{\uparrow}, \rho^{\downarrow}]; \mathbf{x}),
\end{equation}
where $\rho(\mathbf{x}) = \rho^{\uparrow}(\mathbf{x}) + \rho^{\downarrow}(\mathbf{x})$ is the total density.
The first term stands for the Coulomb external potential
\begin{equation}
    V_{\text{ext}}(\mathbf{x}) = -\sum_{j=1}^{N_{\text{nu}}} \frac{Z_{j}}{\vert \mathbf{x} - \mathbf{X}_{j} \vert},
\end{equation}
which depicts the attraction between nuclei and electrons, where $Z_{j}$ is the charge of the $j$-th nucleus, and $\mathbf{X}_{j}$ is the position of the $j$-th nucleus.
The second term is the Hartree potential
\begin{equation}
    V_{\text{Har}}(\mathbf{x}) = \int_{\Omega} \frac{\rho(\mathbf{x}^{\prime})}{|\mathbf{x} - \mathbf{x}^{\prime}|} \mathrm{d} \mathbf{x}^{\prime},
\end{equation}
representing the Coulomb repulsion among the electrons. 
While a direct approach scales as $O(N_{\text{dof}}^{2})$ and is often computationally prohibitive, evaluating this term via the Poisson equation offers a more efficient alternative. In particular,
\begin{equation}\label{eqn:Hartree}
    \begin{cases}
        \displaystyle - \nabla^{2} V_{\text{Har}}(\mathbf{x}) = 4 \pi \rho(\mathbf{x}), & \mathbf{x} \in \Omega,          \\
        \displaystyle V_{\text{Har}}(\mathbf{x}) = \tilde{V}_{\text{Har}} (\mathbf{x}), & \mathbf{x} \in \partial \Omega,
    \end{cases}
\end{equation}
where the boundary condition is given by the multipole expansion,
\begin{equation}
    \tilde{V}_{\text{Har}}(\mathbf{r}) \approx \frac{m_{0}}{\vert \mathbf{r} - \mathbf{r}_{0} \vert} + \frac{\mathbf{p}_{0}^{T} (\mathbf{r} - \mathbf{r}_{0})}{\vert \mathbf{r} - \mathbf{r}_{0} \vert^{3}} + \frac{1}{2} \frac{(\mathbf{r} - \mathbf{r}_{0})^{T} \mathbf{Q}_{0} (\mathbf{r} - \mathbf{r}_{0})}{\vert \mathbf{r} - \mathbf{r}_{0} \vert^{5}},
\end{equation}
and $m_{0}$, $\mathbf{p}_{0}$ and $\mathbf{Q}_{0}$ are defined as
\begin{equation}
    \begin{aligned}
        m_{0} =          & \int_{\Omega} \rho(\mathbf{r}) \mathrm{d} \mathbf{r},                                                                                                                        \\
        \mathbf{p}_{0} = & \int_{\Omega} \rho(\mathbf{r}) (\mathbf{r} - \mathbf{r}_{0}) \mathrm{d} \mathbf{r} = \int_{\Omega} \rho(\mathbf{r}) \mathbf{r} \mathrm{d} \mathbf{r} - m_{0} \mathbf{r}_{0}, \\
        \mathbf{Q}_{0} = & \int_{\Omega} \bigg(\rho(\mathbf{r}) (3 (\mathbf{r} - \mathbf{r}_{0})) (\mathbf{r} - \mathbf{r}_{0})^{T} - \vert \mathbf{r} - \mathbf{r}_{0} \vert^{2} I \mathrm{d}\bigg) \mathbf{r}.
    \end{aligned}
\end{equation}
The last term corresponds to the exchange-correlation potential
\begin{equation}
    V^{\sigma}_{\text{xc}}([\rho^{\uparrow}, \rho^{\downarrow}]; \mathbf{x}) = \frac{\delta E_{\text{xc}}(\rho^{\uparrow}, \rho^{\downarrow})}{\delta \rho^{\sigma}(\mathbf{x})},
\end{equation}
resulting from the Pauli exclusion principle and other non-classical Coulomb interactions. Since analytical expressions are generally unavailable, practical calculations rely on approximations such as the local density approximation (LDA), local spin density approximation (LSDA), and the generalized gradient approximation (GGA) \cite{parr1995density}.

For the spin-unpolarized Kohn-Sham equations, the electrons of the same orbital but with opposite spin state are combined, resulting in the form
\begin{equation}
    \begin{cases}
        \displaystyle \left(- \frac{1}{2} \nabla^{2} + V_{\text{eff}}([\rho], \mathbf{x})\right) \psi_{i}(\mathbf{x}) = \varepsilon_{i} \psi_{i}(\mathbf{x}), & i = 1, \dots, N_{\text{orb}},    \\
        \displaystyle \int_{\Omega} \psi_{i}(\mathbf{x}) \psi_{j}(\mathbf{x}) \mathrm{d} \mathbf{x} = \delta_{i j},                                           & i, j = 1, \dots, N_{\text{orb}},
    \end{cases}
\end{equation}
where $\delta_{i j}$ is the Kronecker delta symbol and $\rho(\mathbf{x}) = \sum_{i=1}^{N_{\text{orb}}} f_{i} \vert \psi_{i}(\mathbf{x}) \vert^{2}$.

Next, the implementation details of our method are presented.
For the basis functions $\{ \phi_{k}(\mathbf{x}) \}_{k=1}^{N_{\text{dof}}}$, multiplying both sides of the equation by $\phi_{k}(\mathbf{x})$ and applying the divergence theorem yields

\begin{equation}
    \frac{1}{2} \int_{\Omega} \nabla \psi_{i}(\mathbf{x}) \nabla \phi_{k}(\mathbf{x}) \mathrm{d} \mathbf{x} + \int_{\Omega} V_{\text{eff}}([\rho^{\uparrow}, \rho^{\downarrow}],\mathbf{x}) \psi_{i}(\mathbf{x}) \phi_{k}(\mathbf{x}) \mathrm{d} \mathbf{x} = \int_{\Omega} \varepsilon_{i} \psi_{i}(\mathbf{x}) \phi_{k}(\mathbf{x}) \mathrm{d} \mathbf{x}.
\end{equation}

\noindent Substituting $\psi_{i}(\mathbf{x}) = \sum_{j=1}^{N_{\text{dof}}} \psi_{i}^{[j]} \phi_{j}(\mathbf{x})$ we obtain

\begin{equation}
    \begin{aligned}
        \sum_{j=1}^{N_{\text{dof}}} \psi_{i}^{[j]} \left(\frac{1}{2} \int_{\Omega} \nabla \phi_{j}(\mathbf{x}) \nabla \phi_{k}(\mathbf{x}) \mathrm{d} \mathbf{x} + \int_{\Omega} V_{\text{eff}}([\rho^{\uparrow}, \rho^{\downarrow}],\mathbf{x}) \phi_{j}(\mathbf{x}) \phi_{k}(\mathbf{x}) \mathrm{d} \mathbf{x}\right) = \varepsilon_{i} \sum_{j=1}^{N_{\text{dof}}} \psi_{i}^{[j]} \int_{\Omega} \phi_{j}(\mathbf{x}) \phi_{k}(\mathbf{x}) \mathrm{d} \mathbf{x}.
    \end{aligned}
\end{equation}

\noindent Therefore, elements of the matrices $S, V, M$ can be defined as

\begin{equation}
    \begin{aligned}
        S_{j, k} = & \int_{\Omega} \nabla \phi_{j}(\mathbf{x}) \nabla \phi_{k}(\mathbf{x}) \mathrm{d} \mathbf{x},              \\
        V_{j, k}[\Psi] = & \int_{\Omega} V_{\text{eff}}([\rho^{\uparrow}, \rho^{\downarrow}],\mathbf{x}) \phi_{j}(\mathbf{x}) \phi_{k}(\mathbf{x}) \mathrm{d} \mathbf{x}, \\
        M_{j, k} = & \int_{\Omega} \phi_{j}(\mathbf{x}) \phi_{k}(\mathbf{x}) \mathrm{d} \mathbf{x}.
    \end{aligned}
\end{equation}

\noindent Consequently, the discretized Kohn-Sham equation can be written as a generalized eigenvalue problem

\begin{equation}\label{eqn:discretized KS}
    \begin{cases}
        \displaystyle H \Psi = \left( \frac{1}{2} S + V[\Psi] \right) \Psi = M \Psi \Lambda, \\
        \displaystyle X^{T} M X = I,
    \end{cases}
\end{equation}

\noindent where $\Lambda = \text{diag}(\varepsilon_{1}, \dots, \varepsilon_{N_{\text{orb}}})$ and $\Psi = (\psi_{1}, \dots, \psi_{N_{\text{orb}}})$ are the smallest $N_{\text{orb}}$ eigenpairs of the generalized eigenvalue problem.

For the large linear system arising from the spectral element method, we utilize the locally optimal block preconditioned conjugate gradient method \cite{knyazev2001toward, knyazev2007block}. Moreover, the Jacobi preconditioner \cite{saad2003iterative, golub2013matrix} is employed to accelerate the convergence.

To handle the nonlinearity in Equation (\ref{eqn:discretized KS}), the self-consistent field (SCF) strategy \cite{ehrenreich1959self} is used, as detailed in Algorithm \ref{algo::scf}. In particular, Pulay's method \cite{pulay1980convergence} is adopted to further accelerate the convergence.


\begin{algorithm}[H]
    \caption{Self-consistent field iteration.}
    \label{algo::scf}
    \KwIn{Initial spin densities: $\rho^{\uparrow}_{*}(\mathbf{x}), \rho^{\downarrow}_{*}(\mathbf{x})$; Energy convergence tolerance: $\delta$; Density mixing operator: $U(\rho_{\mathrm{old}}, \rho_{\mathrm{new}})$.}
    \KwOut{Converged electron densities: $\rho^{\uparrow}(\mathbf{x}), \rho^{\downarrow}(\mathbf{x})$.}
    \BlankLine
    Initialize the electron densities as $\rho^{\uparrow}(\mathbf{x}) \leftarrow \rho_{*}^{\uparrow}(\mathbf{x}), \rho^{\downarrow}(\mathbf{x}) \leftarrow \rho_{*}^{\downarrow}(\mathbf{x})$\;
    Compute the total energy $E$ from $\rho^{\uparrow}(\mathbf{x}), \rho^{\downarrow}(\mathbf{x})$\;
    \Repeat{$\vert E_{\text{prev}} - E \vert \leq \delta$}{
        \tcc{Save the previous densities and energy}
        $\rho^{\uparrow}_{\text{prev}}(\mathbf{x}) \leftarrow \rho^{\uparrow}(\mathbf{x})$, $\rho^{\downarrow}_{\text{prev}}(\mathbf{x}) \leftarrow \rho^{\downarrow}(\mathbf{x})$, $E_{\text{prev}} \leftarrow E$\;
        Compute the Hartree potential by solving Equation (\ref{eqn:Hartree})\;
        Construct the Kohn-Sham potential $V_{\text{eff}}([\rho^{\uparrow}, \rho^{\downarrow}], \mathbf{x})$ via Equation (\ref{eqn:eff pot})\;
        Solve the Kohn-Sham equations (\ref{eqn:discretized KS}) to obtain candidate densities  $\tilde{\rho}^{\uparrow}(\mathbf{x}), \tilde{\rho}^{\downarrow}(\mathbf{x})$\;
        \tcc{Mix the previous and candidate densities}
        $\rho^{\uparrow}(\mathbf{x}) \leftarrow U(\rho_{\text{prev}}^{\uparrow}(\mathbf{x}), \tilde{\rho}^{\uparrow}(\mathbf{x}))$, $\rho^{\downarrow}(\mathbf{x}) \leftarrow U(\rho_{\text{prev}}^{\downarrow}(\mathbf{x}), \tilde{\rho}^{\downarrow}(\mathbf{x}))$\;
    }
\end{algorithm}

\begin{figure}[H]
    \centering
    \subfloat[Density of Oxygen atom.]{\includegraphics[width=0.45\textwidth]{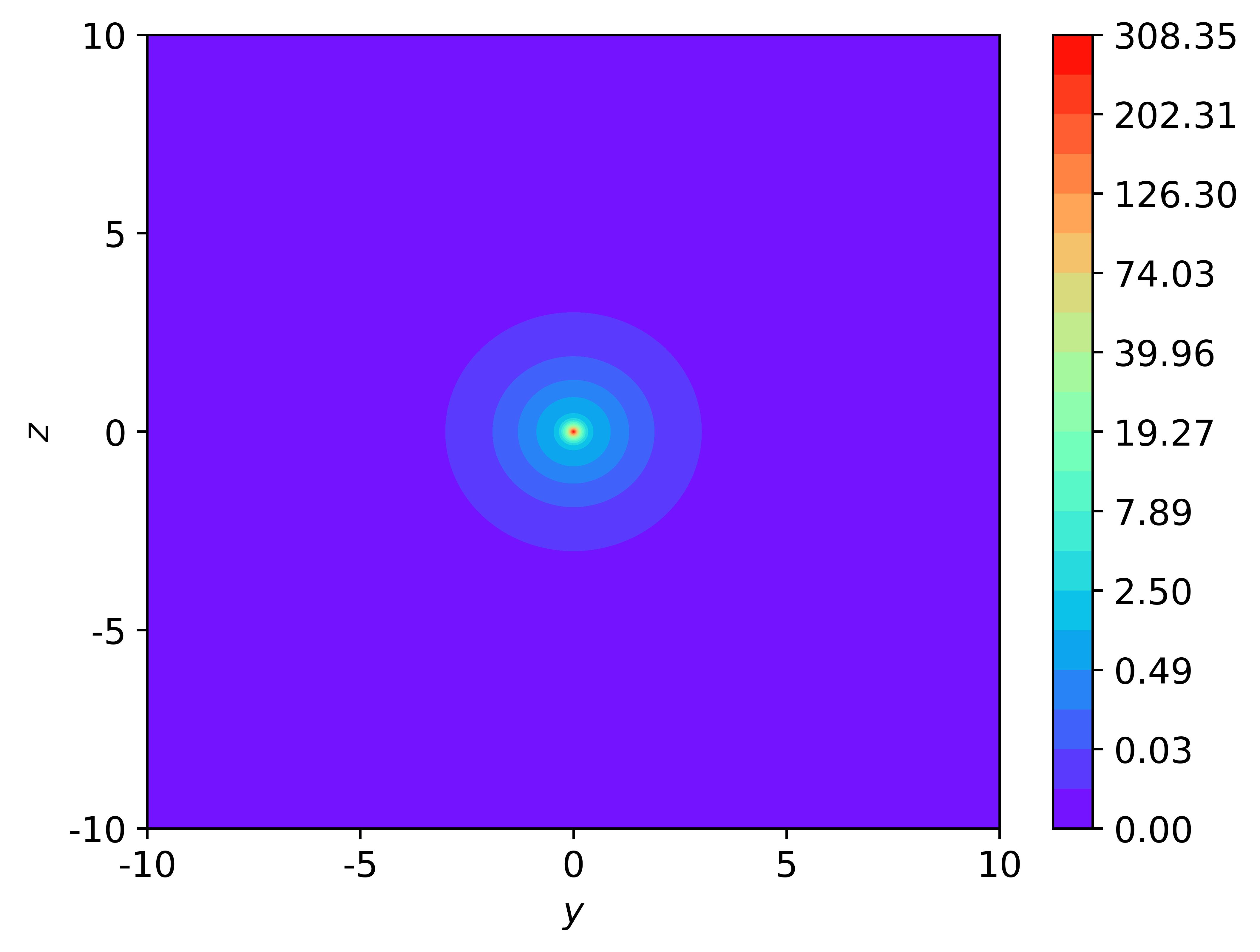} \label{figure::KohnShamDensity::O-pre}}
    \subfloat[Density of Ethene.]{\includegraphics[width=0.45\textwidth]{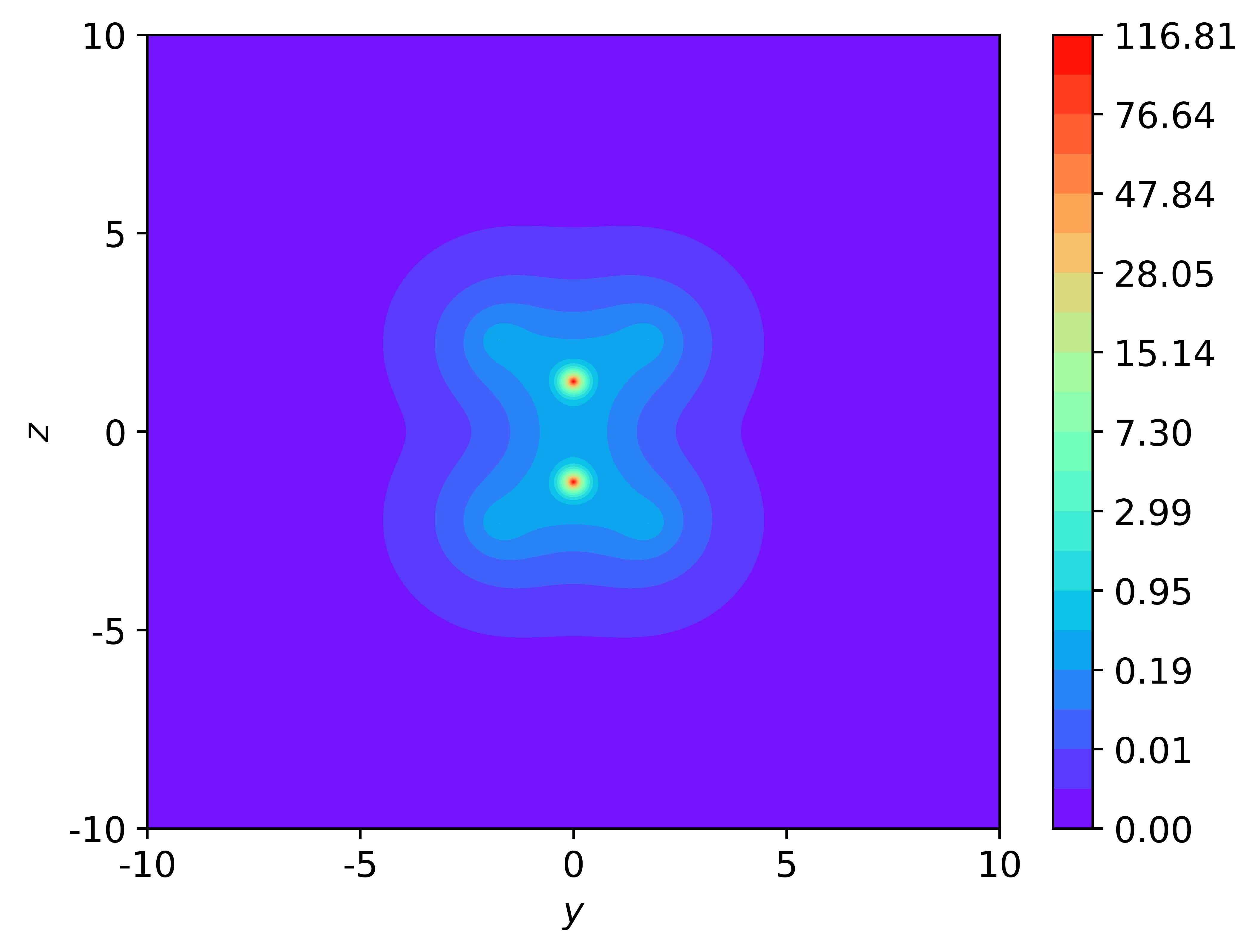} \label{figure::KohnShamDensity::C2H4-pre}}
    \caption{Densities of Oxygen atom (left) and Ethene (right).}
    \label{figure::KohnShamDensity-pre}
\end{figure}

\zw{
Furthermore, due to the Coulomb singularities from nuclei, both localization of the electron density and rapid variations in the orbital wave functions are expected.
This behavior is confirmed by Figure \ref{figure::KohnShamDensity-pre}, which presents the electron densities of the oxygen atom and the ethylene molecule, where pronounced variations near the nuclei can be clearly observed. 
Therefore, we consider the following error indicator for the $i$-th orbital,
\begin{equation}
    \eta_{i, \mathcal{T}} = \sqrt{\left(\frac{h_{\mathcal{T}}}{p}\right)^{2} \int_{\mathcal{T}} \left\vert \left(-\frac{1}{2}\nabla^{2} + V_{\text{eff}}(\mathbf{x}) - \varepsilon_{i} \right) \psi_{i}(\mathbf{x}) \right\vert^{2} \mathrm{d} \mathbf{x} + \frac{1}{2} \sum_{\mathcal{F} \in \partial \mathcal{T}} \frac{h_{\mathcal{F}}}{p} \int_{\mathcal{F}} \left\vert \frac{\partial \psi_{i, \mathcal{T}}}{\partial \mathbf{n}_{\mathcal{F}}} - \frac{\partial \psi_{i, \mathcal{F}}}{\partial \mathbf{n}_{\mathcal{F}}} \right\vert^{2} \mathrm{d} \mathbf{x}},
\end{equation}
with $\psi_{i, \mathcal{T}}$ and $\psi_{i, \mathcal{F}}$ denoting the numerical solution on $\mathcal{T}$ and the adjacent tetrahedron sharing the face $\mathcal{F}$, respectively.}

\subsection{\zw{Numerical Results}}


In this subsection, we present a set of numerical experiments to demonstrate the effectiveness of our method in all-electron electronic structure calculations. In all tests, the exchange-correlation potential is evaluated using the \verb|Libxc| library \cite{marques2012libxc}.
Specifically, we first present benchmark results for atomic and molecular systems to validate the accuracy of the proposed method.
Next, we demonstrate the behavior of the adaptive refinement strategy using specific atomic and molecular systems. 
Subsequently, we investigate the effect of interpolation between adaptive meshes.
Finally, we assess the parallel performance of the proposed method to show its scalability.

\subsubsection{\zw{Ground-State Energy Accuracy}}

We begin with two tables showing that the proposed method achieves satisfactory accuracy in the computed energies. Reference energies for atoms are taken from the NIST database \cite{johnson1999nist, kotochigova2005atomic}, while those for molecules are obtained from the Elk code \cite{elk}. 

In Table \ref{table::Atom}, we present the numerical results for the first two rows of the periodic table, while those for molecules are given in Table \ref{table::Molecular}. In both tables, the energies obtained using LDA or LSDA (denoted as LSD in the tables) are listed in the 3rd column, with the corresponding reference energies shown in the adjacent column. The absolute errors for atoms are reported in the last column, whereas the relative errors are provided in the corresponding columns. The results indicate that our method achieves consistent accuracy across all tested systems.


\begin{table}[H]
    \centering
    \caption{The energies and absolute errors for atoms.}
    \begin{tabular}{|c|c|c|c|c|c|c|c|c|c|}
        \hline
        \multirow{2}{*}{Atom}        & \multirow{2}{*}{$V_{\text{xc}}$} & \multicolumn{2}{c|}{Energy ($\text{a.u.}$)} & \multirow{2}{*}{Abs Error}                         \\
        \cline{3-4}
                                     &                                  & Result                                         & Ref Data                   &                       \\
        \hline
        \multirow{2}{*}{$\text{H}$}  & LDA                              & $-0.445479$                                    & $-0.445671$                & $1.92 \times 10^{-4}$ \\
                                     & LSD                              & $-0.478535$                                    & $-0.478671$                & $1.36 \times 10^{-4}$ \\
        \hline
        \multirow{2}{*}{$\text{He}$} & LDA                              & $-2.834683$                                    & $-2.834836$                & $1.53 \times 10^{-4}$ \\
                                     & LSD                              & $-2.834480$                                    & $-2.834836$                & $3.56 \times 10^{-4}$ \\
        \hline
        \multirow{2}{*}{$\text{Li}$} & LDA                              & $-7.334602$                                    & $-7.335195$                & $5.93 \times 10^{-4}$ \\
                                     & LSD                              & $-7.343524$                                    & $-7.343957$                & $4.33 \times 10^{-4}$ \\
        \hline
        \multirow{2}{*}{$\text{Be}$} & LDA                              & $-14.447646$                                   & $-14.447209$               & $4.37 \times 10^{-4}$ \\
                                     & LSD                              & $-14.447081$                                   & $-14.447209$               & $1.28 \times 10^{-4}$ \\
        \hline
        \multirow{2}{*}{$\text{B}$}  & LDA                              & $-24.344685$                                   & $-24.344198$               & $4.87 \times 10^{-4}$ \\
                                     & LSD                              & $-24.353260$                                   & $-24.353614$               & $3.54 \times 10^{-4}$ \\
        \hline
        \multirow{2}{*}{$\text{C}$}  & LDA                              & $-37.425950$                                   & $-37.425749$               & $2.01 \times 10^{-4}$ \\
                                     & LSD                              & $-37.469826$                                   & $-37.470031$               & $2.05 \times 10^{-4}$ \\
        \hline
        \multirow{2}{*}{$\text{N}$}  & LDA                              & $-54.024477$                                   & $-54.025016$               & $5.39 \times 10^{-4}$ \\
                                     & LSD                              & $-54.136308$                                   & $-54.136799$               & $4.91 \times 10^{-4}$ \\
        \hline
        \multirow{2}{*}{$\text{O}$}  & LDA                              & $-74.473419$                                   & $-74.473077$               & $3.42 \times 10^{-4}$ \\
                                     & LSD                              & $-74.527253$                                   & $-74.527410$               & $1.57 \times 10^{-4}$ \\
        \hline
        \multirow{2}{*}{$\text{F}$}  & LDA                              & $-99.099430$                                   & $-99.099648$               & $2.18 \times 10^{-4}$ \\
                                     & LSD                              & $-99.113904$                                   & $-99.114192$               & $2.88 \times 10^{-4}$ \\
        \hline
        \multirow{2}{*}{$\text{Ne}$} & LDA                              & $-128.233166$                                  & $-128.233481$              & $3.15 \times 10^{-4}$ \\
                                     & LSD                              & $-128.233668$                                  & $-128.233481$              & $1.87 \times 10^{-4}$ \\
        \hline
    \end{tabular}
    \label{table::Atom}
\end{table}

\begin{table}[H]
    \centering
    \caption{The energies and relative errors for moleculars.}
    \label{table::Molecular}
    \begin{tabular}{|c|c|c|c|c|c|c|c|c|c|}
        \hline
        \multirow{2}{*}{Molecule}       & \multirow{2}{*}{$V_{\text{xc}}$} & \multicolumn{2}{c|}{Energy ($\text{a.u.}$)} & \multirow{2}{*}{Rel Error}                         \\
        \cline{3-4}
                                         &                                  & Result                                         & Ref Data                   &                       \\
        \hline
        \multirow{2}{*}{$\text{LiH}$}    & LDA                              & $-7.919592$                                    & $-7.921490$                & $2.40 \times 10^{-4}$ \\
                                         & LSD                              & $-7.919601$                                    & $-7.921493$                & $2.39 \times 10^{-4}$ \\
        \hline
        \multirow{2}{*}{$\text{Li}_{2}$} & LDA                              & $-14.724862$                                   & $-14.728055$               & $2.17 \times 10^{-4}$ \\
                                         & LSD                              & $-14.724863$                                   & $-14.728053$               & $2.17 \times 10^{-4}$ \\
        \hline
        \multirow{2}{*}{$\text{B}_{2}$}  & LDA                              & $-48.827185$                                   & $-48.844055$               & $3.45 \times 10^{-4}$ \\
                                         & LSD                              & $-48.855186$                                   & $-48.867152$               & $2.45 \times 10^{-4}$ \\
        \hline
        \multirow{2}{*}{$\text{N}_{2}$}  & LDA                              & $-108.700049$                                  & $-108.764456$              & $5.92 \times 10^{-4}$ \\
                                         & LSD                              & $-108.699898$                                  & $-108.764456$              & $5.94 \times 10^{-4}$ \\
        \hline
        \multirow{2}{*}{$\text{CO}$}     & LDA                              & $-112.478383$                                  & $-112.551113$              & $6.46 \times 10^{-4}$ \\
                                         & LSD                              & $-112.478197$                                  & $-112.551112$              & $6.48 \times 10^{-4}$ \\
        \hline
        \multirow{2}{*}{$\text{NO}$}     & LDA                              & $-128.975942$                                  & $-129.064883$              & $6.89 \times 10^{-4}$ \\
                                         & LSD                              & $-128.985359$                                  & $-129.073061$              & $6.79 \times 10^{-4}$ \\
        \hline
        \multirow{2}{*}{$\text{O}_{2}$}  & LDA                              & $-149.303961$                                  & $-149.416139$              & $7.51 \times 10^{-4}$ \\
                                         & LSD                              & $-149.341754$                                  & $-149.448999$              & $7.18 \times 10^{-4}$ \\
        \hline
        \multirow{2}{*}{$\text{F}_{2}$}  & LDA                              & $-198.354512$                                  & $-198.492310$              & $6.94 \times 10^{-4}$ \\
                                         & LSD                              & $-198.354081$                                  & $-198.491766$              & $6.94 \times 10^{-4}$ \\
        \hline
    \end{tabular}
\end{table}

\subsubsection{\zw{$h$-adaptive Solution}}

We then focus on two specific systems, the oxygen atom and ethene, to examine the $h$-adaptive behavior, encompassing the total energy, the final density distribution, and the corresponding adaptive mesh. The stopping criterion for the $h$-adaptive strategy is set to $10^{-4}$.

\begin{figure}[H]
    \centering
    \subfloat[Oxygen atom (with LDA).]{\includegraphics[width=0.45\textwidth]{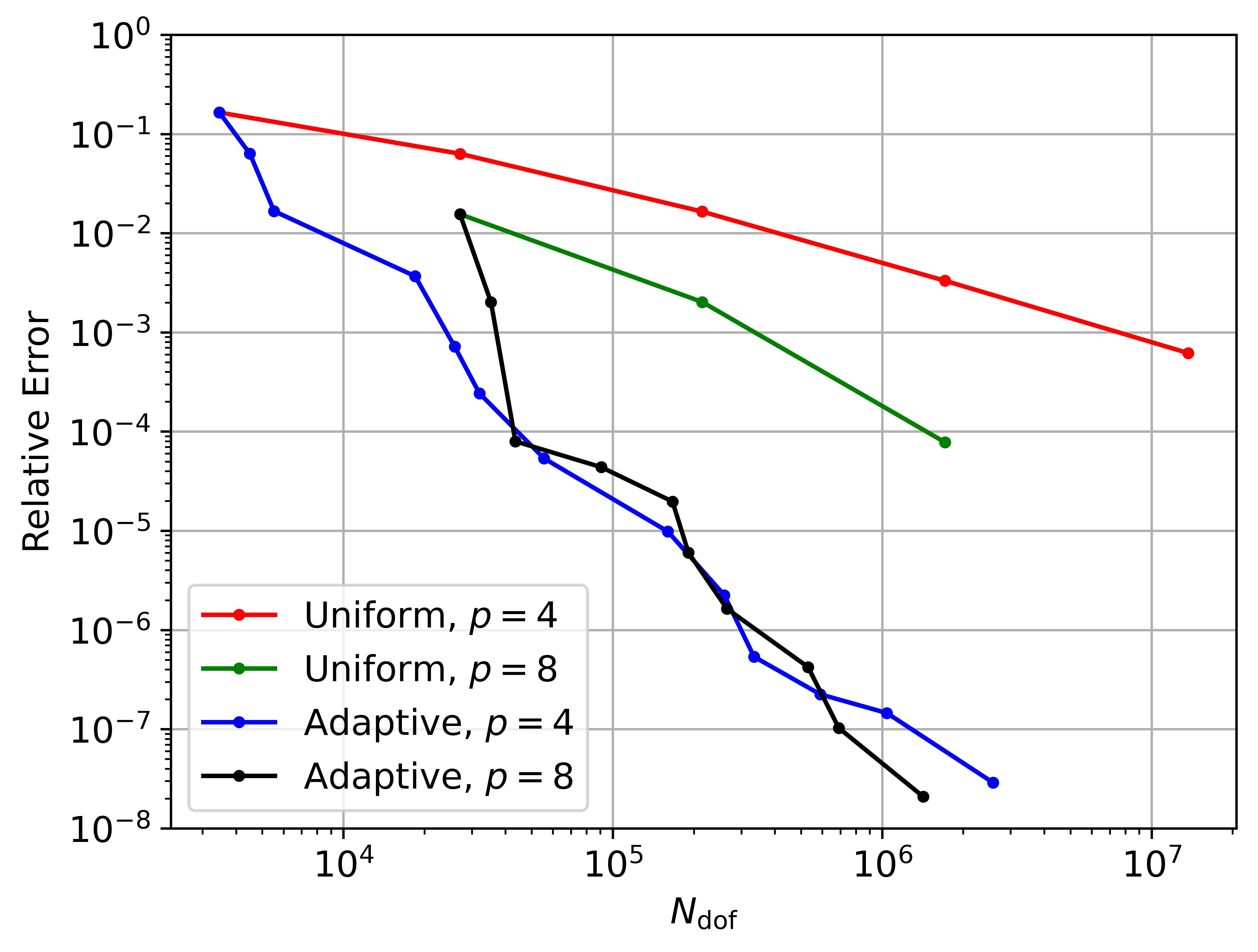}}
    \subfloat[Oxygen atom (with LSDA).]{\includegraphics[width=0.45\textwidth]{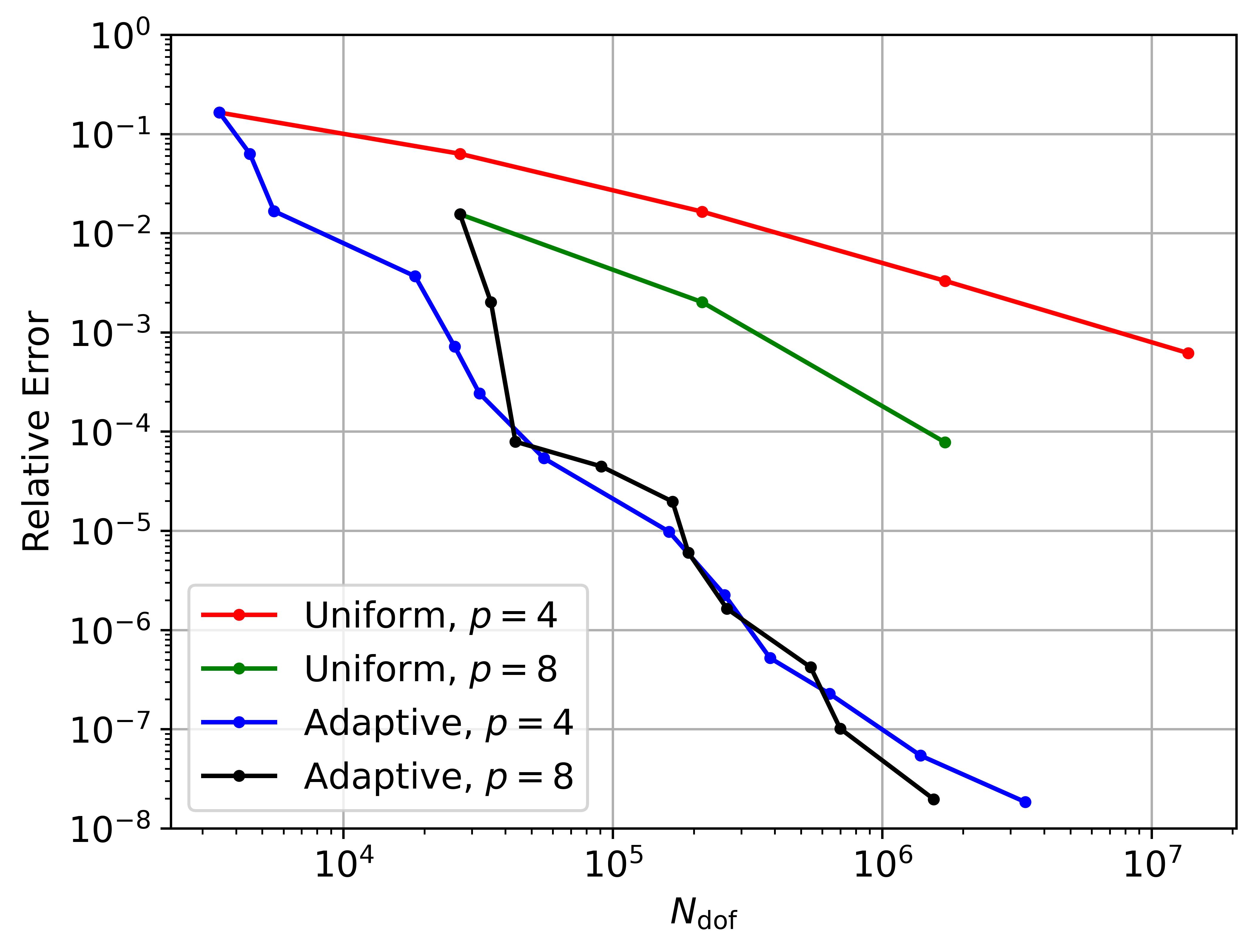}} \
    \subfloat[Ethene (with LDA)]{\includegraphics[width=0.45\textwidth]{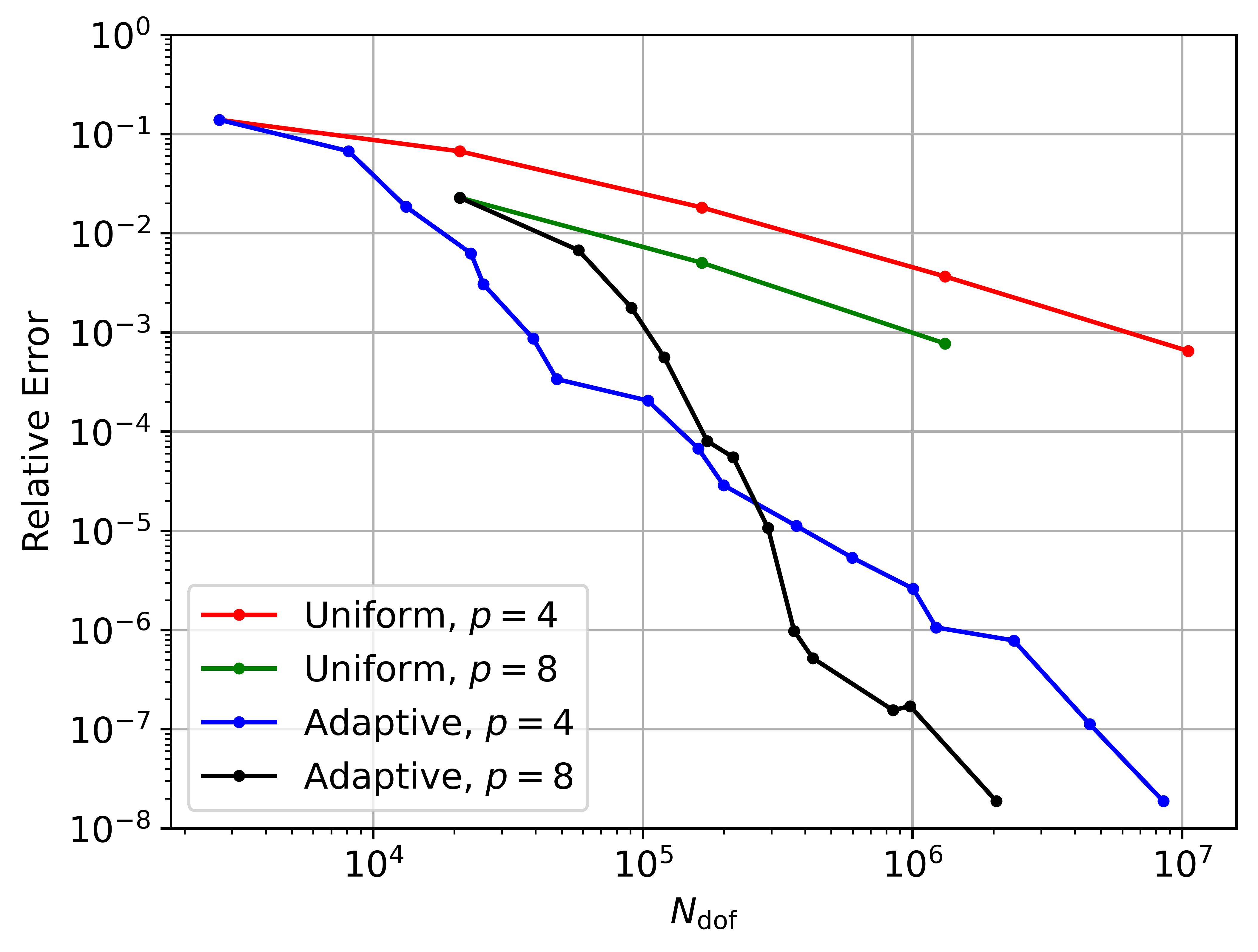}}
    \subfloat[Ethene (with LSDA)]{\includegraphics[width=0.45\textwidth]{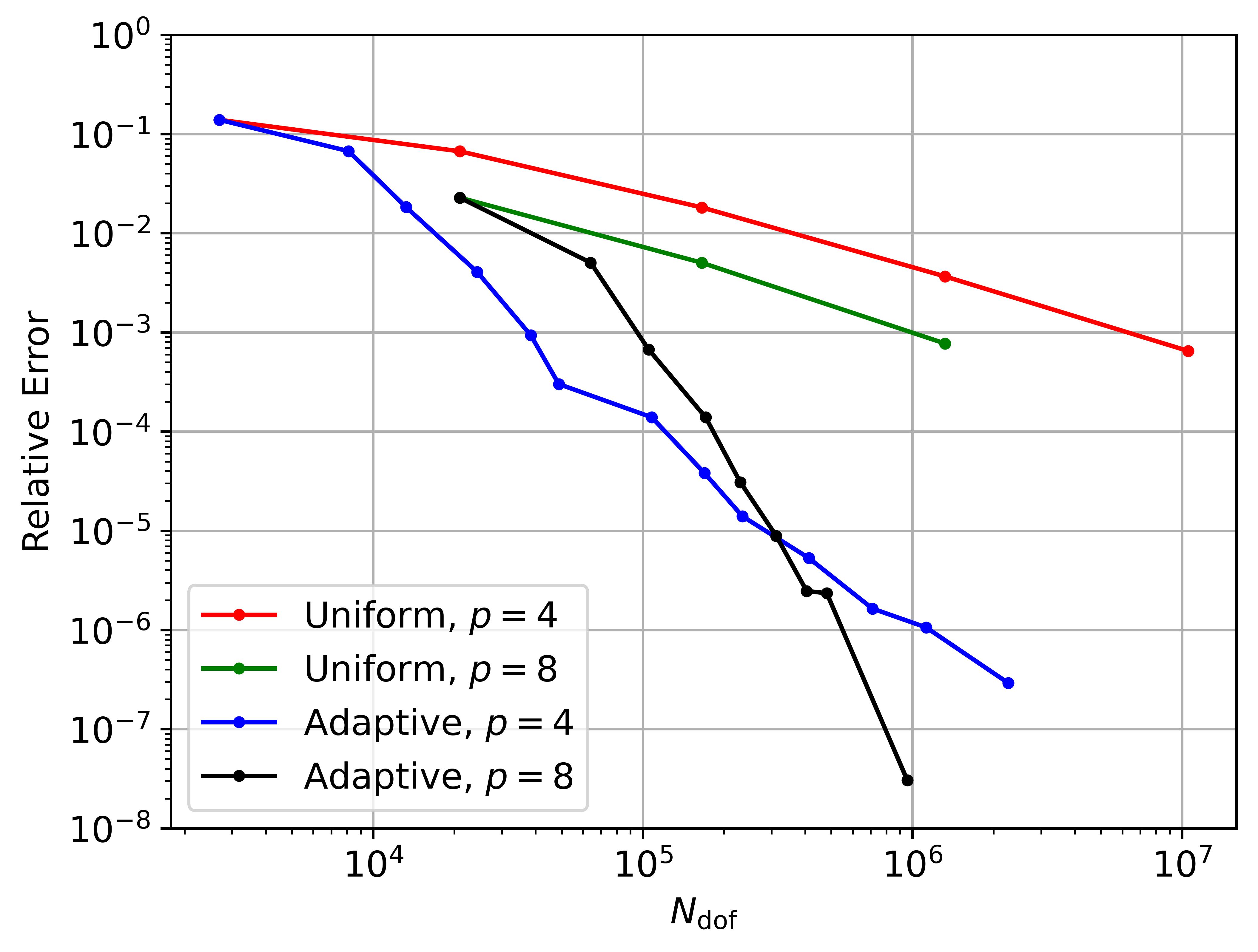}}
    \caption{The error behavior during $h$-refinement for Oxygen atom (top)/ Ethene (bottom) using exchange-correlation function LDA (left)/ LSDA (right).}
    \label{figure::KohnShamDofError}
\end{figure}

In Figure \ref{figure::KohnShamDofError}, we present the evolution of the ground-state energy during the $h$-adaptive process for an oxygen atom and ethene, along with a comparison of different exchange-correlation functionals. In each sub-figure, results for spectral element orders $p=4$ and $p=8$ are shown, as well as those obtained using global refinement. Specifically, the reference energy for the oxygen atom is given by the NIST database. For ethene, since the final energies exhibit a relative difference of $4\times 10^{-8}$ between the $p=4$ and $p=8$ cases, their average is adopted as the reference value. It can be observed from the figures that higher-order basis functions yield better accuracy, demonstrating their advantage in resolving the solution with fewer DoFs.

\begin{figure}[H]
    \centering
    \subfloat[Density of Oxygen atom.]{\includegraphics[width=0.45\textwidth]{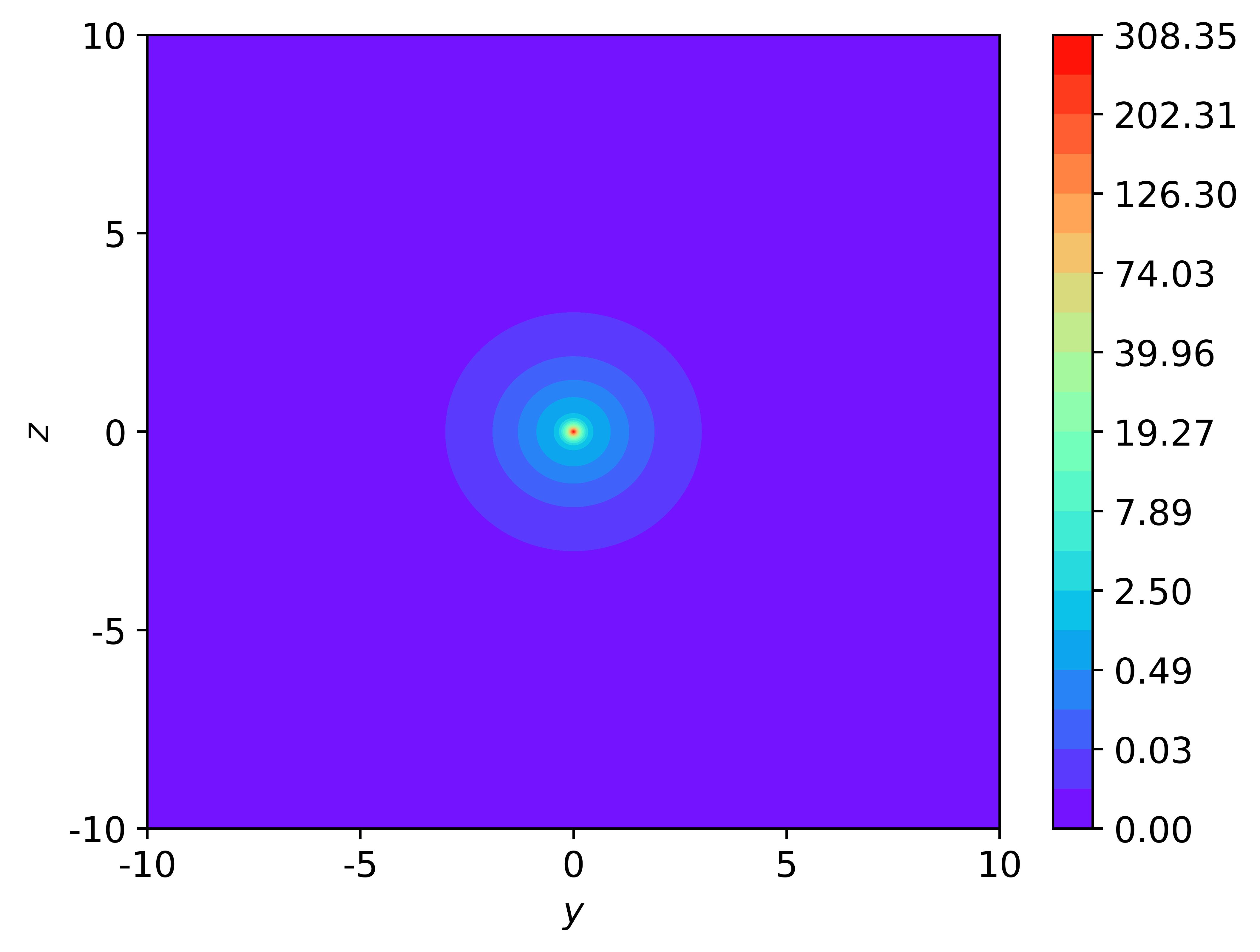} \label{figure::KohnShamDensity::O}}
    \subfloat[Density of Ethene.]{\includegraphics[width=0.45\textwidth]{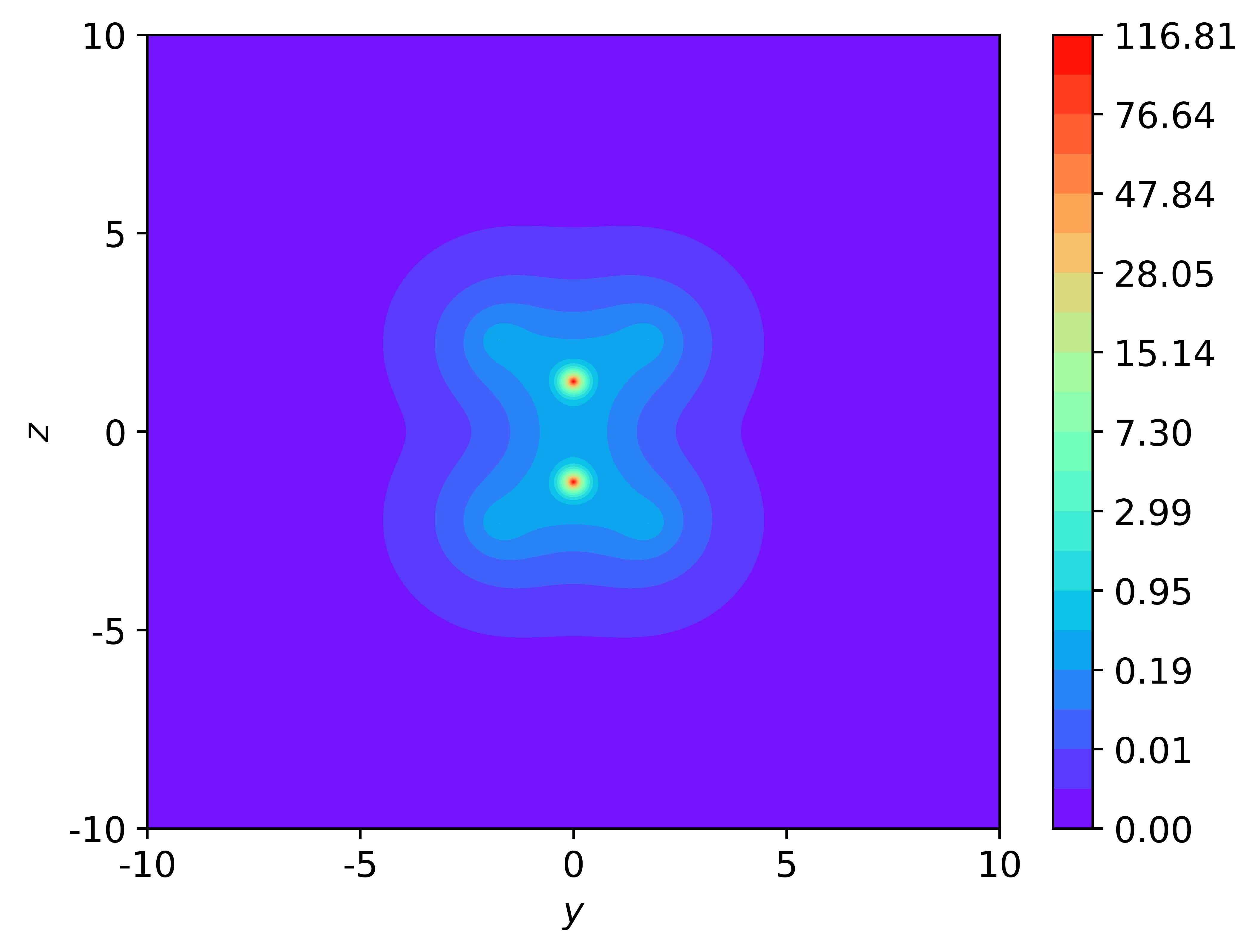} \label{figure::KohnShamDensity::C2H4}}
    \caption{Densities of Oxygen atom (left)/ Ethene (right).}
    \label{figure::KohnShamDensity}
\end{figure}

Next, Figure \ref{figure::KohnShamDensity} exhibits the final electron densities of the oxygen atom and ethene, with the electron density concentrated around the nuclei.

\begin{figure}[H]
    \centering
    \subfloat[Iter 3]{\includegraphics[width=0.25\textwidth]{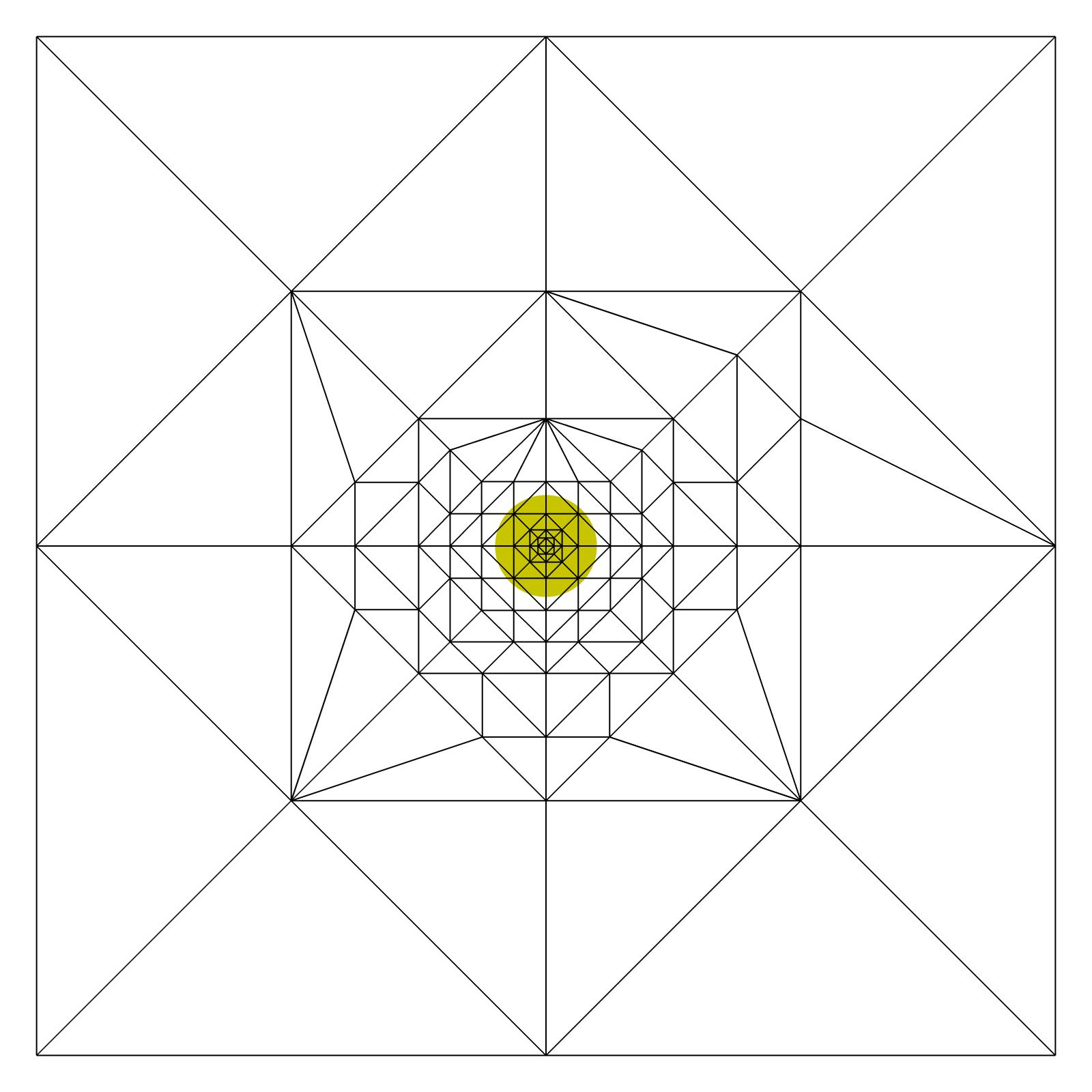}}
    \subfloat[Iter 4]{\includegraphics[width=0.25\textwidth]{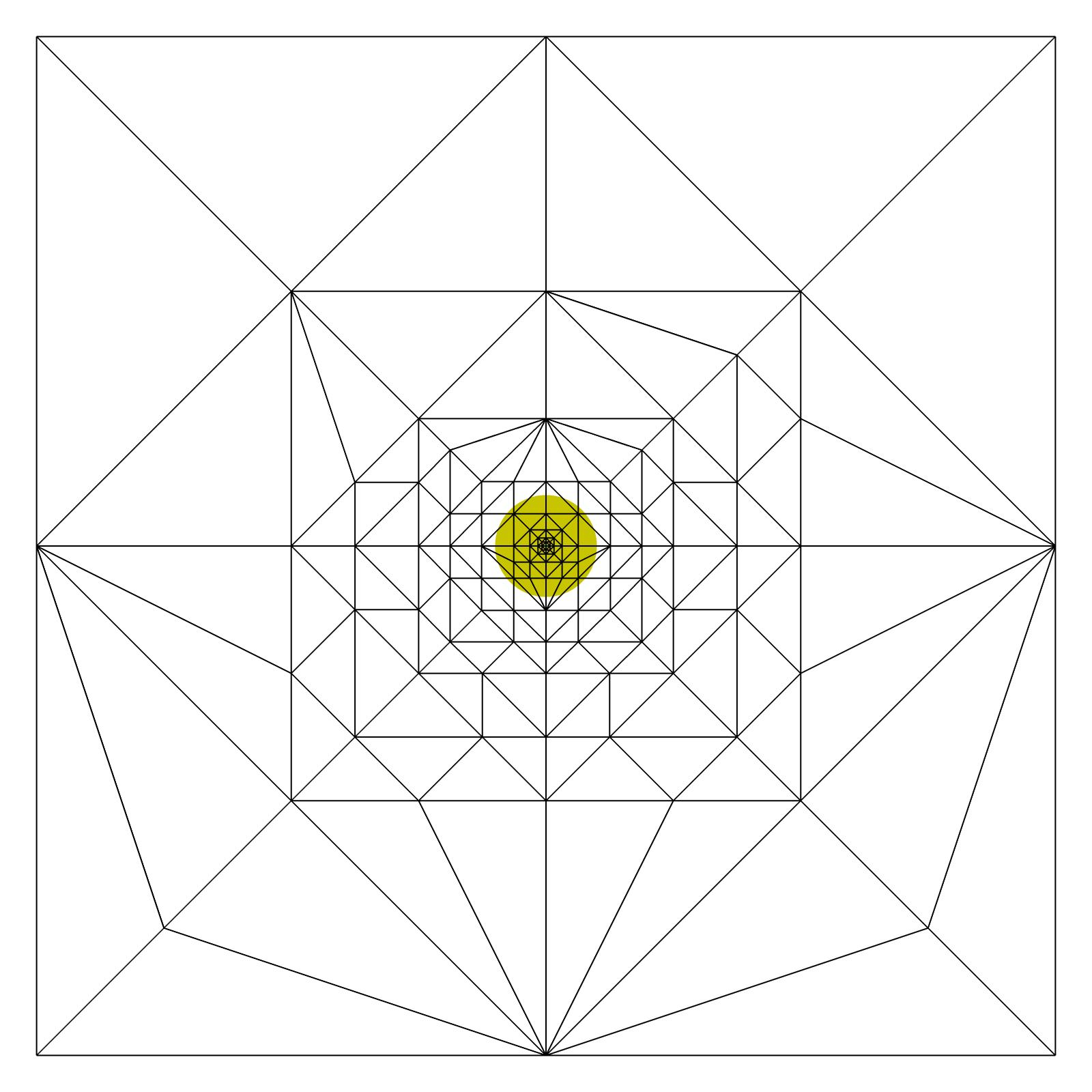}}
    \subfloat[Iter 5]{\includegraphics[width=0.25\textwidth]{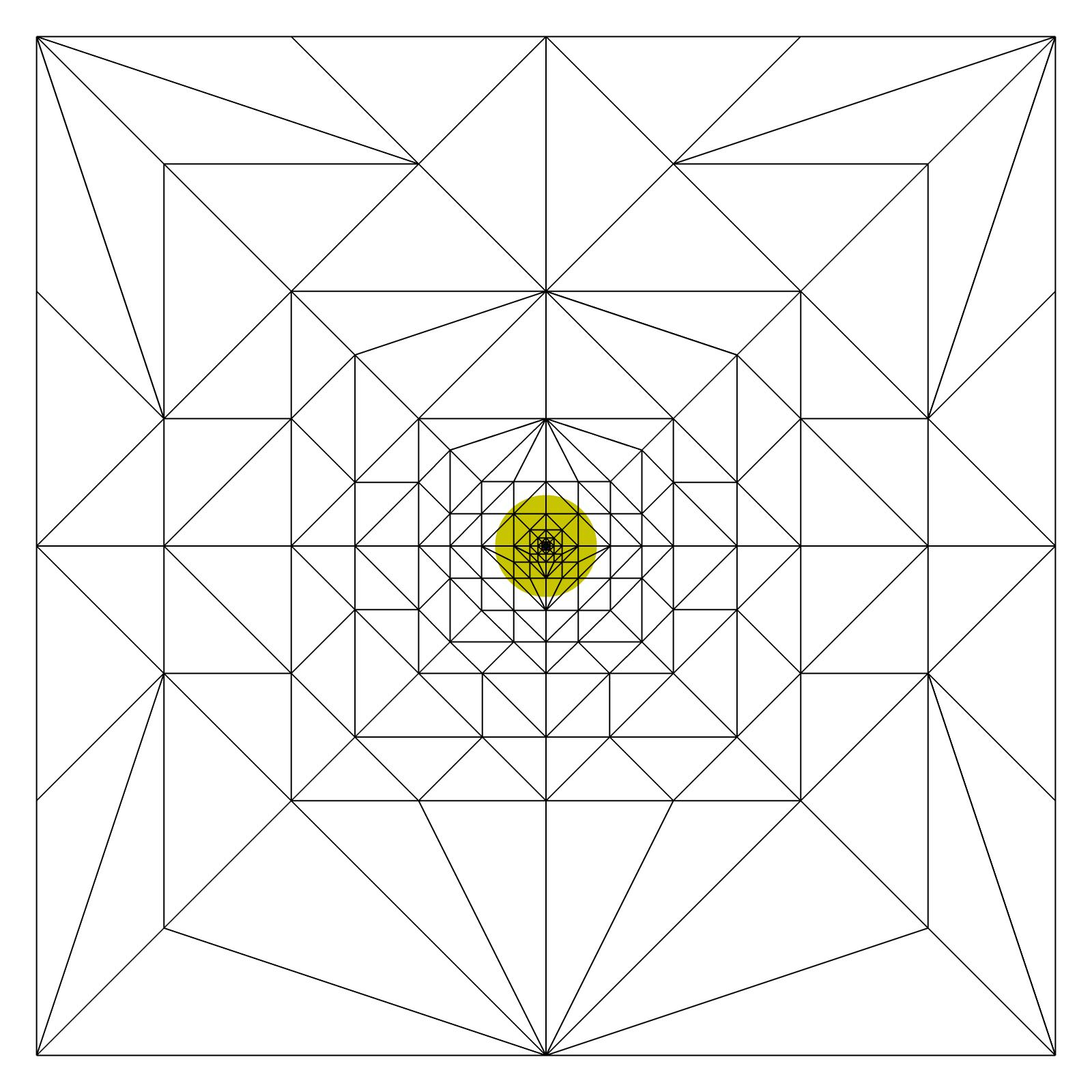}}
    \subfloat[Iter 6]{\includegraphics[width=0.25\textwidth]{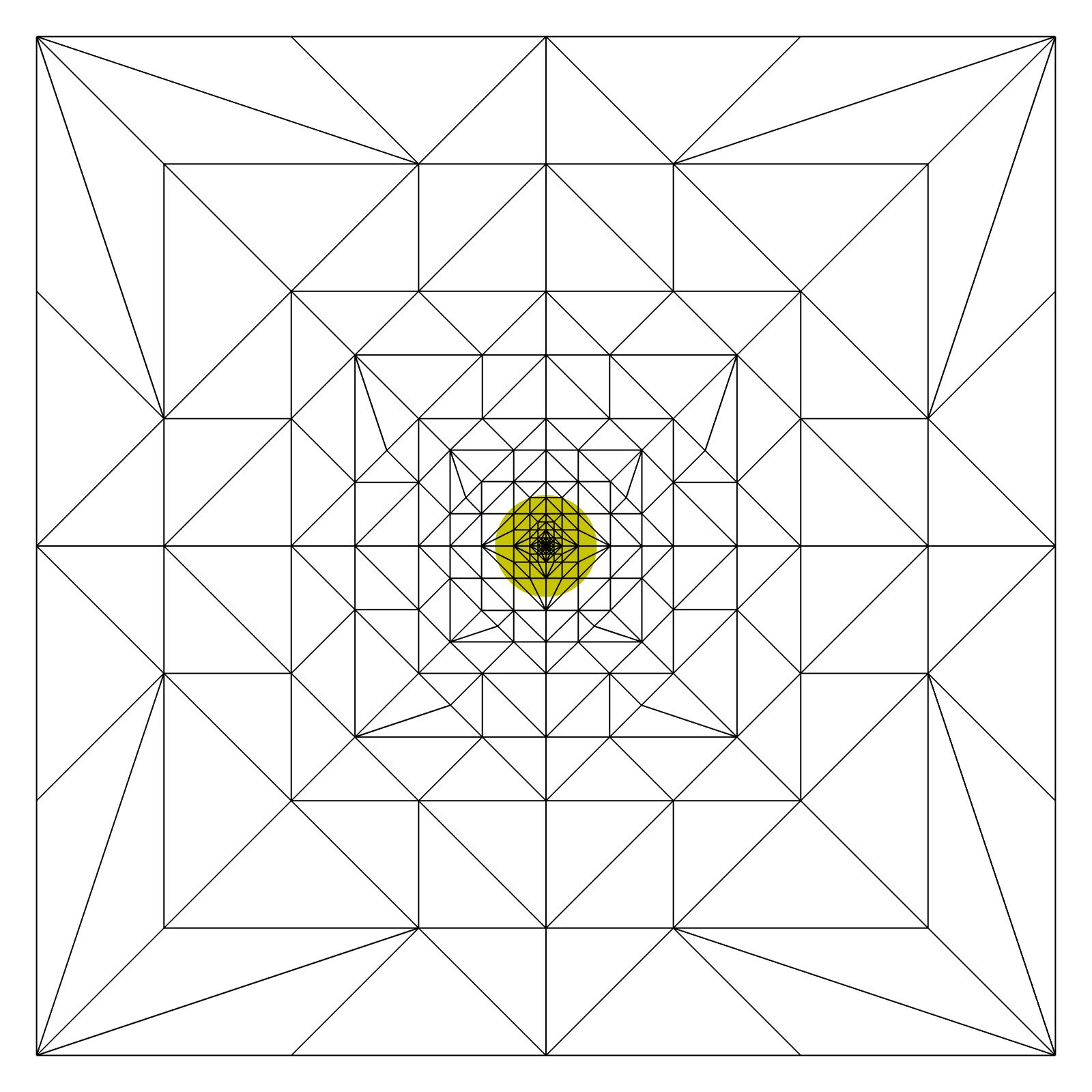}} \
    \subfloat[Iter 7]{\includegraphics[width=0.25\textwidth]{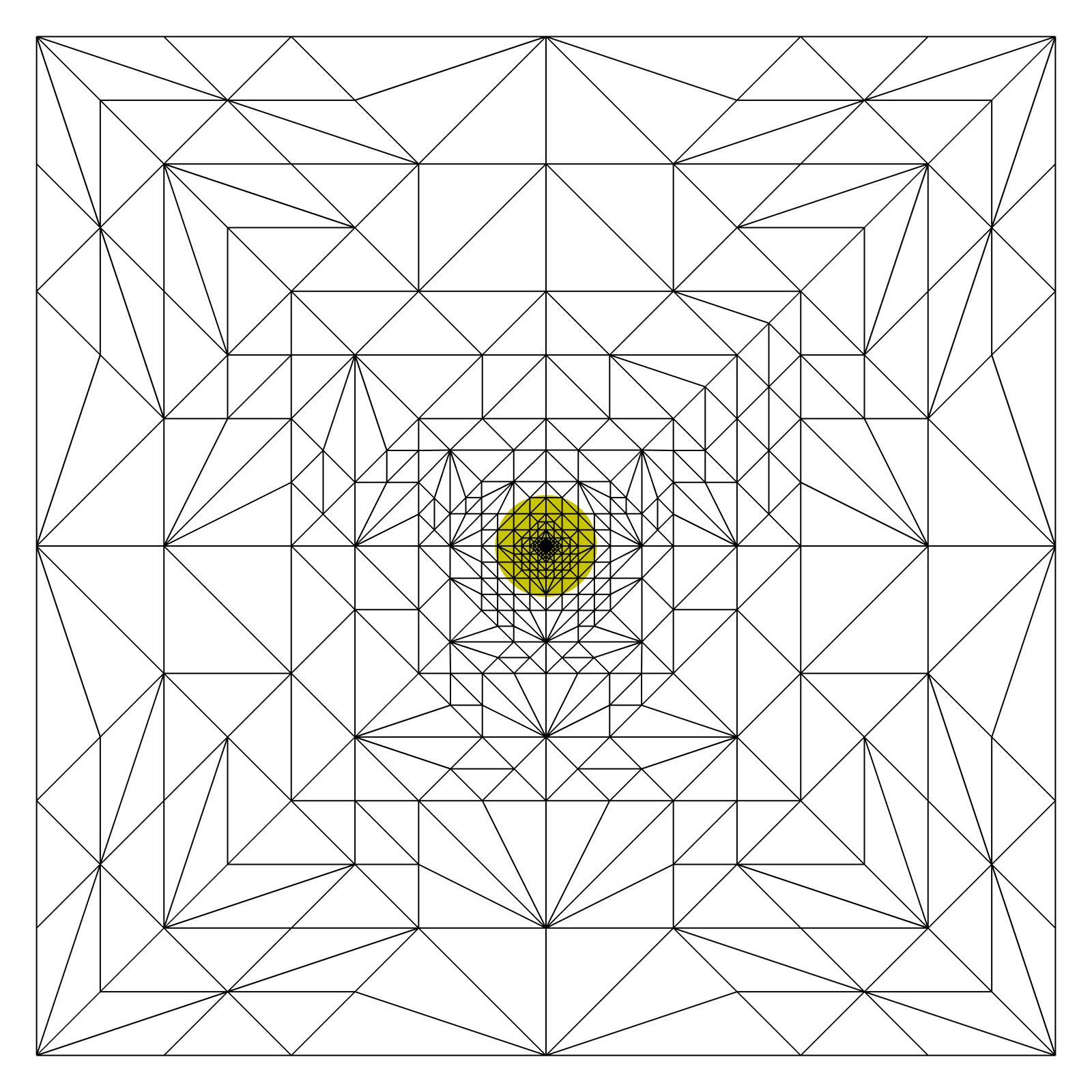}}
    \subfloat[Iter 8]{\includegraphics[width=0.25\textwidth]{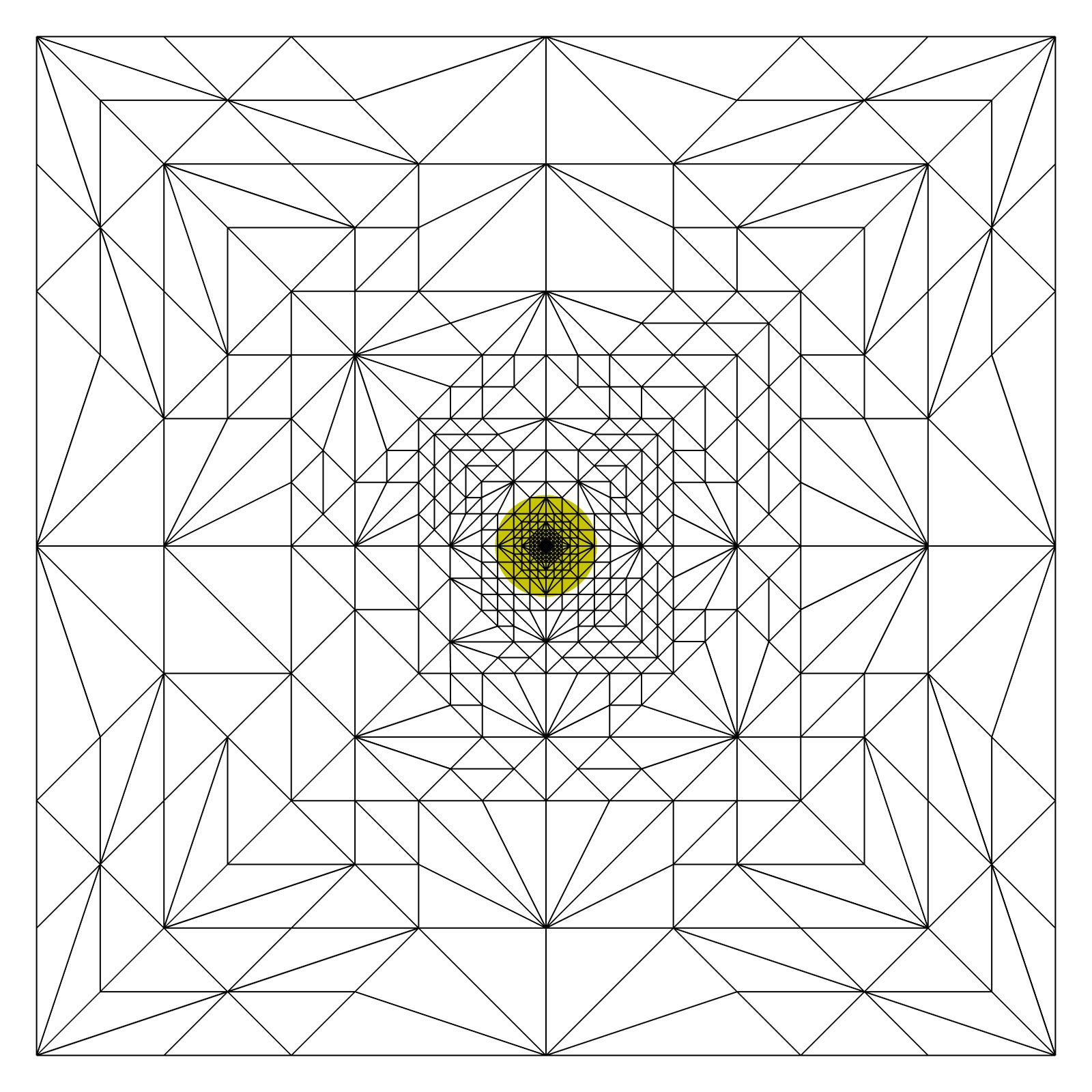}}
    \subfloat[Iter 9]{\includegraphics[width=0.25\textwidth]{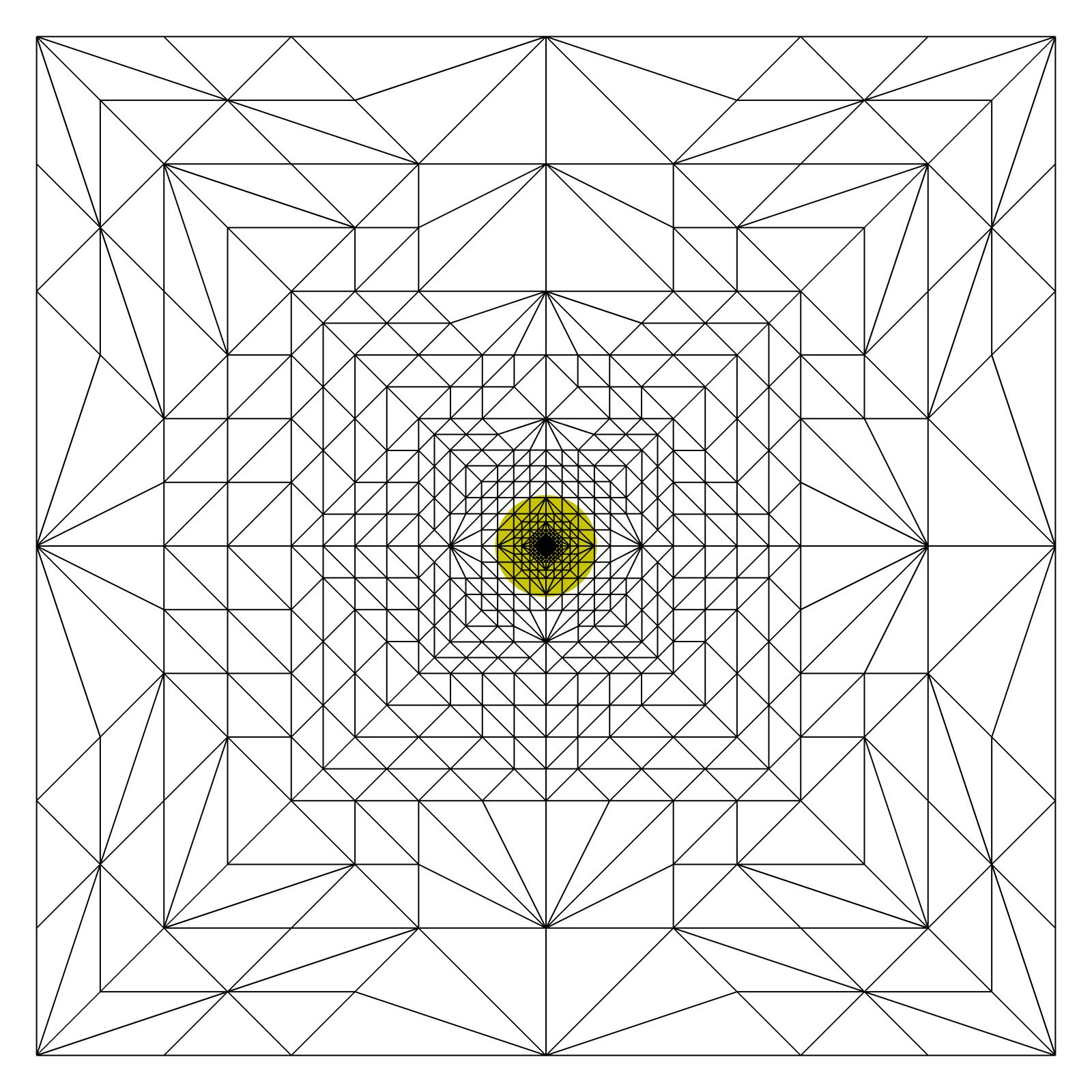}}
    \subfloat[Iter 10]{\includegraphics[width=0.25\textwidth]{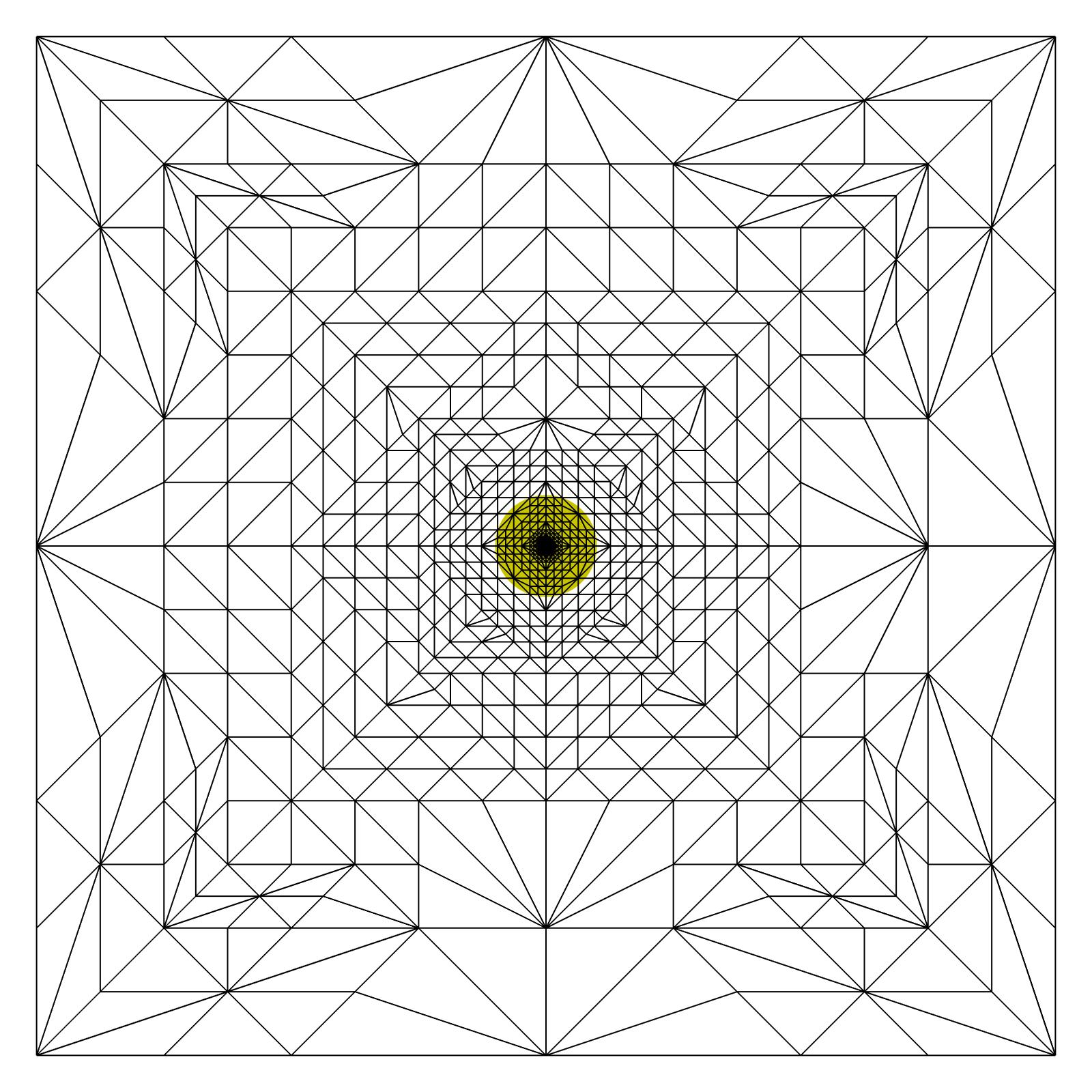}}
    \caption{Adaptive mesh for the Oxygen atom, with the nuclei indicated by disk.}
    \label{figure::KohnShamAdaptiveMesh-O}
\end{figure}

\begin{figure}[H]
    \centering
    \subfloat[Iter 2]{\includegraphics[width=0.25\textwidth]{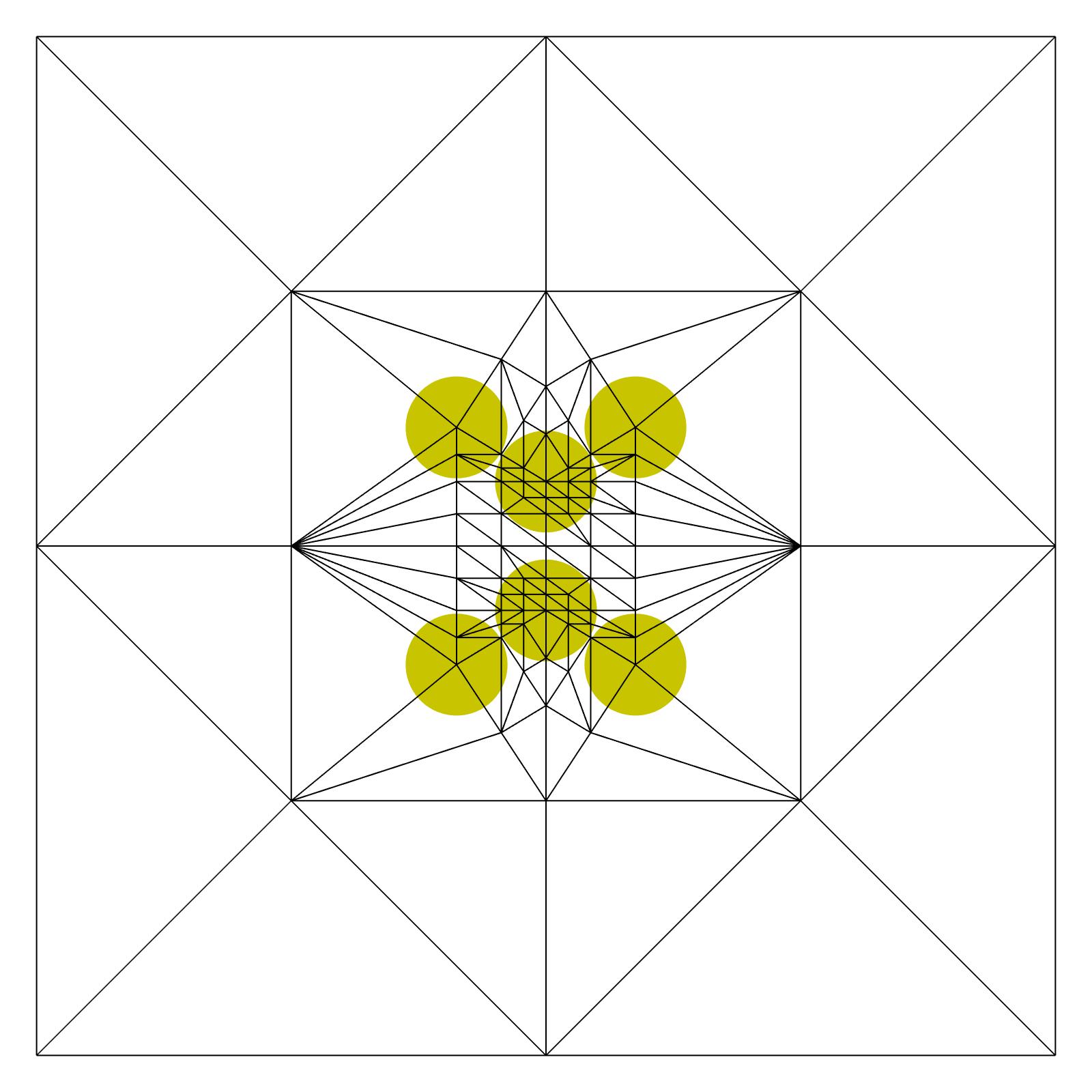}}
    \subfloat[Iter 4]{\includegraphics[width=0.25\textwidth]{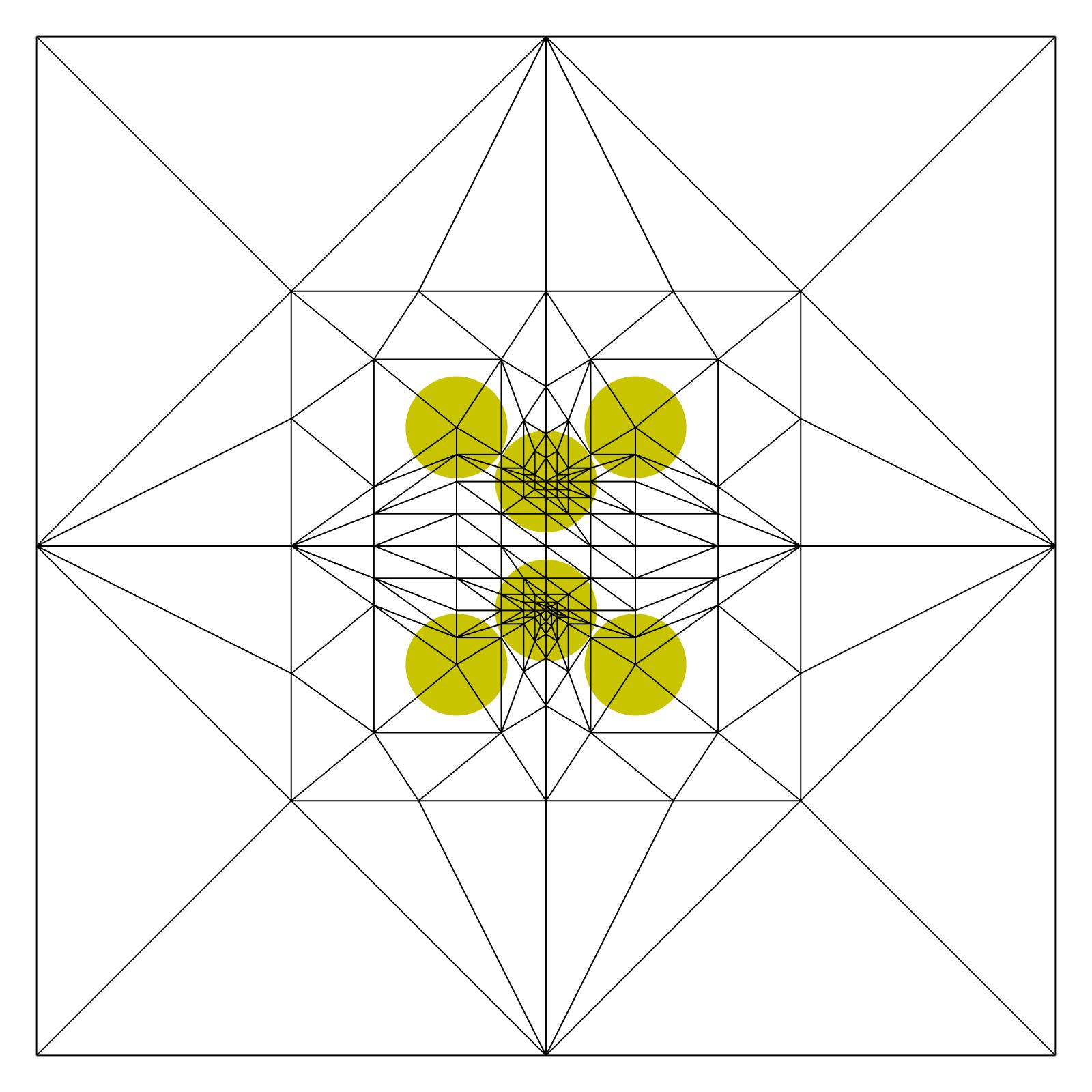}}
    \subfloat[Iter 6]{\includegraphics[width=0.25\textwidth]{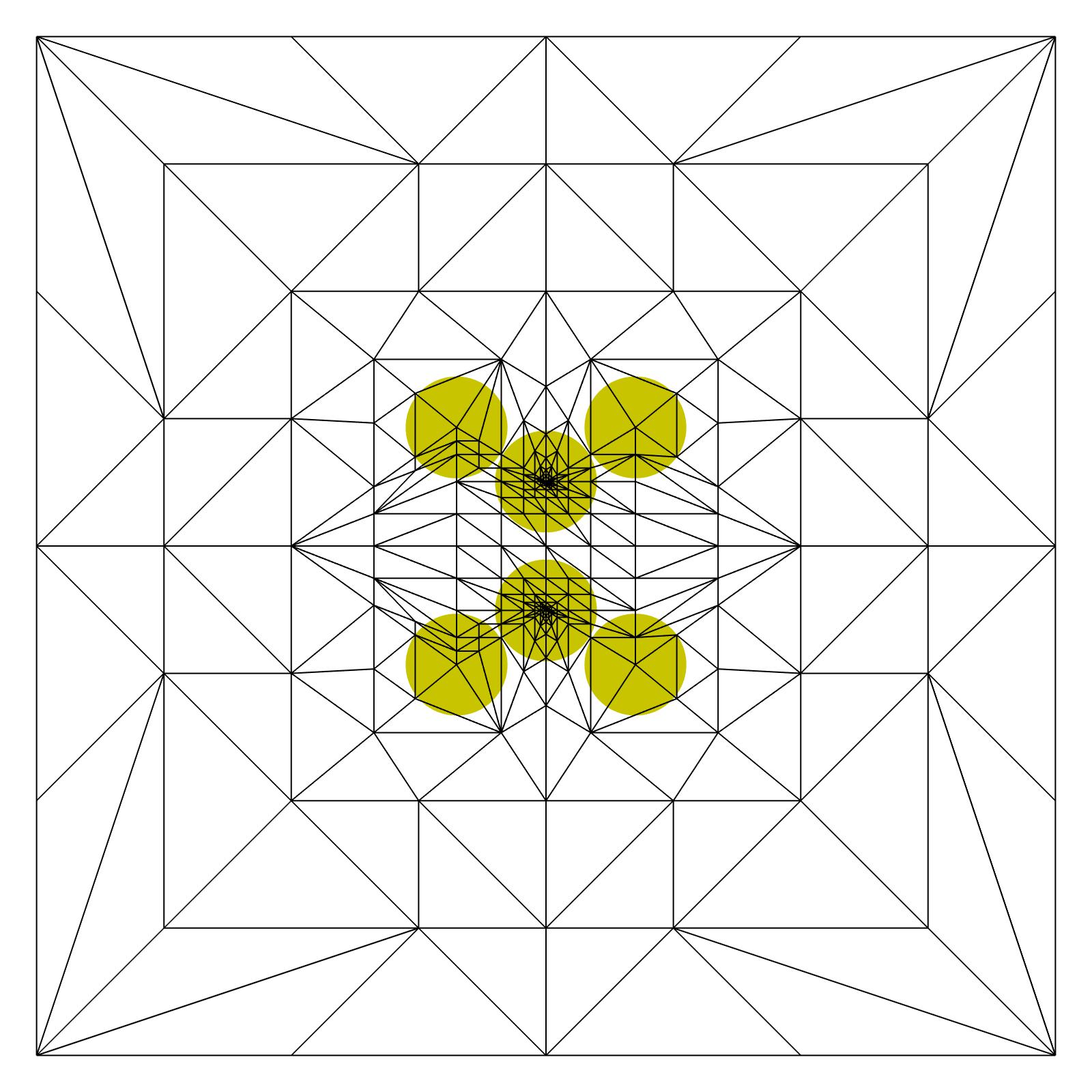}}
    \subfloat[Iter 8]{\includegraphics[width=0.25\textwidth]{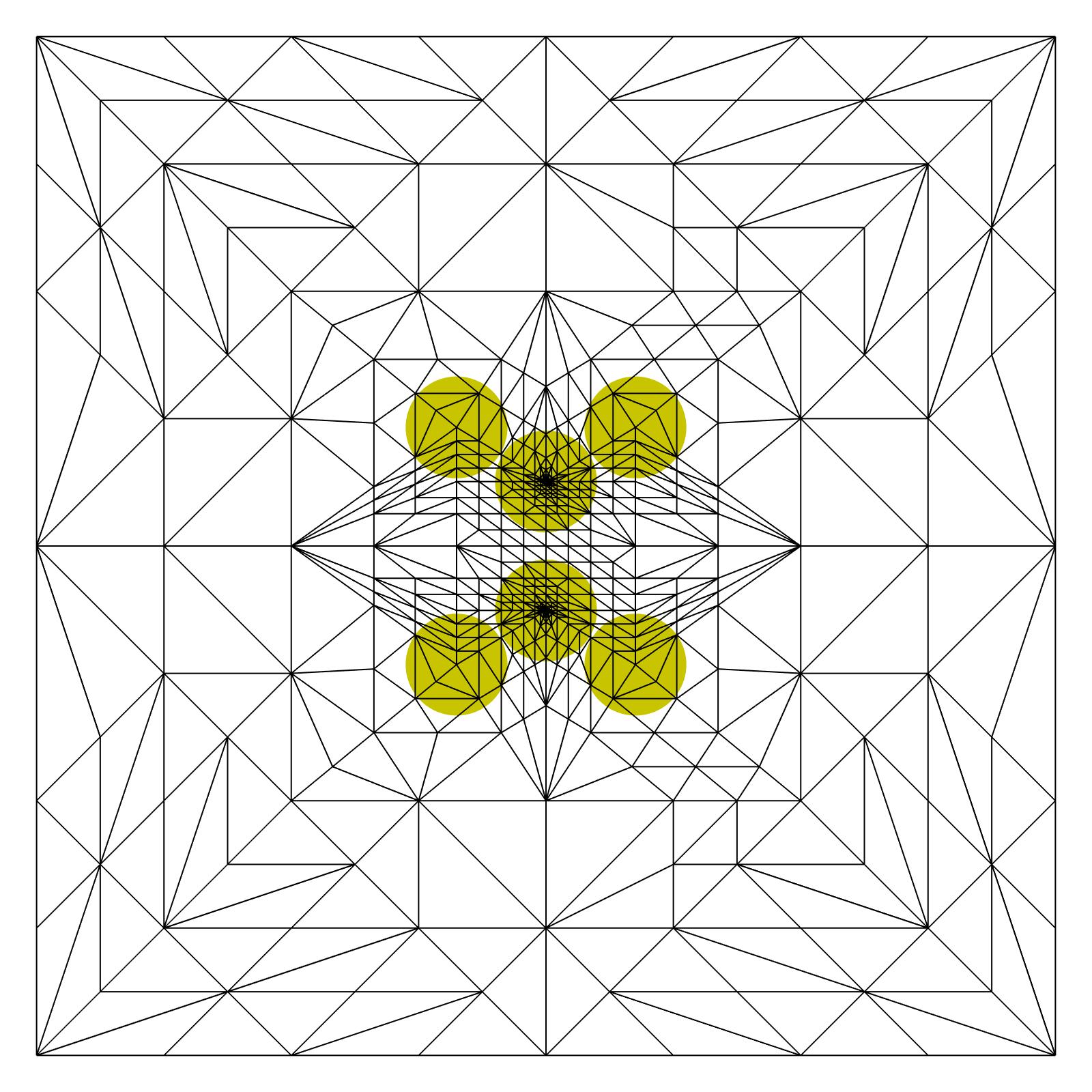}} \
    \subfloat[Iter 10]{\includegraphics[width=0.25\textwidth]{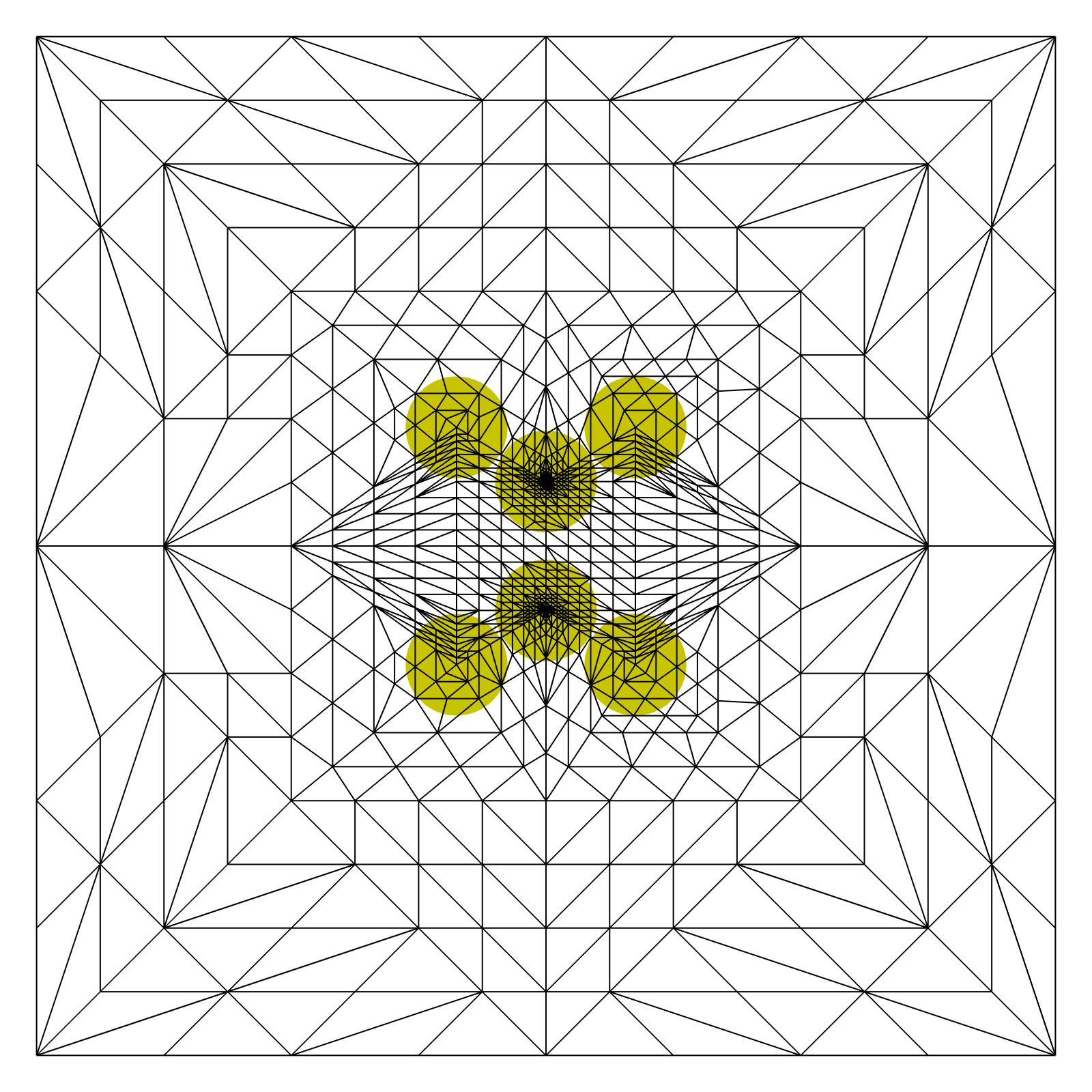}}
    \subfloat[Iter 12]{\includegraphics[width=0.25\textwidth]{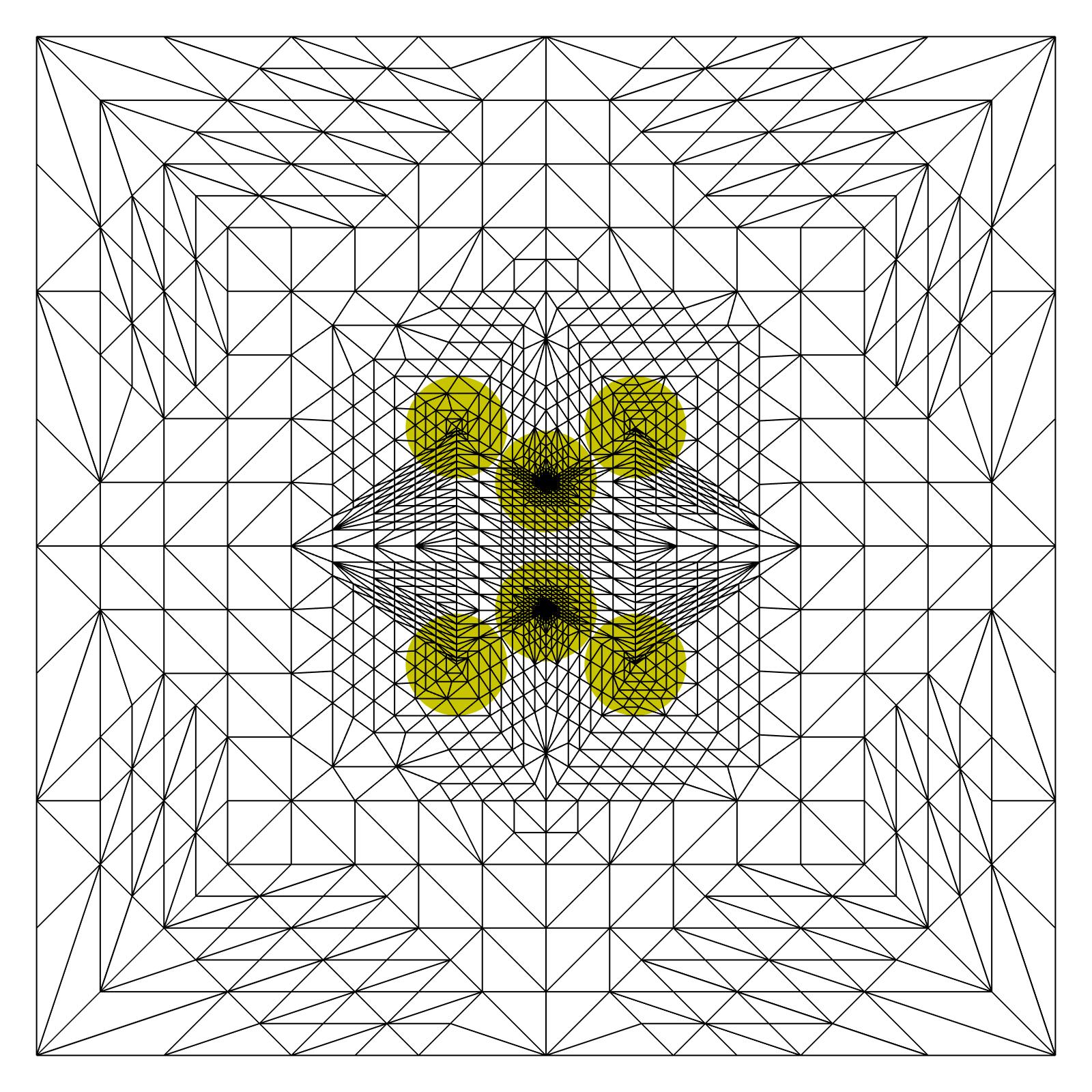}}
    \subfloat[Iter 14]{\includegraphics[width=0.25\textwidth]{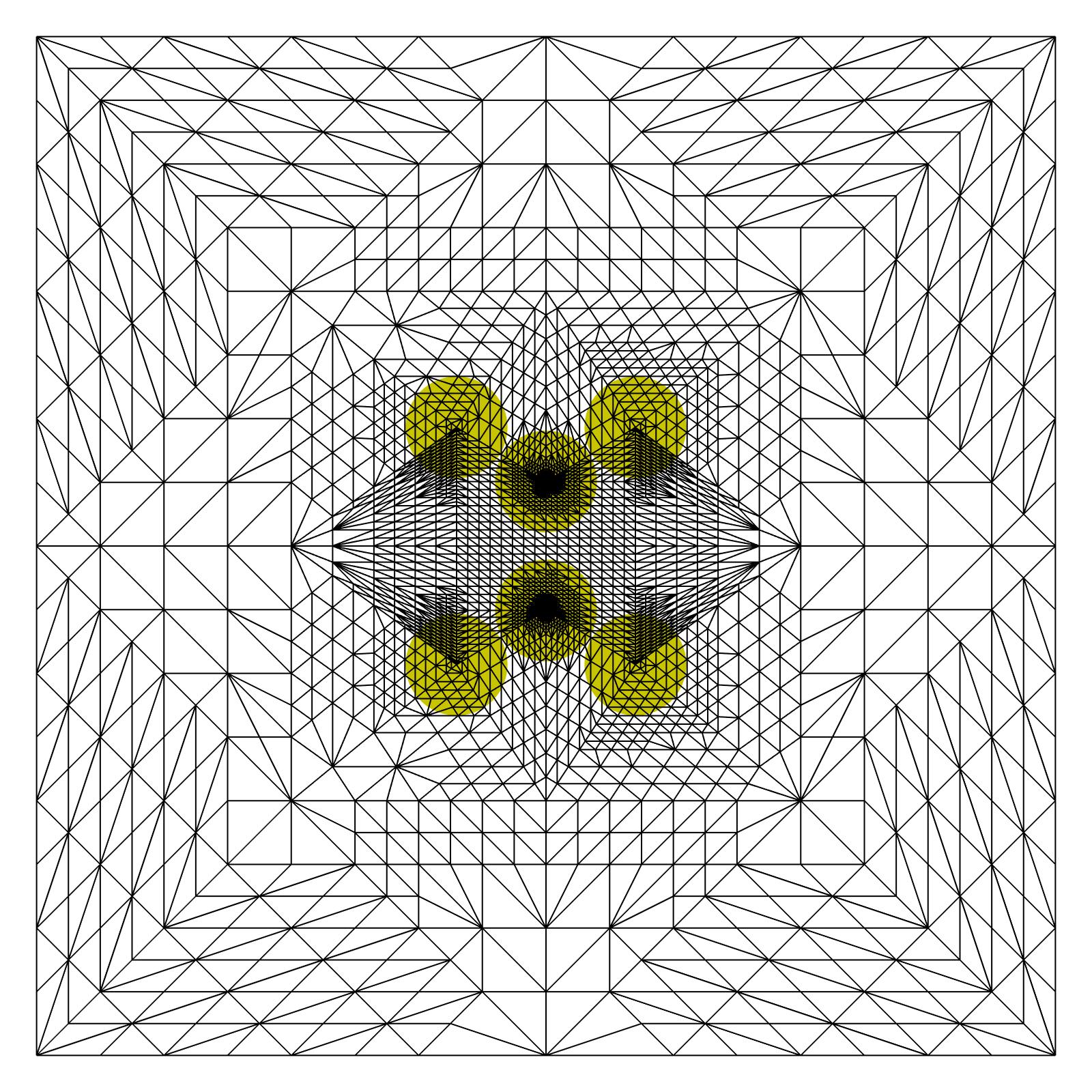}}
    \subfloat[Iter 16]{\includegraphics[width=0.25\textwidth]{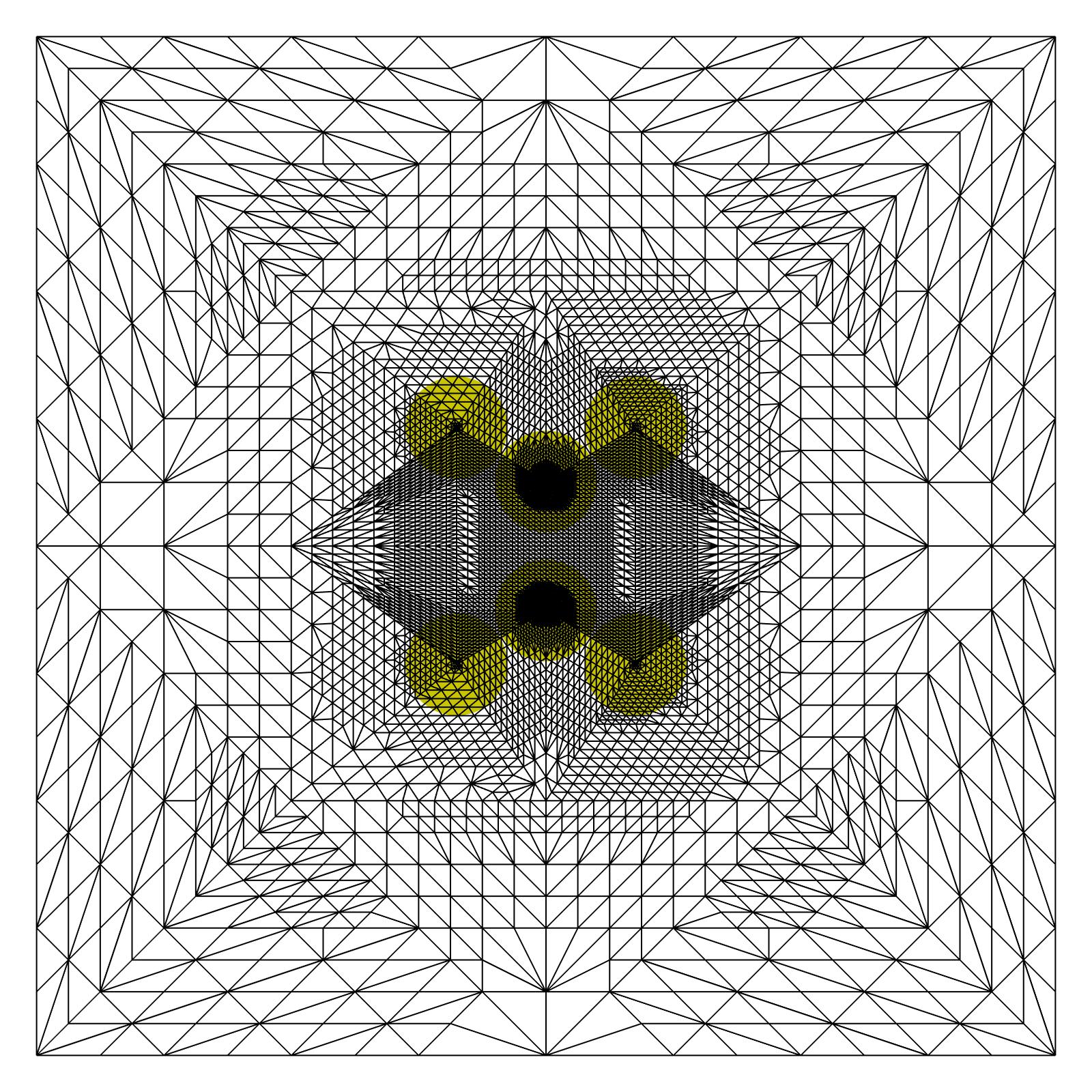}}
    \caption{Adaptive mesh for the Ethene, with the nuclei indicated by disks.}
    \label{figure::KohnShamAdaptiveMesh-C2H4}
\end{figure}

Then, Figures \ref{figure::KohnShamAdaptiveMesh-O} and \ref{figure::KohnShamAdaptiveMesh-C2H4} show the $y$-$z$ plane meshes at different $h$-adaptive iterations for the oxygen atom and ethene, respectively. The mesh is strongly refined near the nuclei due to the rapid variation of the orbital wave functions induced by the nuclear attraction. In addition, for ethene, a denser mesh is observed around the carbon atoms than around the hydrogen atoms, consistent with the stronger nuclear attraction of carbon.

\subsubsection{\zw{Computational Performance}}

Subsequently, we investigate the effect of interpolation in Algorithm \ref{algo::interpolation}. In particular, polynomial basis functions of orders $p=4$ and $p=8$ are tested. We compare the computational time of interpolation with that obtained using a superposition of hydrogen atomic electron densities as the initial guess. In Table \ref{table::TimeCost}, the DoF number is listed in the 2nd column. The iteration counts for SCF and LOBPCG are reported in the 3rd and 4th columns, and the 6th and 7th columns, respectively. The wall time for the artificial initial guess is given in the 5th column, while that for interpolation is presented in the last column. It can be observed that interpolation reduces both the iteration numbers and the wall time, due to the improved quality of the initial guess.


\begin{table}[H]
    \centering
    \caption{The iteration number and computational time (second) for solving Oxygen atom with $64$ cores.}
    \label{table::TimeCost}
    \begin{tabular}{|c|c|c|c|c|c|c|c|c|}
        \hline
        \multirow{3}{*}{$p$} & \multirow{3}{*}{$N_{\text{dof}}$} & \multicolumn{3}{c|}{Superposition of Hydrogen atom} & \multicolumn{4}{c|}{Interpolation} \\
        \cline{3-9}
        & & \multicolumn{2}{c|}{Iteration} & Time ($s$) & \multicolumn{2}{c|}{Iteration} & \multicolumn{2}{c|}{Time ($s$)} \\
        \cline{3-9}
        & & SCF & LOBPCG & LOBPCG & SCF & LOBPCG & LOBPCG & Interpolation \\
        \hline
        $4$ & $12585$  & $10$ & $698$  & $2.51$ & $7$  & $251$ & $0.90$ & $0.50$ \\
        $4$ & $20949$  & $10$ & $807$  & $4.66$ & $7$  & $222$ & $1.30$  & $0.84$ \\
        $4$ & $32009$  & $12$ & $1168$ & $10.11$ & $11$ & $345$ & $3.05$  & $1.22$ \\
        $4$ & $45849$  & $11$ & $1435$ & $18.59$ & $9$  & $334$ & $4.40$  & $1.70$ \\
        $4$ & $64801$  & $13$ & $1787$ & $32.46$ & $8$  & $304$ & $5.57$  & $2.54$ \\
        \hline
        $8$ & $90673$  & $12$ & $2655$ & $92.14$ & $7$  & $254$ & $9.13$  & $2.64$ \\
        $8$ & $166505$ & $13$ & $4129$ & $279.90$ & $8$  & $470$ & $32.60$   & $4.98$ \\
        $8$ & $190853$ & $13$ & $5451$ & $427.45$ & $8$  & $565$ & $44.97$  & $5.29$ \\
        $8$ & $265253$ & $14$ & $5411$ & $581.43$ & $8$  & $409$ & $45.19$  & $7.39$ \\
        $8$ & $530581$ & $13$ & $6731$ & $1494.94$ & $7$  & $515$ & $116.174$  & $15.51$ \\
        \hline
    \end{tabular}
\end{table}


Lastly, we evaluated the parallel scalability of the program with accumulated time costs. The benchmarks were carried out on a server equipped with an AMD EPYC 9534 64-core processor and 1 TB of memory. The tolerances for SCF, PCG and LOBPCG were set to $10^{-6}$, $10^{-10}$ and $10^{-9}$, respectively, and the maximum number of iterations of LOBPCG was set to $1024$.

In Tables \ref{table::TimeCostCore4} and \ref{table::TimeCostCore8}, the computational times of each component is listed separately. In particular, the performance of PCG, LOBPCG, adaptivity, and assembly is similar to that of the Poisson equation. On 64 cores, the adaptivity and assembly modules generally achieve speedups above 15, with one test case attaining an even higher speedup of more than 20.
Furthermore, the speedups of the interpolation step reaches 33.6 and 27.9 for $p=4$ and $p=8$, respectively, suggesting that over 98\% of the computational workload is parallelizable, thus highlighting the strong scalability of the proposed interpolation method.



\begin{table}[H]
    \centering
    \caption{Computational times (seconds) for solving Oxygen atom with different numbers of CPU cores and $p = 4$.}
    \label{table::TimeCostCore4}
    \begin{tabular}{|c|c|c|c|c|c|c|}
        \hline
        Core & PCG & LOBPCG & Interpolation & Adaptivity & Assembly & Total \\
        \hline
        $1$  & $151.22$ & $370.85$ & $2830.65$ & $1396.16$ & $1440.78$ & $6654.90$ \\
        $2$  & $90.86$  & $264.10$ & $2042.02$ & $761.92$  & $802.71$  & $4220.13$ \\
        $4$  & $70.29$  & $225.36$ & $1184.51$ & $414.19$  & $440.05$  & $2479.47$ \\
        $8$  & $69.80$  & $221.39$ & $604.77$  & $240.53$  & $255.81$  & $1482.07$ \\
        $16$ & $54.97$  & $185.72$ & $306.74$  & $144.13$  & $161.23$  & $913.97$ \\
        $32$ & $41.83$  & $163.61$ & $156.64$  & $97.25$   & $115.72$  & $622.14$ \\
        $64$ & $37.64$  & $156.10$ & $84.13$   & $75.91$   & $90.45$   & $486.21$ \\
        \hline
        Speedup & $4.0\times$  & $2.4\times$ & $33.6\times$   & $18.4\times$   & $15.9\times$   & $13.7\times$ \\
        \hline
    \end{tabular}
\end{table}

\begin{table}[H]
    \centering
    \caption{Computational times (seconds) for solving Oxygen atom with different numbers of CPU cores and $p = 8$.}
    \label{table::TimeCostCore8}
    \begin{tabular}{|c|c|c|c|c|c|c|}
        \hline
        Core & PCG & LOBPCG & Interpolation & Adaptivity & Assembly & Total \\
        \hline
        $1$  & $774.56$ & $1564.72$ & $1629.74$ & $919.13$ & $2327.58$ & $7494.35$ \\
        $2$  & $477.09$ & $1057.57$ & $1262.11$ & $498.34$ & $1324.51$ & $4769.61$ \\
        $4$  & $422.82$ & $954.09$  & $833.99$  & $268.96$ & $766.06$  & $3327.08$ \\
        $8$  & $417.01$ & $934.19$  & $439.77$  & $153.02$ & $477.24$  & $2468.50$ \\
        $16$ & $237.61$ & $604.35$  & $224.23$  & $91.72$  & $324.92$  & $1512.08$ \\
        $32$ & $165.44$ & $479.11$  & $115.93$  & $60.11$  & $242.51$  & $1083.75$ \\
        $64$ & $139.60$ & $436.75$  & $58.51$   & $44.17$  & $209.72$  & $906.50$ \\
        \hline
        Speedup & $5.5\times$  & $3.6\times$ & $27.9\times$   & $20.8\times$   & $11.1\times$   & $8.3\times$ \\
        \hline
    \end{tabular}
\end{table}

\section{Conclusion}
\label{chapter::Conclusion}


In this work, we present an $h$-adaptive tetrahedral spectral element method and systematically evaluate its performance using both classical model problems and KSDFT. 
On the one hand, numerical results for the Poisson equation and the Laplacian eigenvalue problem demonstrate the accuracy and efficiency of the proposed framework, where spectral accuracy can be observed, and the adaptive strategy effectively captures the solution behavior.
On the other hand, the proposed framework is successfully applied to the all-electron Kohn-Sham equations, where the nuclear singularities are well resolved.
Scalability results indicate an overall 8-fold speedup. However, whereas the interpolation, adaptivity, and assembly modules closely aligning with Amdahl’s law, the PCG and LOBPCG solvers become the dominant performance bottleneck, limiting the overall acceleration. 

Future work will focus on two directions. 
First, to alleviate the bottleneck in the PCG and LOBPCG solvers, we consider techniques such as domain decomposition to further improve parallel efficiency. 
Second, we plan to extend the present framework to time-dependent problems. High-harmonic generation serves as a promising application, where the stringent accuracy requirements are well aligned with the high-resolution capability of the proposed method. 



\bibliographystyle{unsrt}
\bibliography{reference}

@article{zhan2023novel,
    title     = {A novel tetrahedral spectral element method for Kohn-Sham model},
    author    = {Zhan, Hongfei and Hu, Guanghui},
    journal   = {Journal of Computational Physics},
    volume    = {474},
    pages     = {111831},
    year      = {2023},
    publisher = {Elsevier}
}

@article{parr1995density,
    title     = {Density-Functional Theory of Atoms and Molecules},
    author    = {Parr, Robert G and Weitao, Yang},
    year      = {1995},
    publisher = {Oxford University Press}
}

@article{umemoto2006dissociation,
    title     = {Dissociation of MgSiO3 in the cores of gas giants and terrestrial exoplanets},
    author    = {Umemoto, Koichiro and Wentzcovitch, Renata M and Allen, Philip B},
    journal   = {Science},
    volume    = {311},
    number    = {5763},
    pages     = {983--986},
    year      = {2006},
    publisher = {American Association for the Advancement of Science}
}

@article{motamarri2013higher,
    title     = {Higher-order adaptive finite-element methods for Kohn--Sham density functional theory},
    author    = {Motamarri, Phani and Nowak, Michael R and Leiter, Kenneth and Knap, Jaroslaw and Gavini, Vikram},
    journal   = {Journal of Computational Physics},
    volume    = {253},
    pages     = {308--343},
    year      = {2013},
    publisher = {Elsevier}
}

@article{hohenberg1964inhomogeneous,
    title     = {Inhomogeneous electron gas},
    author    = {Hohenberg, Pierre and Kohn, Walter},
    journal   = {Physical review},
    volume    = {136},
    number    = {3B},
    pages     = {B864},
    year      = {1964},
    publisher = {APS}
}

@article{marques2012libxc,
    title     = {Libxc: A library of exchange and correlation functionals for density functional theory},
    author    = {Marques, Miguel AL and Oliveira, Micael JT and Burnus, Tobias},
    journal   = {Computer physics communications},
    volume    = {183},
    number    = {10},
    pages     = {2272--2281},
    year      = {2012},
    publisher = {Elsevier}
}

@article{umemoto2006namgf3,
    title     = {NaMgF3: A low-pressure analog of MgSiO3},
    author    = {Umemoto, K and Wentzcovitch, RM and Weidner, DJ and Parise, JB},
    journal   = {Geophysical Research Letters},
    volume    = {33},
    number    = {15},
    year      = {2006},
    publisher = {Wiley Online Library}
}

@article{kohn1965self,
    title     = {Self-consistent equations including exchange and correlation effects},
    author    = {Kohn, Walter and Sham, Lu Jeu},
    journal   = {Physical review},
    volume    = {140},
    number    = {4A},
    pages     = {A1133},
    year      = {1965},
    publisher = {APS}
}

@article{ehrenreich1959self,
    title     = {Self-consistent field approach to the many-electron problem},
    author    = {Ehrenreich, H and Cohen, Morrel H},
    journal   = {Physical Review},
    volume    = {115},
    number    = {4},
    pages     = {786},
    year      = {1959},
    publisher = {APS}
}

@article{bao2015real,
    title     = {Real-time adaptive finite element solution of time-dependent Kohn--Sham equation},
    author    = {Bao, Gang and Hu, Guanghui and Liu, Di},
    journal   = {Journal of Computational Physics},
    volume    = {281},
    pages     = {743--758},
    year      = {2015},
    publisher = {Elsevier}
}

@article{dubiner1991spectral,
    title     = {Spectral methods on triangles and other domains},
    author    = {Dubiner, Moshe},
    journal   = {Journal of Scientific Computing},
    volume    = {6},
    pages     = {345--390},
    year      = {1991},
    publisher = {Springer}
}

@book{shen2011spectral,
    title     = {Spectral methods: algorithms, analysis and applications},
    author    = {Shen, Jie and Tang, Tao and Wang, Li-Lian},
    volume    = {41},
    year      = {2011},
    publisher = {Springer Science \& Business Media}
}

@article{pasquetti2004spectral,
    title     = {Spectral element methods on triangles and quadrilaterals: comparisons and applications},
    author    = {Pasquetti, Richard and Rapetti, Francesca},
    journal   = {Journal of Computational Physics},
    volume    = {198},
    number    = {1},
    pages     = {349--362},
    year      = {2004},
    publisher = {Elsevier}
}

@article{guo2006optimal,
    title     = {Optimal spectral-Galerkin methods using generalized Jacobi polynomials},
    author    = {Guo, Ben-Yu and Shen, Jie and Wang, Li-Lian},
    journal   = {Journal of Scientific Computing},
    volume    = {27},
    pages     = {305--322},
    year      = {2006},
    publisher = {Springer}
}

@article{guo2009generalized,
    title     = {Generalized Jacobi polynomials/functions and their applications},
    author    = {Guo, Ben-Yu and Shen, Jie and Wang, Li-Lian},
    journal   = {Applied Numerical Mathematics},
    volume    = {59},
    number    = {5},
    pages     = {1011--1028},
    year      = {2009},
    publisher = {Elsevier}
}

@article{jia2022sparse,
    title     = {Sparse spectral-Galerkin method on an arbitrary tetrahedron using generalized Koornwinder polynomials},
    author    = {Jia, Lueling and Li, Huiyuan and Zhang, Zhimin},
    journal   = {Journal of Scientific Computing},
    volume    = {91},
    number    = {1},
    pages     = {22},
    year      = {2022},
    publisher = {Springer}
}

@article{kotochigova2005atomic,
    title   = {Atomic reference data for electronic structure calculations},
    author  = {Kotochigova, Svetlana and Levine, Zachary H and Shirley, Eric L and Stiles, Mark D and Clark, Charles W},
    journal = {National Institute of Standards and Technology},
    year    = {2005}
}

@article{johnson1999nist,
    title     = {NIST 101. Computational chemistry comparison and benchmark database},
    author    = {Johnson III, Russell D},
    year      = {1999},
    publisher = {Russell D. Johnson III}
}

@article{knyazev2007block,
    title     = {Block locally optimal preconditioned eigenvalue solvers (BLOPEX) in HYPRE and PETSc},
    author    = {Knyazev, Andrew V and Argentati, Merico E and Lashuk, Ilya and Ovtchinnikov, Evgueni E},
    journal   = {SIAM Journal on Scientific Computing},
    volume    = {29},
    number    = {5},
    pages     = {2224--2239},
    year      = {2007},
    publisher = {SIAM}
}

@article{knyazev2001toward,
    title     = {Toward the optimal preconditioned eigensolver: Locally optimal block preconditioned conjugate gradient method},
    author    = {Knyazev, Andrew V},
    journal   = {SIAM journal on scientific computing},
    volume    = {23},
    number    = {2},
    pages     = {517--541},
    year      = {2001},
    publisher = {SIAM}
}

@book{saad2003iterative,
    title     = {Iterative methods for sparse linear systems},
    author    = {Saad, Yousef},
    year      = {2003},
    publisher = {SIAM}
}

@book{golub2013matrix,
    title     = {Matrix computations},
    author    = {Golub, Gene H and Van Loan, Charles F},
    year      = {2013},
    publisher = {JHU press}
}

@article{pulay1980convergence,
    title     = {Convergence acceleration of iterative sequences. The case of SCF iteration},
    author    = {Pulay, P{\'e}ter},
    journal   = {Chemical physics letters},
    volume    = {73},
    number    = {2},
    pages     = {393--398},
    year      = {1980},
    publisher = {Elsevier}
}

@incollection{leitner1995three,
    title     = {Three-dimensional grid adaptation using a mixed-element decomposition method},
    author    = {Leitner, E and Selberherr, S},
    booktitle = {Simulation of Semiconductor Devices and Processes: Vol. 6},
    pages     = {464--467},
    year      = {1995},
    publisher = {Springer}
}

@article{melenk2001residual,
    title     = {On residual-based a posteriori error estimation in hp-FEM},
    author    = {Melenk, Jens Markus and Wohlmuth, Barbara I},
    journal   = {Advances in Computational Mathematics},
    volume    = {15},
    number    = {1},
    pages     = {311--331},
    year      = {2001},
    publisher = {Springer}
}

@misc{elk,
    title        = {{The Elk Code}},
    howpublished = {\url{http://elk.sourceforge.net/}}
}

@article{zhang2015many,
    title     = {How many numerical eigenvalues can we trust?},
    author    = {Zhang, Zhimin},
    journal   = {Journal of Scientific Computing},
    volume    = {65},
    number    = {2},
    pages     = {455--466},
    year      = {2015},
    publisher = {Springer}
}

@book{zienkiewicz1977finite,
    title     = {The finite element method},
    author    = {Zienkiewicz, Olgierd Cecil and Taylor, Robert Leroy and Nithiarasu, Perumal and Zhu, JZ},
    volume    = {3},
    year      = {1977},
    publisher = {Elsevier}
}

@book{hughes2003finite,
    title     = {The finite element method: linear static and dynamic finite element analysis},
    author    = {Hughes, Thomas JR},
    year      = {2003},
    publisher = {Courier Corporation}
}

@book{brenner2008mathematical,
    title     = {The mathematical theory of finite element methods},
    author    = {Brenner, Susanne C and Scott, L Ridgway},
    year      = {2008},
    publisher = {Springer}
}

@article{patera1984spectral,
    title     = {A spectral element method for fluid dynamics: laminar flow in a channel expansion},
    author    = {Patera, Anthony T},
    journal   = {Journal of computational Physics},
    volume    = {54},
    number    = {3},
    pages     = {468--488},
    year      = {1984},
    publisher = {Elsevier}
}

@book{canuto2006spectral,
    title     = {Spectral methods},
    author    = {Canuto, Claudio and Hussaini, M Youssuff and Quarteroni, Alfio and Zang, Thomas A},
    volume    = {285},
    year      = {2006},
    publisher = {Springer}
}

@article{bonito2024adaptive,
    title     = {Adaptive finite element methods},
    author    = {Bonito, Andrea and Canuto, Claudio and Nochetto, Ricardo H and Veeser, Andreas},
    journal   = {Acta Numerica},
    volume    = {33},
    pages     = {163--485},
    year      = {2024},
    publisher = {Cambridge University Press}
}

@article{dorfler1996convergent,
    title     = {A convergent adaptive algorithm for Poisson’s equation},
    author    = {D{\"o}rfler, Willy},
    journal   = {SIAM Journal on Numerical Analysis},
    volume    = {33},
    number    = {3},
    pages     = {1106--1124},
    year      = {1996},
    publisher = {SIAM}
}

@article{2025:arndt.bangerth.ea:deal,
    author  = {Daniel Arndt and Wolfgang Bangerth and Maximilian Bergbauer and Bruno Blais and Marc Fehling and Rene Gassm\"{o}ller
               and Timo Heister and Luca Heltai and Martin Kronbichler and Matthias Maier and Peter Munch and Sam Scheuerman
               and Bruno Turcksin and Siarhei Uzunbajakau and David Wells and Micha{\l} Wichrowski},
    title   = {The deal.II library, Version 9.7},
    journal = {Journal of Numerical Mathematics},
    year    = 2025,
    volume  = 33,
    number  = 4,
    pages   = {403--415},
    doi     = {10.1515/jnma-2025-0115}
}

@article{mfem-2024,
    title     = {High-Performance Finite Elements with {MFEM}},
    author    = {J. Andrej and N. Atallah and J.-P. Bäcker and  J.-S. Camier and
                 D. Copeland and V. Dobrev and Y. Dudouit and T. Duswald and
                 B. Keith and D. Kim and T. Kolev and B. Lazarov and  K. Mittal and
                 W. Pazner and S. Petrides and S. Shiraiwa and M. Stowell and V. Tomov},
    journal   = {The International Journal of High Performance Computing Applications},
    volume    = {38},
    number    = {5},
    pages     = {447-467},
    year      = {2024},
    publisher = {SAGE Publications Sage UK: London, England}
}

@article{hu1999weighted,
    title     = {Weighted essentially non-oscillatory schemes on triangular meshes},
    author    = {Hu, Changqing and Shu, Chi-Wang},
    journal   = {Journal of Computational Physics},
    volume    = {150},
    number    = {1},
    pages     = {97--127},
    year      = {1999},
    publisher = {Elsevier}
}

@article{jiang1996efficient,
    title     = {Efficient implementation of weighted ENO schemes},
    author    = {Jiang, Guang-Shan and Shu, Chi-Wang},
    journal   = {Journal of computational physics},
    volume    = {126},
    number    = {1},
    pages     = {202--228},
    year      = {1996},
    publisher = {Elsevier}
}

@article{li1989multiple,
    title     = {Multiple-harmonic generation in rare gases at high laser intensity},
    author    = {Li, XF and l’Huillier, A and Ferray, M and Lompr{\'e}, LA and Mainfray, G},
    journal   = {Physical Review A},
    volume    = {39},
    number    = {11},
    pages     = {5751},
    year      = {1989},
    publisher = {APS}
}

@article{li2010optimal,
    title   = {Optimal error estimates in Jacobi-weighted Sobolev spaces for polynomial approximations on the triangle},
    author  = {Li, Huiyuan and Shen, Jie},
    journal = {Mathematics of Computation},
    volume  = {79},
    number  = {271},
    pages   = {1621--1646},
    year    = {2010}
}

@article{shan2017triangular,
    title   = {The triangular spectral element method for Stokes eigenvalues},
    author  = {Shan, Weikun and Li, Huiyuan},
    journal = {Mathematics of Computation},
    volume  = {86},
    number  = {308},
    pages   = {2579--2611},
    year    = {2017}
}

@article{mohammad2025advanced,
  title={An advanced algorithm for solving incompressible fluid dynamics: from Navier--Stokes to Poisson equations},
  author={Mohammad, Mutaz and Trounev, Alexander},
  journal={The European Physical Journal Special Topics},
  volume={234},
  number={8},
  pages={2191--2208},
  year={2025},
  publisher={Springer}
}

@article{frkackowiak2011solution,
  title={Solution of the inverse heat conduction problem described by the Poisson equation for a cooled gas-turbine blade},
  author={Fr{\k{a}}ckowiak, A and Wolfersdorf, JV and Cia{\l}kowski, M},
  journal={International Journal of Heat and Mass Transfer},
  volume={54},
  number={5-6},
  pages={1236--1243},
  year={2011},
  publisher={Elsevier}
}

@article{garcia2014survey,
  title={A survey of the parallel performance and accuracy of Poisson solvers for electronic structure calculations},
  author={Garc{\'\i}a-Risue{\~n}o, Pablo and Alberdi-Rodriguez, Joseba and Oliveira, Micael JT and Andrade, Xavier and Pippig, Michael and Muguerza, Javier and Arruabarrena, Agustin and Rubio, Angel},
  journal={Journal of computational chemistry},
  volume={35},
  number={6},
  pages={427--444},
  year={2014},
  publisher={Wiley Online Library}
}

@article{UNKNOWN-CITE,
  title={A general tetrahedral spectral element method and its implementation to Kohn-Sham equation},
  author={Zhan, Hongfei and Jia, Lueling and Li, Huiyuan and Hu, Guanghui},
  journal={Communications in Computational Physics},
  year={2026},
}

@article{ainsworth2003hierarchic,
  title={Hierarchic finite element bases on unstructured tetrahedral meshes},
  author={Ainsworth, Mark and Coyle, Joe},
  journal={International journal for numerical methods in engineering},
  volume={58},
  number={14},
  pages={2103--2130},
  year={2003},
  publisher={Wiley Online Library}
}

@article{wang2007high,
  title={High-order methods for the Euler and Navier--Stokes equations on unstructured grids},
  author={Wang, Zhi Jian},
  journal={Progress in Aerospace Sciences},
  volume={43},
  number={1-3},
  pages={1--41},
  year={2007},
  publisher={Elsevier}
}

@article{visbal2002use,
  title={On the use of higher-order finite-difference schemes on curvilinear and deforming meshes},
  author={Visbal, Miguel R and Gaitonde, Datta V},
  journal={Journal of Computational Physics},
  volume={181},
  number={1},
  pages={155--185},
  year={2002},
  publisher={Elsevier}
}

@book{gottlieb1977numerical,
  title={Numerical analysis of spectral methods: theory and applications},
  author={Gottlieb, David and Orszag, Steven A},
  year={1977},
  publisher={SIAM}
}

@article{gottlieb2001spectral,
  title={Spectral methods for hyperbolic problems},
  author={Gottlieb, David and Hesthaven, Jan S},
  journal={Journal of Computational and Applied Mathematics},
  volume={128},
  number={1-2},
  pages={83--131},
  year={2001},
  publisher={Elsevier}
}

@book{smith1985numerical,
  title={Numerical solution of partial differential equations: finite difference methods},
  author={Smith, Gordon D},
  year={1985},
  publisher={Oxford university press}
}

@book{riviere2008discontinuous,
  title={Discontinuous Galerkin methods for solving elliptic and parabolic equations: theory and implementation},
  author={Rivi{\`e}re, B{\'e}atrice},
  year={2008},
  publisher={SIAM}
}

@book{hesthaven2008nodal,
  title={Nodal discontinuous Galerkin methods: algorithms, analysis, and applications},
  author={Hesthaven, Jan S and Warburton, Tim},
  year={2008},
  publisher={Springer}
}

@article{li2010spectral,
  title={A spectral method on tetrahedra using rational basis functions},
  author={Li, Huiyuan and Wang, L-L},
  journal={International Journal of Numerical Analysis and Modeling},
  volume={7},
  number={2},
  pages={330--355},
  year={2010}
}

@article{cohen2000wavelet,
  title={Wavelet methods in numerical analysis},
  author={Cohen, Albert},
  journal={Handbook of numerical analysis},
  volume={7},
  pages={417--711},
  year={2000},
  publisher={Elsevier}
}

@article{weiss2023spectral,
  title={Spectral element method for 3-D controlled-source electromagnetic forward modelling using unstructured hexahedral meshes},
  author={Weiss, Michael and Kalscheuer, Thomas and Ren, Zhenyong},
  journal={Geophysical Journal International},
  volume={232},
  number={2},
  pages={1427--1454},
  year={2023},
  publisher={Oxford University Press}
}

@article{huang2019spectral,
  title={Spectral-element method with arbitrary hexahedron meshes for time-domain 3D airborne electromagnetic forward modeling},
  author={Huang, Xin and Yin, Changchun and Farquharson, Colin G and Cao, Xiaoyue and Zhang, Bo and Huang, Wei and Cai, Jing},
  journal={Geophysics},
  volume={84},
  number={1},
  pages={E37--E46},
  year={2019},
  publisher={Society of Exploration Geophysicists}
}

@article{hesthaven2000stable,
  title={Stable spectral methods on tetrahedral elements},
  author={Hesthaven, Jan S and Teng, CH},
  journal={SIAM Journal on Scientific Computing},
  volume={21},
  number={6},
  pages={2352--2380},
  year={2000},
  publisher={SIAM}
}

@article{schneider2022large,
  title={A large-scale comparison of tetrahedral and hexahedral elements for solving elliptic PDEs with the finite element method},
  author={Schneider, Teseo and Hu, Yixin and Gao, Xifeng and Dumas, Jeremie and Zorin, Denis and Panozzo, Daniele},
  journal={ACM Transactions on Graphics (TOG)},
  volume={41},
  number={3},
  pages={1--14},
  year={2022},
  publisher={ACM New York, NY}
}

@incollection{henderson1999adaptive,
  title={Adaptive spectral element methods for turbulence and transition},
  author={Henderson, Ronald D},
  booktitle={High-order methods for computational physics},
  pages={225--324},
  year={1999},
  publisher={Springer}
}

@article{mavriplis1994adaptive,
  title={Adaptive mesh strategies for the spectral element method},
  author={Mavriplis, Catherine},
  journal={Computer methods in applied mechanics and engineering},
  volume={116},
  number={1-4},
  pages={77--86},
  year={1994},
  publisher={Elsevier}
}

@inproceedings{hsu1997adaptive,
  title={Adaptive meshes for the spectral element method},
  author={Hsu, Li-Chieh and Mavriplis, Catherine},
  booktitle={9th International Conference on Domain Decomposition Methods},
  pages={374--381},
  year={1997}
}

@article{galvao2008hp,
  title={hp-Adaptive least squares spectral element method for hyperbolic partial differential equations},
  author={Galv{\~a}o, {\'A}rp{\'a}d and Gerritsma, Marc and De Maerschalck, Bart},
  journal={Journal of Computational and Applied Mathematics},
  volume={215},
  number={2},
  pages={409--418},
  year={2008},
  publisher={Elsevier}
}

@article{valenciano2000h,
  title={An h--p adaptive spectral element method for Stokes flow},
  author={Valenciano, Jos{\'e} and Owens, Robert G},
  journal={Applied Numerical Mathematics},
  volume={33},
  number={1-4},
  pages={365--371},
  year={2000},
  publisher={Elsevier}
}

@inproceedings{moxey2017towards,
  title={Towards p-adaptive spectral/hp element methods for modelling industrial flows},
  author={Moxey, D and Cantwell, CD and Mengaldo, G and Serson, D and Ekelschot, D and Peir{\'o}, J and Sherwin, SJ and Kirby, RM},
  booktitle={Spectral and High Order Methods for Partial Differential Equations ICOSAHOM 2016: Selected Papers from the ICOSAHOM conference, June 27-July 1, 2016, Rio de Janeiro, Brazil},
  pages={63--79},
  year={2017},
  organization={Springer}
}

@article{benito2003h,
  title={An h-adaptive method in the generalized finite differences},
  author={Benito, Juan Jos{\'e} and Urena, F and Gavete, L and Alvarez, R},
  journal={Computer methods in applied mechanics and engineering},
  volume={192},
  number={5-6},
  pages={735--759},
  year={2003},
  publisher={Elsevier}
}

@article{hu1998h,
  title={H-adaptive FE analysis of elasto-plastic non-homogeneous soil with large deformation},
  author={Hu, Yuxia and Randolph, MF},
  journal={Computers and Geotechnics},
  volume={23},
  number={1-2},
  pages={61--83},
  year={1998},
  publisher={Elsevier}
}

@article{zienkiewicz1991adaptivity,
  title={Adaptivity and mesh generation},
  author={Zienkiewicz, OC and Zhu, JZ},
  journal={International Journal for Numerical Methods in Engineering},
  volume={32},
  number={4},
  pages={783--810},
  year={1991},
  publisher={Wiley Online Library}
}

@article{salagame1997simple,
  title={A simple p-adaptive refinement procedure for structural shape optimization},
  author={Salagame, Raviprakash R and Belegundu, Ashok D},
  journal={Finite elements in analysis and design},
  volume={24},
  number={3},
  pages={133--155},
  year={1997},
  publisher={Elsevier}
}

@article{demkowicz1989toward,
  title={Toward a universal hp adaptive finite element strategy, Part 1. Constrained approximation and data structure},
  author={Demkowicz, Leszek and Oden, J Tinsely and Rachowicz, Waldemar and Hardy, Oliver},
  journal={Computer Methods in Applied Mechanics and Engineering},
  volume={77},
  number={1-2},
  pages={79--112},
  year={1989},
  publisher={Elsevier}
}

@incollection{mitchell2011survey,
  title={A survey of hp-adaptive strategies for elliptic partial differential equations},
  author={Mitchell, William F and McClain, Marjorie A},
  booktitle={Recent advances in computational and applied mathematics},
  pages={227--258},
  year={2011},
  publisher={Springer}
}

@article{budd2009adaptivity,
  title={Adaptivity with moving grids},
  author={Budd, Chris J and Huang, Weizhang and Russell, Robert D},
  journal={Acta Numerica},
  volume={18},
  pages={111--241},
  year={2009},
  publisher={Cambridge University Press}
}

@article{browne2014fast,
  title={Fast three dimensional r-adaptive mesh redistribution},
  author={Browne, Philip A and Budd, Chris J and Piccolo, Chiara and Cullen, M},
  journal={Journal of Computational Physics},
  volume={275},
  pages={174--196},
  year={2014},
  publisher={Elsevier}
}

@article{CANTWELL2015205,
	title = {Nektar++: An open-source spectral/hp element framework},
	journal = {Computer Physics Communications},
	volume = {192},
	pages = {205-219},
	year = {2015},
	issn = {0010-4655},
	doi = {https://doi.org/10.1016/j.cpc.2015.02.008},
	url = {https://www.sciencedirect.com/science/article/pii/S0010465515000533},
	author = {C.D. Cantwell and D. Moxey and A. Comerford and A. Bolis and G. Rocco and G. Mengaldo and D. {De Grazia} and S. Yakovlev and J.-E. Lombard and D. Ekelschot and B. Jordi and H. Xu and Y. Mohamied and C. Eskilsson and B. Nelson and P. Vos and C. Biotto and R.M. Kirby and S.J. Sherwin},
}

@article{MOXEY2020107110,
	title = {Nektar++: Enhancing the capability and application of high-fidelity spectral/hp element methods},
	journal = {Computer Physics Communications},
	volume = {249},
	pages = {107110},
	year = {2020},
	issn = {0010-4655},
	doi = {https://doi.org/10.1016/j.cpc.2019.107110},
	url = {https://www.sciencedirect.com/science/article/pii/S0010465519304175},
	author = {David Moxey and Chris D. Cantwell and Yan Bao and Andrea Cassinelli and Giacomo Castiglioni and Sehun Chun and Emilia Juda and Ehsan Kazemi and Kilian Lackhove and Julian Marcon and Gianmarco Mengaldo and Douglas Serson and Michael Turner and Hui Xu and Joaquim Peiró and Robert M. Kirby and Spencer J. Sherwin},
}

@article{mahariq2017application,
  title={On the application of the spectral element method in electromagnetic problems involving domain decomposition},
  author={Mahariq, Ibrahim},
  journal={Turkish Journal of Electrical Engineering and Computer Sciences},
  volume={25},
  number={2},
  pages={1059--1069},
  year={2017}
}

@article{kopriva2002computation,
  title={Computation of electromagnetic scattering with a non-conforming discontinuous spectral element method},
  author={Kopriva, David A and Woodruff, Stephen L and Hussaini, M Yousuff},
  journal={International journal for numerical methods in engineering},
  volume={53},
  number={1},
  pages={105--122},
  year={2002},
  publisher={Wiley Online Library}
}

@article{kanungo2019real,
  title={Real time time-dependent density functional theory using higher order finite-element methods},
  author={Kanungo, Bikash and Gavini, Vikram},
  journal={Physical Review B},
  volume={100},
  number={11},
  pages={115148},
  year={2019},
  publisher={APS}
}

@techreport{amd_4th_gen_epyc_2024,
    author       = {{Advanced Micro Devices, Inc. (AMD)}},
    title        = {4th Gen AMD EPYC Processor Architecture},
    institution  = {Advanced Micro Devices, Inc.},
    year         = {2024},
    type         = {White paper},
    url          = {https://www.amd.com/content/dam/amd/en/documents/products/epyc/4th-gen-amd-epyc-processor-architecture-whitepaper.pdf},
    urldate      = {2026-05-11},
    language     = {english}
}

@article{cerveny2019nonconforming,
  title={Nonconforming mesh refinement for high-order finite elements},
  author={Cerveny, Jakub and Dobrev, Veselin and Kolev, Tzanio},
  journal={SIAM Journal on Scientific Computing},
  volume={41},
  number={4},
  pages={C367--C392},
  year={2019},
  publisher={SIAM}
}

@article{zhang2011high,
  title={High-order, multidimensional, and conservative coarse--fine interpolation for adaptive mesh refinement},
  author={Zhang, Qinghai},
  journal={Computer methods in applied mechanics and engineering},
  volume={200},
  number={45-46},
  pages={3159--3168},
  year={2011},
  publisher={Elsevier}
}

@article{ortiz1991adaptive,
  title={Adaptive mesh refinement in strain localization problems},
  author={Ortiz, M and Quigley Iv, JJ},
  journal={Computer Methods in Applied Mechanics and Engineering},
  volume={90},
  number={1-3},
  pages={781--804},
  year={1991},
  publisher={Elsevier}
}

@inproceedings{bussetta2011comparison,
  title={Comparison of Data Transfer Methods between Meshes in the Frame of the Arbitrary Lagrangian Eulerien Formalism},
  author={Bussetta, Philippe and Boman, Romain and Ponthot, Jean-Philippe},
  booktitle={Fifth International Conference on Advanced COmputationalMethods in ENgineering (ACOMEN 2011)},
  year={2011}
}

@article{costin2013numerical,
  title={Numerical study of radial basis function interpolation for data transfer across discontinuous mesh interfaces},
  author={Costin, Will J and Allen, Christian B},
  journal={International Journal for Numerical Methods in Fluids},
  volume={72},
  number={10},
  pages={1076--1095},
  year={2013},
  publisher={Wiley Online Library}
}

@article{garon2020mesh,
  title={Mesh adaptation based on transfinite mean value interpolation},
  author={Garon, Andr{\'e} and Delfour, Michel C},
  journal={Journal of Computational Physics},
  volume={407},
  pages={109248},
  year={2020},
  publisher={Elsevier}
}

\clearpage

\appendix

\section{Red-Green Refinement}
\label{appendix:red-green-refinement}

This appendix provides a detailed description of the red-green refinement strategy used in the adaptive mesh refinement procedure. Starting from a mesh consisting of conforming tetrahedral elements, the refinement strategy follows a queue-driven process that ensures a conforming mesh after each adaptive cycle.

\begin{figure}[!htbp]
    \centering
    \subfloat[Red refinement]{\includegraphics[width=0.3\linewidth]{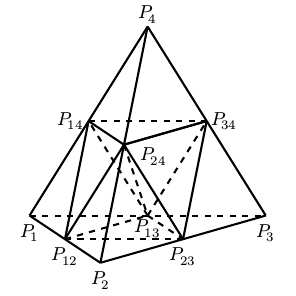} \label{figure::red-green-refinement::a}}
    \subfloat[Green refinement on one edge]{\includegraphics[width=0.3\linewidth]{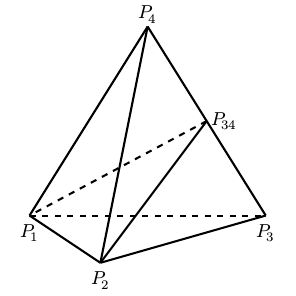} \label{figure::red-green-refinement::b}}
    \subfloat[Green refinement on one face]{\includegraphics[width=0.3\linewidth]{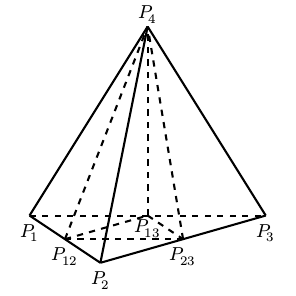} \label{figure::red-green-refinement::c}}
    \caption{The red- and green-refinement templates used in this work: (a) red refinement, (b) and (c) transitional templates for green refinements.}
    \label{figure::red-green-refinement}
\end{figure}

For a tetrahedron marked for refinement or forced to be refined during the closure step, red refinement subdivides it into eight smaller tetrahedra. This is achieved by inserting vertices at the midpoints of all six edges and connecting them appropriately. As shown in Figure \ref{figure::red-green-refinement} \subref{figure::red-green-refinement::a}, points $P_{1}, P_{2}, P_{3}, P_{4}$ are vertices of the tetrahedron, and the point $P_{i j}$'s are the edge midpoints. The red-refinement first construct four corner tetrahedra, which are similar to the original tetrahedron. The remaining region forms an octahedron with the six midpoints as its vertices. Then the octahedron is further subdivided into four tetrahedra by selecting a pair of opposite vertices. One common choice is to connect the pair with the smallest distance, so that the resulting eight tetrahedra are congruent and preserve shape regularity.

Green refinement applies transitional templates to resolve hanging vertices. To guarantee the mesh quality, all green-refined elements are marked as temporary and are subject to coarsening. In this work, two types of templates are applied, depending on the configuration of hanging nodes on the element:
\begin{itemize}
    \item Edge-based template: As shown in Figure \ref{figure::red-green-refinement} \subref{figure::red-green-refinement::b}, when a single edge is subdivided, the tetrahedron is subdivided by connecting the hanging vertex to the opposite vertices;
    \item Face-based template: As shown in Figure \ref{figure::red-green-refinement} \subref{figure::red-green-refinement::c}, when multiple edges on the same face are subdivided, the tetrahedron is refined into four tetrahedra by subdividing the face accordingly.
\end{itemize}

\end{document}